\documentclass{amsart}
\usepackage{amsmath,ifthen,amsthm,amsopn,amssymb,amsfonts}
\usepackage{srcltx,,hyperref}
\usepackage{amscd}
\usepackage[all,cmtip]{xy}

\usepackage{amssymb,latexsym}
\usepackage{color}

\usepackage{amsthm}
\usepackage{eucal}

\newcommand{\showcomments}{yes}
\renewcommand{\showcomments}{no}

\newsavebox{\commentbox}
\newenvironment{com}%
{\ifthenelse{\equal{\showcomments}{yes}}%
{\footnotemark
        \begin{lrbox}{\commentbox}
        \begin{minipage}[t]{1.25in}\raggedright\sffamily\tiny
        \footnotemark[\arabic{footnote}]}
{\begin{lrbox}{\commentbox}}}%
{\ifthenelse{\equal{\showcomments}{yes}}%
{\end{minipage}\end{lrbox}\marginpar{\usebox{\commentbox}}}
{\end{lrbox}}}

\title[Arithmetic manifolds associated to orthogonal groups]{Hodge type theorems for arithmetic manifolds associated to orthogonal groups}

\author{Nicolas Bergeron} 
\thanks{N.B. is a member of the Institut Universitaire de France.}

\address{Institut de Math\'ematiques de Jussieu \\
Unit\'e Mixte de Recherche 7586 du CNRS \\
Universit\'e Pierre et Marie Curie \\
4, place Jussieu 75252 Paris Cedex 05, France \\}
\email{bergeron@math.jussieu.fr}
\urladdr{http://people.math.jussieu.fr/~bergeron}

\author{John Millson}
\thanks{J.M. was partially supported by NSF grant DMS-0907446}
\address{Department of Mathematics\\
University of Maryland\\
College Park, Maryland
20742, USA }
\email{jjm@math.umd.edu}
\urladdr{http://www-users.math.umd.edu/~jjm/}

\author{Colette Moeglin}
\address{Institut de Math\'ematiques de Jussieu \\
Unit\'e Mixte de Recherche 7586 du CNRS \\
4, place Jussieu 75252 Paris Cedex 05, France \\}
\email{moeglin@math.jussieu.fr}
\urladdr{http://www.math.jussieu.fr/~moeglin}

 \DeclareFontFamily{OT1}{rsfs}{}
 
\DeclareFontShape{OT1}{rsfs}{n}{it}{<-> rsfs10}{}
\DeclareMathAlphabet{\mathscr}{OT1}{rsfs}{n}{it}

\newcommand{\C}{\mathbb{C}}

\newcommand{\Z}{\mathbb{Z}}

\newcommand{\p}{\mathfrak{p}}
\newcommand{\q}{\mathfrak{q}}
\newcommand{\g}{\mathfrak{g}}
\newcommand{\Mp}{\mathrm{Mp}}

\DeclareFontFamily{OT1}{rsfs}{}

\DeclareFontShape{OT1}{rsfs}{n}{it}{<-> rsfs10}{}
\DeclareMathAlphabet{\mathscr}{OT1}{rsfs}{n}{it}
\newcommand{\Q}{\mathbb{Q}}

\renewcommand{\k}{\mathfrak{k}}

\newcommand{\R}{\mathbb{R}}
\newcommand{\SO}{\mathrm{SO}}
\newcommand{\Sp}{\mathrm{Sp}}

\newcommand{\OO}{\mathrm{O}}
\newcommand{\GSpin}{\mathrm{GSpin}}

\swapnumbers
\newtheorem{thm}[subsection]{Theorem}  
\newtheorem{lem}[subsection]{Lemma}         
\newtheorem*{lem*}{Lemma}         
\newtheorem{prop}[subsection]{Proposition}
\newtheorem*{prop*}{Proposition}

\newtheorem{conj}[subsection]{Conjecture}

\newtheorem{cor}[subsection]{Corollary}

\theoremstyle{definition}
\newtheorem{defn}[subsection]{Definition}

\newtheorem*{ques}{Question}

\numberwithin{equation}{subsection}

\newcommand{\A}{\mathbb A}

\newcommand{\GL}{\mathrm{GL}}
\newcommand{\SL}{\mathrm{SL}}

\newcommand{\PGL}{\mathrm{PGL}}

\def\adots{\mathinner{\mkern2mu\raise1pt\hbox{.}
\mkern3mu\raise4pt\hbox{.}\mkern1mu\raise7pt\hbox{.}}}

\begin{document}

\begin{com}
{\bf \normalsize COMMENTS\\}
ARE\\
SHOWING!\\
\end{com}

\begin{abstract}  
We show that special cycles generate a large part of the 
cohomology of locally symmetric spaces associated to orthogonal groups.
We prove in particular that classes of totally geodesic submanifolds generate the cohomology groups
of degree $n$ of compact congruence $p$-dimensional hyperbolic manifolds ``of simple type'' as long as $n$ is strictly smaller than $\frac{p}{3}$. 
We also prove that for connected Shimura varieties associated to $\OO (p,2)$ the Hodge conjecture
is true for classes of degree $< \frac{p+1}{3}$. The proof of our general theorem makes use of the recent endoscopic classification of automorphic representations of orthogonal groups by \cite{ArthurBook}.  As such our results are conditional on the hypothesis made in this book, whose proofs have only appear on preprint form so far; see the second paragraph of subsection \ref{org2} below. 
\end{abstract}

\maketitle
\tableofcontents

\section{Introduction} 

\subsection{} Let $D$ be the $p$-dimensional hyperbolic space and let $Y_{\Gamma} = \Gamma \backslash D$ be a compact hyperbolic manifold. 

Thirty-five years ago one of us (J.M.) proved, see  \cite{Millson2}, that if $Y_{\Gamma}$ was a compact hyperbolic manifold of simple arithmetic type (see subsection \ref{simpletype} below) then there was a congruence covering $Y'_{\Gamma} \to Y_{\Gamma}$ such that $Y'_{\Gamma}$ contained a nonseparating embedded totally  geodesic hypersurface $F'$.  Hence the associated homology class $[F']$ was nonzero and the first Betti number of
$Y'_{\Gamma}$ was nonzero.  Somewhat later the first author refined this result to

\begin{prop} \label{prop:millson}
Assume $Y_{\Gamma}$ is arithmetic and contains an (immersed) totally geodesic codimension one
submanifold 
\begin{equation} \label{map}
F \rightarrow Y_{\Gamma}.
\end{equation} 
Then, there exists
a finite index subgroup $\Gamma' \subset \Gamma$ such that the map \eqref{map} lifts to an embedding
$F \hookrightarrow Y_{\Gamma '}$ and the dual class $[F] \in H^1 (Y_{\Gamma '})$ is non zero.
\end{prop}

This result was the first of a series of results on non-vanishing of cohomology
classes in hyperbolic manifolds, see \cite{BC1} for the best known results in that direction. 
In this paper we investigate to which extent classes dual to totally geodesic submanifolds generate the
whole cohomology. We work with {\it congruence} hyperbolic manifolds.

\subsection{}\label{simpletype} First we recall the general definition of congruence hyperbolic manifolds of simple type. 
Let $F$ be a totally real field and $\A$ the ring of adeles of $F$. Let $V$
be a nondegenerate quadratic space over $F$ with $\dim_F V =m$. We assume that $G=\SO(V)$
is compact at all but one infinite place. We denote by $v_0$ the infinite place where $\SO(V)$ is non
compact and assume that $G(F_{v_0}) = \SO (p,1)$. 

Consider the image $\Gamma = \Gamma_K$ in $\SO(p,1)_0$ of the 
intersection $G(F) \cap K$, where $K$ is a compact-open subgroup of $G(\A_f)$
the group of finite ad\`elic points of $G$. According to a classical theorem of Borel and Harish-Chandra, it is a lattice in $\SO(p,1)_0$. It is a cocompact
lattice if and only if $G$ is anisotropic over $F$. If $\Gamma$ is sufficiently deep, i.e. $K$ is a
sufficiently small compact-open subgroup of $G(\A_f)$, then $\Gamma$ is moreover torsion-free.

The special orthogonal group $\SO(p)$ is a maximal compact subgroup of $\SO(p,1)_0$, and the quotient $\SO(p,1)_0 / \SO (p)$ --- the associated symmetric space --- is isometric to
the $n$-dimensional hyperbolic space $D$. 

A {\it compact congruence hyperbolic manifolds of
simple type} is a quotient $Y_K = \Gamma \backslash D$ with $\Gamma = \Gamma_K$ a torsion-free congruence subgroup obtained as above.
$\Gamma \backslash D$
is an $n$-dimensional {\it congruence} hyperbolic manifold. In general, a hyperbolic manifold is  {\it arithmetic} if it shares a common finite  cover with a congruence hyperbolic manifold.

\subsection{} Compact congruence hyperbolic manifolds of
simple type contains many (immersed) totally geodesic codimension one submanifolds to which 
Proposition \ref{prop:millson} applies. In fact: to any 
 totally positive definite sub-quadratic space $U \subset V$
of dimension $n \leq p$ we associate a totally geodesic (immersed) submanifold $c(U,K)$ of codimension $n$ in $Y_K$. Set $H= \SO (U^{\perp})$ so that $H(F_{v_0}) = \SO (p-n , 1)$. 
There is a natural morphism $H \rightarrow G$.
Recall that we can realize $D$ as the set of negative lines in $V_{v_0}$. We then let $D_H$ be the subset of $D$ consisting of those lines which lie in $U^{\perp}_{v_0}$. Let 
$\Gamma_U$ be the image of $H(F) \cap K$ in $\SO (p-n, 1)_0$. The cycle $c(U, K)$ is the image 
of the natural map 
$$\Gamma_U \backslash D_H \rightarrow \Gamma \backslash D.$$
It defines a cohomology class $[c(U,K)] \in H^n (Y_K , \Q)$. 

The following theorem can thus be thought as a converse to Proposition \ref{prop:millson}.

\begin{thm} \label{thm:intro1}
Suppose $n < \frac{p}{3}$. Let $Y_K$ be a
$p$-dimensional compact congruence hyperbolic manifold of simple type. Then $H^n (Y_K, \Q)$
is spanned by the Poincar\'e duals of classes of totally geodesic (immersed) submanifolds of codimension $n$.
\end{thm}

\medskip

\noindent
{\it Remark.} Note that this result is an analogue in constant negative curvature of the results that totally geodesic flat subtori span the homology groups of flat tori (zero curvature) and totally-geodesic subprojective spaces span the homology with $\Z/2$-coefficients for the real projective spaces (constant positive curvature).

\subsection{ } There are results for local coefficients analogous to Theorem \ref{thm:intro1} that are important for the deformation theory of locally homogeneous structures on hyperbolic manifolds.   Let $\rho:\Gamma \to \SO(p,1)$ be the inclusion.  Suppose $G$ is  $\SO(p+1,1)$ resp.
$\GL(p+1, \R)$. In each case we have a natural inclusion $\iota : \SO(p,1) \to G$. For the case $G = \SO(p+1,1)$ the image of $\iota$ is the subgroup leaving the first basis vector of $\R^{p+2}$ fixed,
in the second the inclusion is the ``identity''.
The representation $\widetilde{\rho} = \iota \circ \rho: \Gamma \to G$ is no longer rigid
(note that $\widetilde{\rho}(\Gamma)$ has infinite covolume in $G$).  Though there was some
earlier work this was firmly established in the  early 1980's by Thurston, who discovered
the ``Thurston  bending deformations'', which are nontrivial deformations $\widetilde{\rho}_t, t \in \R$
associated to {\it embedded} totally geodesic hypersurfaces $C_U=c(U,K)$ where $\dim (U) =1$, see \cite{JohnsonMillson} \S 5, for
an algebraic description of these deformations.   It is known that the Zariski tangent space to the real algebraic variety $\mathrm{Hom} (\Gamma, G)$ of  representations at the point  $\widetilde{\rho}$
is the space of one cocycles  $Z^1(\Gamma, \R^{p+1})$ in the first case and $Z^1(\Gamma, \mathcal{H}^2(\R^{p+1}))$ in the second case.  Here  $\mathcal{H}^2(\R^{p+1}))$ denotes the space of harmonic (for the Minkowski metric) degree two polynomials on $\R^{p+1}$. Also trivial deformations correspond to $1$-coboundaries.  Then Theorem 5.1 of \cite{JohnsonMillson} proves that
the tangent vector to the curve $\widetilde{\rho}_t$ at $t=0$ is cohomologous to the Poincar\'e
dual of  the embedded hypersurface $C_U$ equipped with the coefficient $u$ in the first case and
the harmonic projection of $u \otimes u$ in the second  where $u$ is a suitable vector in $U$
determining the parametrization of the curve $\widetilde{\rho}_t$.

We then have 

\begin{thm} \label{thm:intro2}
Suppose $p \geq 4$. Let $\Gamma= \Gamma_K$ be a cocompact congruence lattice of simple type in 
$\SO(p,1)_0$. Then $H^1 (\Gamma , \R^{p+1})$ resp. $H^1 (\Gamma,\mathcal{H}^2(\R^{p+1}))$ are spanned by 
the Poincar\'e duals of (possibly non-embedded) totally-geodesic hypersurfaces with coefficients in $\R^{p+1}$ resp. $\mathcal{H}^2(\R^{p+1})$ .
\end{thm}
 
\medskip
\noindent
{\it Remark.}
First, we remind the reader that the above deformation spaces of representations are locally homeomorphic to deformation spaces of locally homogeneous structures.  In the first case, a hyperbolic structure on a compact manifold $M$ is a fortiori a flat conformal structure and a neighborhood
of $\widetilde{\rho}$ in the first space of representations (into $\mathrm{SO}(p+1,1)$) is homeomorphic to a neighborhood of $M$ in the space of (marked) flat conformal structures.
In the second case (representations into $\PGL(p+1,\R)$) a neighborhood of $\rho$ is homeomorphic to a neighborhood of the hyperbolic manifold in the space of (marked) flat real
projective structures. Thus it is of interest  to describe a neighborhood of $\widetilde{\rho}$ in these two cases. 
 By the above theorem we know the infinitesimal deformations of $\widetilde{\rho}$ are spanned modulo coboundaries by the Poincar\'e duals of totally-geodesic hypersurfaces with coefficients.  The first obstruction (to obtaining a curve of structures or equivalently a curve of representations) can be nonzero, see \cite{JohnsonMillson} who showed that
the first obstruction is obtained by intersecting the representing totally geodesic hypersurfaces with coefficients.   By Theorem \ref{thm:intro6}  we can compute the first obstruction as the restriction of the
wedge of  {\it holomorphic} vector-valued one-forms on $Y^{\C}$. This suggests the higher obstructions will be zero and in fact  the deformation spaces will be cut out from the above first cohomology groups by the vector-valued quadratic equations given by the first obstruction (  see \cite{GoldmanMillson}).

\subsection{} Theorems \ref{thm:intro1} and \ref{thm:intro2} bear a strong ressemblance to the famous
{\it Hodge conjecture} for complex projective manifolds: Let $Y$ be a projective complex manifold. Then  every rational cohomology class of type $(n,n)$ on $Y$ is a linear combination with rational
coefficients of the cohomology classes of complex subvarieties of $Y$. 

Hyperbolic manifolds are not complex (projective) manifolds, so that Theorem \ref{thm:intro1} is not obviously related to the Hodge conjecture. We may nevertheless consider the congruence locally
symmetric varieties associated to orthogonal groups $\OO(p,2)$. These are connected Shimura varieties and as such are projective complex manifolds. As in the case of real hyperbolic manifolds, 
one may associate special algebraic cycles to orthogonal subgroups $\OO(p-n, 2)$ of $\OO(p,2)$. 

The proof of the following theorem now follows the same lines as the proof of Theorem \ref{thm:intro1}.

\begin{thm} \label{thm:intro3}
Let $Y$ be a connected compact Shimura variety associated to the orthogonal group $\OO(p,2)$. Let $n$ be an integer $<  \frac{p+1}{3}$. Then 
every rational cohomology class of type $(n,n)$ on $Y$ is a linear combination with rational
coefficients of the cohomology classes Poincar\'e dual to complex subvarieties of $Y$. 
\end{thm}
Note that the complex dimension of $Y$ is $p$. Hodge theory provides $H^{2n} (Y)$ with a pure Hodge structure of weight $2n$ and we more precisely prove that $H^{n,n} (Y)$ is defined over ${\Bbb Q}$ and  that every rational cohomology class of type $(n,n)$ on $Y$ is a linear combination with rational coefficients of the cup product with some power of the Lefschetz class of cohomology classes
associated to the special algebraic cycles corresponding to orthogonal subgroups.\footnote{We should note that $H^{2n} (Y)$ is of pure $(n,n)$-type when $n < p/4$ but this is no longer the case in general when $n \geq p/4$.}


It is very important to extend 
Theorem \ref{thm:intro3} to the noncompact case 
since, for the case  $\OO(2,19) \cong \OO(19,2)$ (and the correct  isotropic rational quadratic forms of signature $(2,19)$ depending on a parameter $g$, the genus), the resulting noncompact  locally symmetric spaces are now the moduli spaces of quasi-polarized $K3$ surfaces $\mathcal{K}_g$ of genus $g$.
 In Theorem \ref{thm:intro4} of this paper we state a  theorem for general orthogonal groups, $\OO(p,q)$  that as a special case  extends  Theorem \ref{thm:intro3} to   the {\it noncompact} case by proving   that the cuspidal  projections, see subsection \ref{cuspidalprojection},  onto  the cuspidal Hodge summand of type $(n,n)$  of the Poincar\'e-Lefschetz duals of special cycles span the {\it cuspidal} cohomology of Hodge type $(n,n)$ for $n < \frac{p+1}{3}$.  

In \cite{BLMM} (with out new collaborator Zhiyuan Li),  we  extend this last result from the cuspidal cohomology to the reduced $L^2$-cohomology, see \cite{BLMM}, Theorem 0.3.1. and hence, by a result of Zucker, to  the entire cohomology groups $H^{n,n}(M)$. 
In \cite{BLMM}, we apply this extended result  to the special case  $\OO(19,2)$ and  prove that the Noether-Lefschetz divisors (the special cycles 
$c(U,K)$ with $\dim(U) =1$)  generate the Picard variety of the moduli spaces $\mathcal{K}_g, g \geq 2$ thereby    giving  an affirmative solution of the Noether-Lefschetz conjecture formulated by Maulik and Pandharipande in \cite{MP}. 
We note that using work of Weissauer \cite{Weissauer} the result that special cycles span the second homology of  the noncompact Shimura varieties associated to $\OO(n,2)$ for the special case $n=3$,  was proved  earlier  by Hoffman and He \cite{HoffmanHe}. In this  case the Shimura varieties are Siegel modular threefolds and the special algebraic cycles are Humbert surfaces and the main theorem of \cite{HoffmanHe} states that the Humbert surfaces rationally generate the Picard groups of Siegel modular threefolds.  

Finally, we   also point out that in \cite{BMM2} we prove the Hodge conjecture --- as well as its  generalization in the version first  formulated  (incorrectly) by Hodge  --- away from the middle dimensions for Shimura varieties uniformized by complex balls. The main ideas of the proof are the same, although the extension to unitary groups is quite substantial; moreover, the extension  is a  more subtle. In the complex  case sub-Shimura varieties do not provide enough cycles, see \cite{BMM2} for more details.

\subsection{A general theorem} \label{010}
As we explain in Sections \ref{10} and \ref{11}, Theorems \ref{thm:intro1}, \ref{thm:intro2} and \ref{thm:intro3} follow from Theorem \ref{Thm:main9} (see also Theorem \ref{thm:intro4} of this Introduction) which is the main result of our paper. It is concerned with general (i.e. not necessarily compact) arithmetic congruence manifolds $Y=Y_K$ associated to $G=\SO(V)$
as above but such that $G(F_{v_0} )= \SO(p,q)$ with $p+q=m$. It is related to the Hodge conjecture
through a refined Hodge decomposition of the cuspidal cohomology that we now describe.  In the noncompact quotient case it is necesary to first project the Poincar\'e-Lefschetz dual form of the special cycle onto the cusp forms.  We now construct this projection - for another construction not using Franke's Theorem see Section \ref{specialcyclessection}.
\subsection{The cuspidal projection of the class of a special cycle with coefficients} \label{cuspidalprojection}
In case $Y_K$ is not compact the special cycles (with coefficients) are also  not necessarily compact.  However,  they are properly embedded and hence we may  consider them as Borel-Moore cycles or  as  cycles relative to the Borel-Serre boundary of $Y_K$.  The Borel-Serre boundaries   of the special cycles with coefficients have been computed in \cite{FMRes}.  The smooth differential forms on $Y_K$ that are  Poincar\'e-Lefschetz dual to the special cycles are not necessarily $L^2$ but by \cite{Franke} they are cohomologous to forms with automorphic form coefficients  (uniquely up to coboundaries of such forms) which can then be projected onto cusp forms since any automorphic form may be decomposed into an Eisenstein component and a cuspidal component.  We will abuse notation and call the class of the resulting cusp form  the {\it cuspidal projection of the class of the 
special cycle}. Since this cuspidal projection is $L^2$ it has a harmonic projection which can then be used to defined the refined Hodge type(s)  of the cuspidal projection and hence of the original class according to the next subsection. 

\subsection{The refined Hodge decomposition} \label{refinedHodgetype} 
As first suggested by Chern \cite{Chern} the decomposition of exterior powers of the cotangent bundle of $D$ under the action of the holonomy group, i.e. the maximal compact subgroup of $G$, yields a natural notion of {\it refined Hodge decomposition} of the cohomology groups of the associated locally 
symmetric spaces. Recall that 
$$D= \SO_0 (p,q) / (\SO(p) \times \SO(q))$$
and let $\g = \k \oplus \p$ be the corresponding (complexified) Cartan decomposition. 
As a representation of $\SO(p,\C) \times \SO (q, \C)$ the space $\p$ is isomorphic to $V_+ \otimes V_-^*$ where $V_+= \C^p$ (resp. $V_- = \C^q$) is the standard representation of $\SO(p, \C)$ (resp. $\SO(q, \C)$). 
The refined Hodge types therefore correspond to irreducible summands in the decomposition of $\wedge^{\bullet} \p^*$ as a $(\SO(p,\C) \times \SO (q, \C))$-module.
In the case of the group $\mathrm{SU}(n,1)$ (then $D$ is the complex hyperbolic space) it is an exercise to check that one recovers the usual Hodge-Lefschetz decomposition. But in general the decomposition is much finer and in our orthogonal case it is hard to write down the 
full decomposition of $\wedge^{\bullet} \p$ into irreducible modules. Note that, as a $\GL (V_+) \times \GL (V_-)$-module, the decomposition is already quite complicated. We have (see \cite[Equation (19), p. 121]{Fulton}):
\begin{equation} \label{GLdec}
\wedge^R (V_+ \otimes V_-^*) \cong \bigoplus_{\mu \vdash R} S_{\mu} (V_+) \otimes S_{\mu^*} (V_-)^*.
\end{equation}
Here we sum over all partition of $R$ (equivalently Young diagram of size $|\mu|=R$) and $\mu^*$ is the conjugate partition (or transposed Young diagram). 

It nevertheless follows from work of Vogan-Zuckerman \cite{VZ} that very few of the irreducible
submodules of $\wedge^{\bullet} \p^*$ can occur as refined Hodge types of non-trivial coholomogy classes.
The ones which can occur (and do occur non-trivially for some $\Gamma$) are understood in terms of cohomological representations of $G$. We review cohomological representations of $G$ in section \ref{sec:4}. We recall in particular how to associate to each cohomological representation $\pi$ of $G$ a {\it strongly primitive} refined Hodge type. This refined Hodge type
correspond to an irreducible representation of $\SO(p) \times \SO(q)$ which is uniquely determined by 
some special kind of partition $\mu$ as in \eqref{GLdec}, see \cite{MSMF} where these special partitions are called {\it orthogonal}. The first degree where these refined Hodge types can occur is $R=|\mu|$. 
We will use the notation $H^{\mu}$ for the space of the cohomology in degree $R= |\mu|$ corresponding
to this special Hodge type. 

Note that since $\wedge^{\bullet} \mathfrak{p} = \wedge^{\bullet} ( V_+ \otimes V_ -^*)$ the group $\mathrm{SL}(q) = \mathrm{SL}(V_ -)$ acts on $\wedge^{\bullet} \mathfrak{p^*}$. 
In this paper we will be mainly concerned with elements of $(\wedge^{\bullet} \mathfrak{p^*} )^{\SL (q)}$ --- that is elements that are trivial 
on the $V_-$-side. Note that in general $(\wedge^{\bullet} \mathfrak{p^*} )^{\SL (q)}$ is {\it strictly} contained in $(\wedge^{\bullet} \mathfrak{p^*} )^{\SO (q)}$.  Recall that if $q$ is even there exists an invariant element 
$$e_q \in (\wedge^{q} \mathfrak{p^*} )^{\SO (p  ) \times \SL (q) },$$ 
the {\it Euler class/form} (see subsection \ref{eulerform} for the definition). We define $e_q=0$ if $q$ is odd. 
We finally note that if $\mu$ is the partition $q+ \ldots + q$ ($n$ times) then $S_{\mu^*} (V_-)$ is the {\it trivial} representation of $\SL (V_-)$ and the special Hodge type associated to 
$\mu$ occurs in $(\wedge^{nq} \mathfrak{p^*} )^{\SL (q)}$; in that  case we use the notation $\mu = n \times q$.

\subsection{The refined Hodge decomposition of special cycles with coefficients} We also consider general local systems of coefficients.
Let $\lambda$ be a dominant weight for $\SO (p,q)$ with at most $n\leq p/2$ 
nonzero entries and let $E(\lambda)$ be the associated finite dimensional representation of $\SO(p,q)$. 

The reader will verify that the subalgebra $\wedge^{\bullet} (\p^*)^{\mathrm{SL}(q)}$ of $\wedge^{\bullet} (\p^*)$ is invariant under $K_{\infty}= \SO ( p ) \times \SO (q)$. 
Hence, we may form the  associated subbundle 
$$F= D \times_{K_{\infty}} ( \wedge^{\bullet} ( \p^*)^{\mathrm{SL}(q)} \otimes E(\lambda))$$ 
of the bundle 
$$D \times_{K_{\infty}} ( \wedge^{\bullet} ( \p^*) \otimes E(\lambda))$$
of exterior powers of the cotangent bundle of $D$ twisted by $E(\lambda)$. The space of sections of $F$ is invariant under the Laplacian
and hence under harmonic projection, compare \cite[bottom of p. 105]{Chern}. In case $E(\lambda)$ is trivial  the space of sections of $F$ is a subalgebra of the algebra of  differential forms.
 
We denote by $H^{\bullet}_{\mathrm{cusp}} (Y , E(\lambda))^{\rm SC}$ the corresponding subspace (subalgebra if $E(\lambda)$ is trivial)  of $H^{\bullet}_{\mathrm{cusp}} (Y , E(\lambda))$. Note that 
when $q=1$ we have $H^{\bullet}_{\mathrm{cusp}} (Y , E(\lambda))^{\rm SC} = H^{\bullet}_{\mathrm{cusp}} (Y , E(\lambda))$ and when $q=2$ we have 
$$H^{\bullet}_{\mathrm{cusp}} (Y , E(\lambda))^{\rm SC} = \oplus_{n=1}^p H^{n,n}_{\mathrm{cusp}} (Y , E(\lambda)).$$

As above we may associate to $n$-dimensional totally positive sub-quadratic spaces of $V$ special cycles of codimension $nq$ in $Y$ with coefficients 
in the finite dimensional representation $E(\lambda)$. They yield classes 
in $H^{nq}_{\mathrm{cusp}} (Y , E(\lambda))$. In fact we shall show that these classes belong to the subspace $H^{nq}_{\mathrm{cusp}} (Y , E(\lambda))^{\rm SC}$ and it follows from Proposition \ref{PropLit} that
\begin{equation} \label{subringSC}
H^{\bullet}_{\mathrm{cusp}} (Y , E (\lambda))^{\rm SC} = \oplus_{r=0}^{[p/2]} \oplus_{k=0}^{p-2r} e_q^{k} H^{r \times q}_{\mathrm{cusp}} (Y , E(\lambda )).
\end{equation}
(Compare with the usual Hodge-Lefschetz decomposition.) We call $H^{n \times q}_{\mathrm{cusp}} (Y , E (\lambda ))$ 
the {\it primitive part} of $H^{nq}_{\mathrm{cusp}} (Y , E (\lambda))^{\rm SC}$.
We see then that if $q$ is odd the above special classes have {\it pure refined Hodge type} and if $q$ is even each such class is the sum of at most $n+1$
refined Hodge types. In what follows we will consider the primitive part of the special cycles i.e. their projections into the subspace associated to the refined Hodge type $n \times q$: 

The notion of refined Hodge type is a local Riemannian geometric one.  However since we are dealing with locally symmetric spaces there
is an equivalent global definition in terms of  automorphic representations. 
Let $A_{\mathfrak{q}} (\lambda)$ be 
the cohomological Vogan-Zuckerman $(\g , K)$-module $A_{\mathfrak{q}} (\lambda)$ where $\mathfrak{q}$ is a $\theta$-stable parabolic subalgebra of $\g$ whose associated Levi subgroup $L$ is isomorphic to $\mathrm{U}(1)^n \times \SO(p-2n, q)$. 
Let $H^{nq}_{\mathrm{cusp}} (Y , E(\lambda))_{A_{\mathfrak{q}} (\lambda )}$ denote the space of cuspidal harmonic $nq$-forms such that the corresponding automorphic representations of the
adelic orthogonal group have distinguished (corresponding to the noncompact factor) infinite component equal to the unitary representation corresponding to $A_{\mathfrak{q}}(\lambda)$. Then we have 
$$H^{nq}_{\mathrm{cusp}} (Y , E(\lambda))_{A_{\mathfrak{q}} (\lambda )} = H^{n \times q}_{\mathrm{cusp}} (Y , E(\lambda)).$$

Now Theorem \ref{Thm:main9} reads as:

\begin{thm} \label{thm:intro4}
Suppose $p>2n$ and $m-1>3n$. Then the space $H^{n \times q}_{\mathrm{cusp}} (Y , E(\lambda))$ is spanned by the cuspidal projections of classes of special cycles.
\end{thm}

\medskip
\noindent
{\it Remark.} See Subsection \ref{cuspidalprojection} for the definition of the cuspidal projection of the class of a special cycle.  In  case we make  the slightly stronger assumption $p> 2n+1$ it is proved in \cite{FMRes}, see Remark 1.2, that the form of Funke-Millson is {\it square integrable}.  We can then immediately project it into the space of 
cusp forms and arrive at the cuspidal projection without first passing to an automorphic representative. 
\medskip

In degree $R< \min \left( m-3 , pq/4 \right)$ one may deduce from the Vogan-Zuckerman classification of cohomological representations that $H^{R}_{\mathrm{cusp}} (Y , E(\lambda))$ is generated by cup-products of invariant forms with primitive
subspaces $H^{n \times q}_{\mathrm{cusp}} (Y , E(\lambda))$ or $H^{p \times n}_{\mathrm{cusp}} (Y , E(\lambda))$. 
Exchanging the role of $p$ and $q$ we may therefore apply Theorem~\ref{thm:intro4} to prove the following:

\begin{cor} \label{cor:intro}
Let $R$ be an integer $<\min \left( m-3  , pq/4 \right)$. Then the full cohomology group
$H^R_{\mathrm{cusp}} (Y, E )$ is generated by cup-products of classes of totally geodesic cycles and invariant forms.
\end{cor}

Beside proving Theorem \ref{thm:intro4} we also provide strong evidence for the following:
\begin{conj} 
If $p=2n$ or $m-1 \leq 3n$ the space $H^{n\times q}_{\mathrm{cusp}} (Y , E(\lambda))$ is {\rm not} spanned by projections of classes of special cycles.
\end{conj}

In the special case $p=3$, $q=n=1$ we give an example of a cuspidal class 
of degree one in the cohomology of Bianchi hyperbolic manifolds which does not belong to the subspace spanned by classes of special cycles, see Proposition \ref{prop:dim3}.

\subsection{Organisation of the paper} \label{organization}
The proof of Theorem \ref{Thm:main9} is the combination of three main steps. 

The first step is  the work
of Kudla-Millson \cite{KM3} --- as extended by Funke and Millson \cite{FM}. It relates the subspace
of the cohomology of locally symmetric spaces associated to orthogonal groups generated by 
special cycles to certain cohomology classes associated to the ``special theta lift'' using vector-valued  adelic Schwartz functions with a fixed vector-valued component at infinity.  More precisely, the special theta lift restricts the general theta lift to 
Schwartz functions that have  at the distinguished  infinite place where the orthogonal group is noncompact  the fixed 
Schwartz function $\varphi_{nq,[\lambda]}$ taking values in   the vector space $S_{\lambda}(\C^n)^* \otimes \wedge^{nq}(\mathfrak{p})^* \otimes S_{[\lambda]}(V)$, see \S \ref{par:5.3}, at infinity.  The Schwartz functions at the other infinite places are Gaussians (scalar-valued) and at the finite places are scalar-valued and otherwise arbitrary.  The main point is that $\varphi_{nq,[\lambda]}$ is a relative Lie algebra cocycle for the orthogonal group allowing one to interpret the special theta lift cohomologically.    

The second step, accomplished in  Theorem \ref{StepTwo} and depending essentially on  Theorem \ref{Thm:5.8}, is to show that the intersections of the images of the general theta lift
and the special theta lift just described with the subspace of  the cuspidal automorphic forms that have infinite component the Vogan-Zuckerman representation  $A_{\mathfrak{q}}(\lambda)$ coincide (of course the first intersection is potentially larger). In other words the special theta lift accounts for all the cohomology of type
$A_{\mathfrak{q}}(\lambda)$ that may be obtained from theta lifting. This is the analogue of the main result of
the paper of Hoffman and He \cite{HoffmanHe} for the special case of $\SO(3,2)$ and our arguments are very similar to theirs.  Combining the first two steps, we show that, in low degree (small $n$),  {\it all} cuspidal cohomology classes of degree $nq$ and type $A_{\mathfrak{q}}(\lambda)$ 
that can by obtained from the general theta lift coincide with the span of the special chomology classes dual to the special cycles of Kudla-Millson and Funke-Millson. 

The third  step (and it is here that we use Arthur's classification \cite{ArthurBook}) is to show  that in low degree (small $n$)  any cohomology class in $H^{nq}_{\mathrm{cusp}} (Y , E(\lambda))_{A_{\mathfrak{q}} (\lambda )}$ can be obtained as a projection of the class of a theta series. In other words, we prove the low-degree cohomological surjectivity of the general  theta lift (for cuspidal classes of the refined Hodge type $A_{\mathfrak{q}} (\lambda )$). 
In particular in the course of the proof we obtain the following (see Theorem \ref{thm:7.13}):
\begin{thm}\label{thm:intro5}
Assume that $V$ is anisotropic and let $n$ be an integer such that $p>2n$ and  $m-1>3n$. 
Then the global theta correspondence induces an isomorphism between the space of cuspidal holomorphic Siegel 
modular forms, of weight $S_{\lambda} (\C^n)^* \otimes \C_{- \frac{m}{2}}$ at $v_0$ and weight $\C_{-\frac{m}{2}}$ at all the others infinite 
places, on the {\rm connected} Shimura variety associated to the symplectic group $\Sp_{2n} |_F$ and the space 
$$H^{nq} (\mathrm{Sh}^0 (G) ,E(\lambda))_{A_{\mathfrak{q}} (\lambda )} = \lim_{\substack{\rightarrow \\ K}} H^{nq} (Y_K , E(\lambda) )_{A_{\mathfrak{q}} (\lambda )}.$$
\end{thm}
Combining the two steps we find that in low degree the space  
$$\lim_{\substack{\rightarrow \\ K}} H^{nq} (Y_K , E(\lambda) )_{A_{\mathfrak{q}} (\lambda )}$$ is spanned by images of  duals of  special cycles.  From this we deduce (again for small $n$) that  $H^{nq}_{\mathrm{cusp}} (Y , E(\lambda))_{A_{\mathfrak{q}} (\lambda )}$ is spanned by totally geodesic cycles.

The injectivity part of the  previous theorem is not new. It follows from Rallis inner product formula \cite{Rallis}. 
In our case it is due to Li, see \cite[Theorem 1.1]{Li2}. The surjectivity is the subject of
\cite{Moeglin97a,GJS} that we summarize in section \ref{sec:1}. In brief
a cohomology class (or more generally any automorphic form) is in the image of the theta lift if its partial
$L$-function has a pole far on the right. This condition may be thought of as asking that the automorphic
form --- or rather its lift to $\GL(N)$ --- is very non-tempered in all but a finite number of places. 
To apply this result we have to relate this global condition to the local condition that
our automorphic form is of a certain cohomological type at infinity. 

\subsection{} \label{org2} This is where the deep theory of Arthur comes into play. We summarize Arthur's theory in Section \ref{sec:2}. Very briefly: Arthur classifies
automorphic representations of classical groups into global packets. Two automorphic representations belong to the same packet if their partial $L$-functions are the same i.e. if the local components of the two automorphic representations are isomorphic almost everywhere. Moreover in loose terms: Arthur shows that if an automorphic form is very non tempered at one place then it is very non tempered everywhere. To conclude we therefore have to study the cohomological representations
at infinity and show that those we are interested in are very non-tempered, this is the main issue of section \ref{sec:AP}.  Arthur's work on the endoscopic classification of representations of classical groups relates the automorphic representations of the orthogonal groups to the automorphic representations of $\GL (N)$ twisted by some outer automorphism $\theta$. Note however that the relation is made through the stable trace formula for the orthogonal groups (twisted by an outer automorphism in the even case) and the stable trace formula for the twisted (non connected) group $\GL (N) \rtimes \langle \theta \rangle$. 

Thus, as pointed above,  our work uses the hypothesis made in Arthur's book. The twisted trace formula has now been stabilized (see \cite{waldspurgerseoul} and \cite{stabilisationX}). As opposed to the case of unitary groups considered in \cite{BMM2}, there is one more hypothesis to check. Indeed: in Arthur's book there is also an hypothesis about the twisted transfer at the Archimedean places which, in the case of orthogonal groups, is only partially proved by Mezo. This is used by Arthur to find his precise multiplicity formula. We do not use this precise multiplicity formula but we use the fact that a discrete twisted automorphic representation of a twisted $\GL (N)$ is the transfer from a stable discrete representation of a unique endoscopic group. So we still have to know that:  
{\it at a real place, the transfer of the stable distribution which is the sum of the discrete series in one Langlands packet is the twisted trace of an elliptic representation of $\GL(N)$ normalized using a Whittaker functional as in Arthur's book.} 

Mezo  \cite{Mezo2} has proved this result up to a constant which could depend on the Langlands' packet. Arthur's \cite[\S 6.2.2]{ArthurBook} suggests a local-global method to show that this constant is equal to $1$. This is worked out by the third author in \cite{MoeglinNoel} which will eventually be part of a joint work with N. Arancibia and D. Renard.

\subsection{} Part 4 is devoted to applications. Apart from those already mentioned, we deduce from our results and recent results of Cossutta \cite{Cossutta} and Cossutta-Marshall \cite{CossuttaMarshall} an
estimate on the growth of the small degree Betti numbers in congruence covers of hyperbolic manifolds 
of simple type. We also deduce from our results an application to the non-vanishing of certain periods of automorphic forms.

We finally note that the symmetric space $D$ embeds as a totally geodesic and totally real submanifold in the Hermitian symmetric space $D^{\C}$ associated to the unitary group $\mathrm{U} (p,q)$. Also there exists a representation $E(\lambda)^{\C}$ of $\mathrm{U} (p,q)$ whose restriction
of $\mathrm{O}(p,q)$ contains the irreducible representation $E(\lambda)$.  Hence there is a $\mathrm{O}(p,q)$ homomorphism from 
$E(\lambda)^{\C}|\mathrm{O}(p,q)$
to $E(\lambda)$.
As explained in \S \ref{5.4} the form $\varphi_{nq, [\lambda]}$ is best understood as the restriction of a holomorphic  form $\varphi_{nq,0, [\lambda]}$ on $D^{\C}$. Now any $Y=Y_K$ as in \S \ref{010} embeds as 
a totally geodesic and totally real submanifold in a connected Shimura variety $Y^{\C}$ modelled on $D^{\C}$. And the proof of Theorem \ref{thm:intro4} implies:

\begin{thm} \label{thm:intro6}
Suppose $p>2n$ and  $m-1>3n$. Then the space $H^{n \times q}_{\rm cusp} (Y,E(\lambda))$ is spanned by the restriction 
of {\rm holomorphic} forms in $H^{nq,0}_{\rm cusp} (Y^{\C} , E(\lambda)^{\C})$. 
\end{thm}

As holomorphic forms are easier to deal with, we hope that this theorem may help to shed light
on the cohomology of the non-Hermitian manifolds $Y$.

\subsection{More comments} General arithmetic manifolds associated to $\SO(p,q)$ are of two types: The simple type and the non-simple type. In this paper we only deal with the former, i.e. arithmetic manifolds associated to a quadratic space $V$ of signature $(p,q)$ at one infinite place and definite at all other infinite places. Indeed: the manifolds 
constructed that way contain totally geodesic submanifolds associated to subquadratic spaces. 
But when $m=p+q$ is even there are other constructions of arithmetic lattices  in $\SO(p,q)$ commensurable with the group of units of an appropriate skew-hermitian form over a quaternion field (see e.g. \cite[Section 2]{LiMillson}. Note that when $m= 4,8$ there are further constructions that we won't discuss here. 

Arithmetic manifolds of non-simple type do not contain as many  totally geodesic cycles as those of simple type and Theorem \ref{thm:intro4} cannot hold. For example the real hyperbolic manifolds constructed in this way in \cite{LiMillson} do not contain codimension~$1$ totally geodesic submanifolds. We should nevertheless point out that there is a general method to produce nonzero cohomology classes for these manifolds: As first noticed by Raghunathan and Venkataramana, these manifolds can be embedded as totally geodesic and totally real submanifolds in unitary arithmetic manifolds of {\it simple type}, see \cite{RV}. On the latter a general construction due to Kazhdan and extended by Borel-Wallach \cite[Chapter VIII]{BorelWallach} produces nonzero holomorphic cohomology classes as theta series. Their restrictions to the totally real submanifolds we started with can produce nonzero cohomology classes and Theorem \ref{thm:intro6} should still hold in that case. This would indeed follow from our proof modulo the natural extension of \cite[Theorem 1.1 (1)]{GJS} for unitary groups of skew-hermitian form over a quaternion field.

\medskip

{\it We  would like to thank Jeffrey Adams for helpful conversations about this paper.
The second author would like to thank Stephen Kudla and Jens Funke for their collaborations which formed a critical input to this paper.}

\part{Automorphic forms}

\section{Theta liftings for orthogonal groups: some background} \label{sec:1}

\subsection{Notations} 
Let $F$ be a number field and $\A$ the ring of adeles of $F$. Let $V$
be a nondegenerate quadratic space over $F$ with $\dim_F V =m$. 

\subsection{The theta correspondence}
Let $X$ be a symplectic $F$-space with $\dim_F X =2p$. We consider the tensor product $X\otimes V$. 
It is naturally a symplectic $F$-space and we let $\Sp (X \otimes V)$ be the corresponding symplectic 
$F$-group. Then $(\OO(V) , \Sp (X))$ forms a reductive
dual pair in $\Sp (X \otimes V)$, in the sense of Howe \cite{Howe}. We denote by $\Mp (X)$ the metaplectic double cover of $\Sp (X)$ if $m$ is odd or simply $\Sp (X)$ is $m$ is even. 

For a non-trivial additive character $\psi$ of $\A / F$, we may define the oscillator representation 
$\omega_{\psi}$. 
It is an automorphic representation of the metaplectic double cover $\widetilde{\Sp} (X \otimes V)$ 
of $\Sp = \Sp (X \otimes V)$, which is realized in the Schr\"odinger model. The maximal compact subgroup of
$\Sp (X \otimes V)$ is $\mathrm{U} = \mathrm{U}_{pm}$, the unitary group in $pm$ variables. 
We denote by $\widetilde{\mathrm{U}}$ its preimage in $\widetilde{\Sp} (X \otimes V)$. 
The associated space of smooth vectors of $\omega$ is the Bruhat-Schwartz space $\mathcal{S} (V(\A)^p)$. The $(\mathfrak{sp} , \widetilde{\mathrm{U}})$-module associated to $\omega$ is made explicit by the realization of $\omega$ known as the Fock model that we will brefly review in \S \ref{Fock}. Using it,
one sees that the $\widetilde{\mathrm{U}}$-finite vectors in $\omega$ is the subspace 
$\mathbf{S}(V(\A)^p) \subset \mathcal{S} (V(\A)^p)$ obtained by replacing, at each infinite place, the Schwartz space by the {\it polynomial Fock space} $\mathbf{S}(V^p) \subset \mathcal{S} (V^p)$, i.e. the image of holomorphic polynomials on $\C^{pm}$ under the intertwining map from the Fock model of the oscillator
representation to the Schr\"odinger model.

\subsection{} \label{par:1.3} We denote by $\OO_m (\A)$, $\Mp_{2p} (\A)$ and $\widetilde{\Sp}_{2pm} (\A)$
the adelic points of respectively $\OO(V)$, $\Mp (X)$ and $\widetilde{\Sp} (X \otimes V)$.
The global metaplectic group $\widetilde{\Sp}_{2pm} (\A )$ acts in $\mathcal{S} (V(\A)^p)$ via $\omega$ and preserves the dense subspace $\mathbf{S}(V(\A)^p)$. 
For each $\phi \in \mathcal{S} (V(\A)^p)$ we form the theta function 
\begin{equation}
\theta_{\psi , \phi} (x) = \sum_{\xi \in V(F)^p} \omega_{\psi} (x ) (\phi) (\xi) 
\end{equation}
on $\widetilde{\Sp}_{2pm} (\A)$. There is a natural homomorphism
$$\OO_m (\A) \times \Mp_{2p} (\A) \rightarrow \widetilde{\Sp}_{2pm} (\A)$$
which is described with great details in \cite{JiangSoudry}. We pull the oscillator representation 
$\omega_{\psi}$ back to $\OO_m (\A) \times \Mp_{2p} (\A)$.  Then 
$(g,g') \mapsto \theta_{\psi, \phi} (g',g)$ is a smooth, slowly increasing function on 
$\OO (V) \backslash \OO_m (\A) \times 
\Mp (X) \backslash \Mp_{2p} ( \A)$; see \cite{Weil,Howe}.

\subsection{The global theta lifting} \label{par:1.4}
We denote by $\mathcal{A}^c (\Mp (X))$ the set of irreducible cuspidal automorphic representations 
of $\Mp_{2p} (\A)$, which occur as irreducible subspaces in the space of cuspidal automorphic
functions in $L^2 ( \Mp (X) \backslash \Mp_{2p} (\A))$. For a $\pi ' \in \mathcal{A}^c (\Mp (X))$, the integral
\begin{eqnarray} \label{theta}
\theta_{\psi, \phi}^f (g) = \int_{\Mp (X)  \backslash  \Mp_{2p} (\A)} \theta_{\psi, \phi} (g,g') f(g') dg' , 
\end{eqnarray}
with $f \in H_{\pi '}$ (the space of $\pi '$), defines an automorphic function on 
$\OO_m (\A)$~:  the integral \eqref{theta} 
is well defined, and determines a slowly increasing function on $\OO (V) \backslash \OO_m (\A)$.
We denote by $\Theta_{\psi, X}^V (\pi ')$ the space of the automorphic representation generated by all 
$\theta_{\psi, \phi}^f (g)$ as $\phi$ and $f$ vary, and call $\Theta_{\psi,X}^V (\pi ')$ the $\psi$-theta lifting 
of $\pi '$ to $\OO_m (\A)$. Note that, since $\mathbf{S}(V(\A)^p)$ is dense in $\mathcal{S} (V(\A)^p)$ we may as well let $\phi$ vary in the subspace $\mathbf{S}(V(\A)^p)$.

We can similarly define $\mathcal{A}^c (\OO (V))$ and 
$\Theta_{\psi, V}^X$ the $\psi$-theta correspondence from $\OO (V)$ to $\Mp (X)$. 

\subsection{} It follows from \cite{Moeglin97a} and from \cite[Theorem 1.3]{JiangSoudry} that if 
$\Theta_{\psi,X}^V (\pi ')$ contains non-zero cuspidal automorphic functions on $\OO_m (\A)$ then
the representation of $\OO_m (\A)$ in $\Theta_{\psi,X}^V (\pi ')$ is irreducible (and cuspidal). 
We also denote by $\Theta_{\psi,X}^V (\pi ')$ the corresponding element of $\mathcal{A}^c (\OO (V))$
In that case it moreover follows from \cite{Moeglin97a} and \cite[Theorem 1.1]{JiangSoudry} that 
$$\Theta_{\psi^{-1} , V}^X (\Theta_{\psi , X}^V (\pi ') )= \pi' .$$

We say that a representation $\pi \in \mathcal{A}^c (\OO (V))$ {\it is in the image of the cuspidal $\psi$-theta correspondence from a smaller group} if there exists a symplectic space $X$ with $\dim X \leq m$ 
and a representation $\pi ' \in \mathcal{A}^c (\Mp (X))$ such that
$$\pi =\Theta_{\psi,X}^V (\pi ').$$

\subsection{} \label{1.6}The main technical point of this paper is to prove that if $\pi  \in \mathcal{A}^c (\OO (V))$
is such that its {\it local} component at infinity is ``sufficiently non-tempered'' (this has to be made precise) then the {\it global} representation $\pi$ is in the image of the cuspidal $\psi$-theta correspondence from a smaller group. 

As usual we encode local components of $\pi$ into an $L$-function. In fact we only consider its 
{\it partial $L$-function} $L^S (s, \pi)= \prod_{v \not\in S} L(s , \pi_v)$ where $S$ is a sufficiently big finite set of places such that 
$\pi_v$ is unramified for each $v \not\in S$. For such a $v$ we define the local factor $L (s , \pi_v)$
by considering the Langlands parameter of $\pi_v$. 

\medskip
\noindent
{\it Remark.} We will loosely identifies the partial $L$-function of $\pi$ and that of its restriction to $\SO(V)$. However we should note that the restriction of $\pi_v$ to the special orthogonal group may be reducible: If $v \not\in S$ we may associate to the Langlands parameter of $\pi_v$ representations from the principal series of the special orthogonal group. Each of these has a unique unramified subquotient and the restriction of $\pi_v$ to the special orthogonal group is then the sum of the non-isomorphic subquotients.\footnote{There are at most two such non-isomorphic subquotients.} Anyway: the local $L$-factor is the same for each summand of the restriction as it only depends on the Langlands parameter.

\medskip

We may generalize these definitions to form the partial $L$-functions $L^S (s , \pi \times \eta)$ for
any automorphic character $\eta$.

Now the following proposition is a first important step toward the proof that a ``sufficiently non-tempered'' 
automorphic representation is in the image of the cuspidal $\psi$-theta correspondence from a smaller group. 
It is symmetric to \cite[Theorem 7.2.5]{KudlaRallis} and is the subject of \cite{Moeglin97a} and \cite[Theorem 1.1 (1)]{GJS}; it is revisited and generalized in \cite{GanTakeda}.

\begin{prop} \label{Prop:GJS}
Let $\pi \in  \mathcal{A}^c (\OO (V))$ and let $\eta$ be a quadratic character of $F^* \backslash \A^*$. 
Let $a$ be a nonnegative integer with $a+1 \equiv m \ \mathrm{mod} \ 2$. We assume that the partial 
$L$-function $L^S (s , \pi \times \eta)$ is holomorphic in the half-plane $\mathrm{Re} (s) > \frac12 (a+1)$
and has a pole in $s= \frac12 (a+1)$. Let $p= \frac12 (m-a-1)$ and $X$ be a symplectic $F$-space with $\dim X = 2p$.

Then there exists an automorphic sign character $\epsilon$ of $\OO_m (\A)$ such that the $\psi^{-1}$-theta lifting 
of $(\pi \otimes \eta) \otimes \epsilon$ to $\Mp_{2p} (\A)$ does not vanish. 
\end{prop}

The big second step to achieve our first goal will rely on Arthur's theory.

\section{Arthur's theory} \label{sec:2}

\subsection{Notations} \label{notations2}
Let $F$ be a number field, $\A$ its ring of adeles and $\Gamma_F  = {\rm Gal} (\overline{{\mathbb Q}} / F)$. Let $V$ be a nondegenerate quadratic space over $F$ with $\dim_F V =m$. 
We let $G$ be the special orthogonal group $\SO (V)$ over $F$. We set $\ell=[m/2]$ and $N=2\ell$. 

The group $G$ is an inner form of a quasi-split form $G^*$. As for now Arthur's work only deals with quasi-split groups. We first describe the group $G^*$ according to the parity of $m$ and briefly recall
the results of Arthur we shall need. We recall from the introduction that Arthur's work rely on extensions to the twisted case of two results which have only been proved so far in the case of connected groups: The first is the
stabilization of the twisted trace formula for the two groups $\GL(N)$ and $\SO(2n)$, see \cite[Hypothesis 3.2.1]{ArthurBook}. The second is Shelstad's strong spectral transfer of tempered archimedean characters. Taking these for granted we will explain how to deal with non-quasi-split groups in the next section.

\subsection{} We first assume that $m=N+1$ is odd. Then the  special orthogonal group $G$  is an inner form of the {\it split} form $G^* = \SO (m)$ over $F$ associated to the symmetric bilinear form
whose matrix is
$$J= \left(
\begin{array}{ccc}
0 & & 1 \\
 & \adots & \\
1 &  & 0 
\end{array} \right).$$ 
The (complex) dual group of $G^*$ is $G^{\vee} = \Sp (N, \C)$ and ${}^L \- G = G^{\vee} \times \Gamma_F$.

\subsection{} We now assume that $m=N$ is even. We let $\SO(N)$ be the {\it split} orthogonal 
group over $F$ associated to the symmetric bilinear form
whose matrix is $J$.
The {\it quasi-split} forms of $\SO (N)$ are parametrized by morphisms $\Gamma_F \rightarrow \Z/2\Z$,
which by class field theory correspond to characters $\eta$ on $F^* \backslash {\mathbb A}^*$ such that
$\eta^2 =1$ --- {\it quadratic Artin characters}. We denote by $\SO (N, \eta)$ the outer twist of the split
group $\SO (N)$ determined by $\eta$: the twisting is induced by the action of $\Gamma_F$ on the
Dynkin diagram via the character $\eta$. 

When $m=N$ is even, there exists a quadratic Artin character $\eta$ such that $G$ is an inner form
of the {\it quasi-split} group $G^* = \SO(N,\eta)$. The (complex) dual group of $G^*$ is then 
$G^{\vee} = \SO (N, \C)$ and ${}^L \- G = G^{\vee} \rtimes \Gamma_F$, where $\Gamma_F$ acts on
$G^{\vee}$ by an order $2$ automorphism --- trivial on the kernel of $\eta$ --- and fixes a splitting, see 
\cite[p. 79]{BC} for an explicit description.

\medskip
\noindent
{\it Remark.} Let $v$ be an infinite real place of $F$ such that $G(F_v) \cong \SO (p,q)$ with $m=p+q$
even so that $m=2\ell$. Then $\eta_v$ is trivial if and only if $(p-q)/2$ is even. We are lead to the following dichotomy for
real orthogonal groups: if $(p-q)/2$ is odd, $\SO (p,q)$ is an inner form of $\SO (\ell -1 , \ell +1)$ and 
if $(p-q)/2$ is even, $\SO (p,q)$ is an inner form of $\SO (\ell , \ell)$ (split over $\R$).

\subsection{Global Arthur parameters}
In order to extends the classification \cite{MW}
of the discrete automorphic spectrum of $\GL(N)$ to the classical groups, Arthur represents the
discrete automorphic spectrum of $\GL (N)$ by a set of formal tensor products
$$\Psi = \mu \boxtimes R,$$
where $\mu$ is an irreducible, unitary, cuspidal automorphic representation of $\GL (d)$ and $R$
is an irreducible representation of $\SL_2 (\C)$ of dimension $n$, for positive integers $d$ and $n$ such that $N=dn$. For any such $\Psi$, we form the induced representation
$$\mathrm{ind} (\mu | \cdot |^{\frac12 (n-1)} , \mu | \cdot |^{\frac12 (n-3)} , \ldots , \mu | \cdot |^{\frac12 (1-n)} )$$
(normalized induction from the standard parabolic subgroup of type $(d, \ldots , d)$). We then write 
$\Pi_{\Psi}$ for the unique irreducible quotient of this representation. 

We may more generally associate an automorphic representation  $\Pi_{\Psi}$ of $\GL(N)$ to a 
formal sum of formal tensor products:
\begin{equation} \label{psi}
\Psi = \mu_1 \boxtimes R_1 \boxplus \ldots \boxplus \mu_r \boxtimes R_r
\end{equation}
where each $\mu_j$ is an irreducible, unitary, cuspidal automorphic representation of $\GL (d_i) / F$, 
$R_j$ is an irreducible representation of $\SL_2 (\C)$ of dimension $n_j$ and $N=n_1d_1 + 
\ldots + n_r d_r$. 

Now consider the outer automorphism:
$$\theta : x \mapsto J {}^t x^{-1} J = J {}^t x^{-1} J^{-1}  \ \ \ (x \in \GL (N)).$$
There is an action $\Pi_{\Psi} \mapsto \Pi_{\Psi}^{\theta}$ on the set of representations $\Pi_{\Psi}$.
If $\Psi$ is as in \eqref{psi}, set
$$\Psi^{\theta} =  \mu_1^{\theta} \boxtimes R_1 \boxplus \ldots \boxplus \mu_r^{\theta} \boxtimes R_r.$$
Then $\Pi_{\Psi}^{\theta} = \Pi_{\Psi^{\theta}}$. 

Arthur's main result \cite[Theorem 1.5.2]{ArthurBook} (see also \cite[Theorem 30.2]{Arthur}) then parametrizes the discrete automorphic spectrum of
$G^*$ by formal sum of formal tensor products $\Psi$ as in \eqref{psi} such that:
\begin{enumerate}
\item the $\mu_j \boxtimes R_j$ in \eqref{psi} are all distinct, and
\item for each $j$, $\mu_j^{\theta} \boxtimes R_j = \mu_j \boxtimes R_j$.
\end{enumerate}

\subsection{Local Arthur parameters} \label{S35}
Assume that $k = F_v$ is local and let $W_k'$ be its Weil-Deligne group. We can similarly define packets over $k$. We define a local packet
over $k$ as a formal sum of formal tensor product \eqref{psi} where each $\mu_j$ is now a tempered
irreducible representation of $\GL(d_j, k)$ that is square integrable modulo the center.\footnote{Because we do not know that the extension to $\GL(N)$ of Ramanujan's conjecture is valid, we do not know that the local components of automorphic representations of $\GL(N)$ are indeed tempered. So that
in principle the $\mu_j$ are not necessarily tempered: their central characters need not be unitary. This requires a minor generalization that Arthur addresses in \cite[Remark 3 p. 247]{Arthur}. Anyway the approximation to Ramanujan's conjecture proved by Luo, Rudnick and Sarnak \cite{LRS} is enough for our purposes and it makes notations easier to assume that each $\mu_j$ is indeed
tempered.} The other components $R_j$ remain irreducible representations of $\SL_2 (\C)$. 
To each $\mu_j \boxtimes R_j$ we associate the unique irreducible quotient $\Pi_i$ of 
$$\mathrm{ind} (\mu_j | \cdot |^{\frac12 (n_j-1)} , \mu_j | \cdot |^{\frac12 (n_j-3)} , \ldots , \mu_j | \cdot |^{\frac12 (1-n_j)} )$$
(normalized induction from the standard parabolic subgroup of type $(d_j, \ldots , d_j)$). We then define
$\Pi_{\Psi}$ as the induced representation
$$\mathrm{ind} (\Pi_1 \otimes \ldots \otimes \Pi_r )$$
(normalized induction from the standard parabolic subgroup of type $(n_1d_1 , \ldots , n_r d_r)$). It is 
irreducible and unitary. Finally, the local parameter consists of those representations such that 
$\Psi^{\theta} = \Psi$. Then $\Pi_{\Psi}$ is theta-stable: $\Pi_{\Psi} \circ \theta \cong \Pi_{\Psi}$.

\subsection{Local Arthur packets} Assume that $k = F_v$ is local and that $G(k) = G^* (k)$ is quasi-split. 
We now recall how Arthur associates a finite packet of representations of $G (k)$ to a local
parameter $\Psi$.

We say that two functions in $C^{\infty}_c (G (k))$ are stably equivalent if they have the same stable
orbital integrals, see e.g. \cite{LanglandsShelstad}. 
Thanks to the recent proofs by Ngo \cite{Ngo} of the fundamental lemma and Waldspurger's work \cite{waldspurger}, we have a natural notion of transfer $f \leadsto f^{G}$ from a test fonction $f \in C^{\infty}_c (\GL(N, k))$ to a representative $f^{G}$ of a stable equivalence
class of functions in $C^{\infty}_c (G(k))$ such that $f$ and $f^{G}$ are {\it associated} i.e. they have
matching stable orbital integrals~\footnote{Here the orbital integrals on the $\GL (N)$-side are 
{\it twisted} orbital integrals.}, see \cite{KS} for more details about twisted transfer. Over 
Archimedean places existence of transfer is due to Shelstad, see \cite{Shelstad08}; we note that in that
case being $K$-finite is preserved by transfer.
When $G = \SO (N , \eta)$ one must
moreover ask that $f^{G}$ is invariant under an outer automorphism $\alpha$ of $G$; we may assume that $\alpha^2=1$.

Let $\Psi$ be a local parameter as above and $\mathcal{H}_{\Pi_{\Psi}}$ be
the space of $\Pi_{\Psi}$. We fix an intertwining operator $A_{\theta} : \mathcal{H}_{\Pi_{\Psi}}
\rightarrow \mathcal{H}_{\Pi_{\Psi}}$ ($A_{\theta}^2 = 1$) intertwining $\Pi_{\Psi}$ and $\Pi_{\Psi} \circ 
\theta$. 

When $G = \SO (N , \eta)$ we identify the irreducible representations $\pi$ of $G (k)$ that 
are conjugated by $\alpha$. Then $\mathrm{trace} \pi (f^{G})$ is well defined when $f^{G}$ is as explained above.

The following proposition follows from \cite[Theorem 2.2.1]{ArthurBook} (see also \cite[Theorem 30.1]{Arthur}).

\begin{prop} \label{P37}
There exists a finite family $\small\prod (\Psi )$ of representations of $G (k)$, and some multiplicities
$m(\pi) >0$ ($\pi \in \small\prod (\Psi)$) such that, for associated $f$ and $f^{G}$~:
\begin{equation} \label{traceident}
\mathrm{trace} (\Pi_{\Psi} (f) A_{\theta}) = \sum_{\pi \in \small\prod (\Psi)} \varepsilon (\pi) m(\pi) \mathrm{trace} \ \pi (f^{G}),
\end{equation}
where each $\varepsilon (\pi)$ is a sign $\in \{ \pm 1\}$.
\end{prop}
We remark that \eqref{traceident} uniquely determines $\small\prod (\Psi)$ as a set of representations-with-multiplicities; it also uniquely determines the signs $\varepsilon (\pi)$.
In fact Arthur explicitely computes these signs for some particular choice of an intertwiner $A_{\theta}$. 

By the local Langlands correspondence, a local parameter $\Psi$ for $G (k)$ can be represented as a
homomorphism 
\begin{equation} \label{Aparam}
\Psi : W_{k}' \times \SL_2 (\C) \rightarrow {}^L \- G.
\end{equation}
Arthur associates to such a parameter the $L$-parameter $\varphi_{\Psi} : W_k ' \rightarrow {}^L \- G$
given by 
$$\varphi_{\Psi} (w) = \Psi \left( w, 
\left(
\begin{array}{cc}
|w|^{1/2} & \\
& |w|^{-1/2} 
\end{array} \right) \right).$$
One key property of the local Arthur's packet $\small\prod (\Psi)$ is that it contains 
all representations of Langlands' $L$-packet associated to $\varphi_{\Psi}$. This is proved by Arthur, see also \cite[Section 6]{Moeglin4}.

Ignoring the minor generalization needed to cover the lack of Ramanujan's conjecture, the global part of Arthur's theory (see \cite[Corollary 3.4.3]{ArthurBook} when $G$ is quasi-split and \cite[Proposition 9.5.2]{ArthurBook} in general) now implies:

\begin{prop} \label{Prop:Arthur}
Let $\pi$ be an irreducible automorphic representation of $G (\A)$ which occurs (discretely) 
as an irreducible subspace of $L^2 (G (F) \backslash G (\A))$. Then there exists a global
Arthur parameter $\Psi$ and a finite set $S$ of places of $F$ containing all Archimedean ones 
such that for all $v \notin S$, the group $G(F_v)=G^* (F_v )$ is quasi-split, the representation $\pi_v$ is unramified and the $L$-parameter of $\pi_v$ is $\varphi_{\Psi_v}$. 
\end{prop}

\medskip
\noindent
{\it Remark.} Proposition \ref{Prop:Arthur} in particular implies that the $\SL_2 (\C)$ part of a 
local Arthur parameter has a {\it global} meaning. This puts serious limitations on the kind of non-tempered
representations which can occur discretely: e.g. an automorphic representation $\pi$ of $G^* (\A)$ which occurs discretely in $L^2 (G^* (F) \backslash G^* (\A))$ and which is non-tempered at one place $v$ is non-tempered at all places. 

\medskip

The above remark explains how Arthur's theory will be used in our proof. This will be made effective through the use of $L$-functions.

\subsection{Application to $L$-functions}  Let $\pi$ be an irreducible automorphic representation of $G (\A)$ which occurs (discretely) 
as an irreducible subspace of $L^2 (G (F) \backslash G (\A))$ and let
$$\Psi = \mu_1 \boxtimes R_1 \boxplus \ldots \boxplus \mu_r \boxtimes R_r$$
be its global Arthur parameter (Proposition \ref{Prop:Arthur}). 
We factor each $\mu_j = \otimes_v \mu_{j,v}$ where $v$ runs over 
all places of $F$. Let $S$ be a finite set of places of $F$ containing the set $S$ of Proposition \ref{Prop:Arthur}, and all $v$ for which either one $\mu_{j,v}$ or $\pi_v$ is ramified.
We can then define the formal Euler product
$$L^S (s, \Pi_{\Psi})= \prod_{j=1}^r \prod_{v \not\in S} L_{v}(s- \frac{n_j -1}{2} , \mu_{j,v}) L_{v}(s- \frac{n_j -3}{2} , \mu_{j,v}) \ldots L _{v}(s- \frac{1-n_j}{2} , \mu_{j,v}).$$
 Note that 
$L^S (s , \Pi_{\Psi})$ is the 
partial $L$-function of a very special automorphic representation of $\GL (M)$ with $M= \sum_{i\in [1,r]}d(n_{i}R_{i})$; it is the product of partial $L$-functions of the square integrable automorphic representations associated to the parameters $\mu_j \boxtimes R_j$. According to Jacquet and Shalika \cite{JacquetShalika} $L^S (s , \Pi_{\Psi})$, which is a product absolutely convergent for $ \mathrm{Re}( s)\gg 0$,  extends to a meromorphic function of $s$. It moreover
follows from Proposition \ref{Prop:Arthur} and the definition of $L^S (s , \pi)$ that:
$$L^S (s , \pi ) = L^S (s , \Pi_{\Psi}).$$

\medskip
\noindent
{\it Remark.} Given an automorphic character $\eta$ we can similarly write $L^S(s, \eta \times \pi)$
as a product of $L$-functions associated to linear groups: replace $\mu_{j}$ by $\eta\otimes \mu_{j}$ in the above discussion (note that each $\mu_{j}$ is self-dual).

\medskip

We can now relates Arthur's theory to Proposition \ref{Prop:GJS}: 

\begin{lem} \label{Lem:L}
Let $\pi \in \mathcal{A}^c (G (k))$ whose global Arthur parameter $\Psi$ is a sum
$$\Psi = \left( \boxplus_{(\rho,b)}\rho \boxtimes R_{b} \right) \boxplus \eta\boxtimes R_{a}$$ 
where $\eta$ is a selfdual (quadratic) automorphic character and for each pair $(\rho,b)$, either
$b<a$ or $b=a$ and $\rho\not\simeq \eta$. Then the partial $L$-function 
$L^S (s,\eta\times \pi)$ --- here $S$ is a finite set of places containing the set $S$ of Proposition \ref{Prop:Arthur} and
all the places where $\eta$ ramifies --- is holomorphic in the half-plane $\mathrm{Re} (s) > (a+1)/2$ and
it has a simple pole in $s=(a+1)/2$.
\end{lem}
\begin{proof} Writing $L^S(s,\eta\times \pi)$ explicitely on a right half-plane of absolute convergence; 
we get a product of $L^S (s-(a-1)/2, \eta \times \eta)$ by factors $L^S (s-(b-1)/2,\eta\times \rho)$. Our hypothesis on $a$ forces $b\leq a$ and if $b=a$, $\rho\not\simeq \eta$.
The conclusion of the lemma follows.
\end{proof}

 
\subsection{Infinitesimal character} \label{212}
Let $\pi_{v_0}$ be the local Archimedean factor of representation $\pi  \in  \mathcal{A}^c (G)$ with 
global Arthur parameter $\Psi$. We may associate to $\Psi$ the parameter $\varphi_{\Psi_{v_0}} : 
\C^* \rightarrow G^{\vee} \subset \GL (N , \C)$ given by 
$$\varphi_{\Psi_{v_0}} (z) = \Psi_{v_0} \left( z , \left(
\begin{array}{cc}
(z \overline{z})^{1/2} & \\
& (z\overline{z})^{-1/2} 
\end{array} \right) \right).$$
Being semisimple, it is conjugate into the maximal torus 
$$T^{\vee} = \{ \mathrm{diag} (x_1 , \ldots , x_{\ell},
x_{\ell}^{-1} , \ldots , x_1^{-1} ) \}$$ 
of $G^{\vee}$. We may therefore write $\varphi_{\Psi_{v_0}} = 
(\eta_1 , \ldots , \eta_{\ell} , \eta_{\ell}^{-1} , \ldots , \eta_1^{-1})$ where each $\eta_j$ is a character $z \mapsto z^{P_i} \overline{z}^{Q_i}$. One easily checks that the vector 
$$\nu_{\Psi} = (P_1 , \ldots , P_{\ell}) \in \C^{\ell} \cong \mathrm{Lie}(T) \otimes \C$$
is uniquely defined modulo the action of the Weyl group $W$ of $G(F_{v_0})$. The following proposition
is detailed in \cite{BC1}.

\begin{prop}
The infinitesimal character of $\pi_{v_0}$ is the image of $\nu_{\Psi}$ in $\C^{\ell} /W$.
\end{prop}

Recall that the infinitesimal character $\nu_{\Psi}$ is said to be {\it regular} if the $P_j$ are all distinct.

\section{A surjectivity theorem for theta liftings}

\subsection{Notations} \label{3.1} 
Let $F$ be a number field and $\A$ be its ring of adeles. 
Fix $\psi$ a non-trivial additive character of $\A /F$. 
Let $V$ be a nondegenerate quadratic space over $F$ with $\dim_F V =m$. 
We set $\ell=[m/2]$ and $N=2\ell$. 

We say that a representation  $\pi \in \mathcal{A}^c (\SO (V))$ --- i.e. an irreducible cuspidal automorphic representation of $\SO (V)$ --- is in the image of the cuspidal $\psi$-theta correspondence from a smaller
group if there exists a symplectic space $X$ with $\dim X \leq N$ and a lift $\tilde{\pi}$ of $\pi$ to $\OO (V)$ such that $\tilde{\pi}$ is the image of a cuspidal automorphic form of $\Mp (X)$ by the $\psi$-theta correspondence. We note that then the restriction of 
$\tilde{\pi}$ to $\SO (V)$ is isomorphic to $\pi$. In fact --- but we don't need it here --- the representation
$\tilde{\pi}$ should be uniquely determined here. At each finite place this follows from Kudla-Rallis' theta dichotomy \cite{KudlaRallis2} and there is little doubt that this still holds at Archimedean places (see Kashiwara-Vergne \cite{KashiwaraVergne} for the case of compact places).  

Here we combine propositions \ref{Prop:GJS} and \ref{Prop:Arthur} to prove the following result
which is the main automorphic ingredient in our work. 

We say that a global representation 
$\pi \in \mathcal{A}^c (\SO (V))$ is {\it highly non-tempered} if its global Arthur
parameter $\Psi$ contains a factor $\eta \boxtimes R_a$ where $\eta$ is a quadratic character and $3a > m-1$.

\begin{thm} \label{Thm1}
Let $\pi \in \mathcal{A}^c (\SO (V))$. Assume that $\pi$ is highly non-tempered and that $\pi_{v_0}$ has a regular infinitesimal character. Then there exists an automorphic quadratic character $\chi$ such that $\pi\otimes \chi$ is in the image of the cuspidal $\psi$-theta correspondence from a smaller group associated to a symplectic space of dimension $m-a-1$.
\end{thm}
\begin{proof} Let $G= \SO (V)$ and let 
$\tilde{\pi} \in \mathcal{A}^c (\OO (V))$ be a lift of $\pi$. Recall from Remark \ref{1.6} that the partial 
$L$-function $L^S ( s ,\tilde{\pi}) = L^S (s , \pi)$.

\begin{lem} \label{LL}
The global Arthur parameter of $\pi$ is a sum 
$$\Psi = \left( \boxplus_{(\rho,b)}\rho \boxtimes R_{b} \right) \boxplus \eta\boxtimes R_{a}$$ 
where $\eta$ is a selfdual (quadratic) automorphic character and each pair $(\rho,b)$ 
consists of a (selfdual) cuspidal automorphic representation $\rho$ of some 
$\GL(d_{\rho})$ and a positive integer $b<a$ such that 
$\sum_{(\rho,b)}bd_{\rho}+a=N$. 
\end{lem}
\begin{proof} The global Arthur parameter of $\pi$ is a sum
$\Psi = \eta \boxtimes R_a \boxplus \boxplus_{(\rho,b)}\rho \boxtimes R_{b}$ (without multiplicity). Then 
$$
a+\sum_{(\rho ,b )}d_{\rho}b=N.
$$
We want to prove that in the second case there is no $b\geq a$ occuring in the second sum. 
Indeed: suppose by contradiction that such a $(\rho  , b)$ occurs. We first note that $a \equiv m-1$ mod $2$ as the representations 
$\eta\otimes R_{a}$ must be 
orthogonal (resp. symplectic) if the dual group of $\SO(V)$ is orthogonal (resp. symplectic). In particular $a$ is odd if $m$ is even. The hypothesis $3a > m-1$ therefore implies that $3a >N$ so that$d_{\rho} = 1$. According to \S \ref{212} 
(since $\eta$ is a quadratic character) $(\eta , a)$ contributes to the infinitesimal character of $\pi_{v_0}$ by 
$\left( \frac{a-1}{2} , \frac{a-3}{2} , \ldots , 0 \right)$ or $\left( \frac{a-1}{2} , \frac{a-3}{2} , \ldots , \frac12 \right)$ according to the parity of $a$. Now since $d_{\rho }=1$ the automorphic representation $\rho$
is also a quadratic character and $(\rho , b)$ contributes to the infinitesimal character of $\pi_{v_0}$ by 
$\left( \frac{b-1}{2} , \frac{b-3}{2} , \ldots , 0 \right)$ or $\left( \frac{b-1}{2} , \frac{b-3}{2} , \ldots , \frac12 \right)$ according to the parity of $b$. Finally since $a , b \equiv m-1$ mod $2$, both
$a$ and $b$ have the same parity and the infinitesimal character of $\pi_{v_0}$ cannot be regular in contradiction with our hypothesis. 
\end{proof}

\medskip
\noindent
{\it Remarks.} 1. The stronger hypothesis $a> \ell = N/2$ (without any hypothesis on the infinitesimal
character of $\pi_{v_0}$) directly implies that $b<a$ if $(\rho,b)$ appears in the parameter of 
$\pi$. 

2. The proof shows that the conclusion of Theorem \ref{Thm1} still holds if the \emph{local} Arthur parameter $\Psi_{v_0}$ contains a factor $\eta \boxtimes R_a$ where 
$\eta$ is a quadratic character and $3a > m-1$.

\subsection{} Let $\eta$ be the automorphic character given by Lemma \ref{LL}. 
Then Lemma \ref{Lem:L} implies that for some finite set of places $S$, the partial $L$-function 
$L^S( s, \pi \times \eta)$ is holomorphic in the half-plane $\mathrm{Re} (s) > (a+1)/2$ and has a 
simple pole in $s= (a+1)/2$. Adding a finite set of places to $S$ we may assume that this also
holds for the partial $L$-function $L^S (s , \tilde{\pi} \times \eta )$. 

Now let $p= \frac12 (m-a-1)$ and $X$ be a symplectic $F$-space with $\dim X = 2p$.
Proposition \ref{Prop:GJS} implies that there exists an automorphic sign
character $\epsilon$ of $\OO_m (\A)$ such that the $\psi^{-1}$-theta lifting 
of $(\tilde{\pi} \otimes \eta) \otimes \epsilon$ to $\Mp_{2p} (\A)$ does not vanish. 

\subsection{} Here we prove that 
$\pi ' := \Theta_{\psi^{-1} , V}^X ((\tilde{\pi} \otimes \eta) \otimes \epsilon)$ is cuspidal.
Let $\pi_0 '$ be the first (non-zero) occurrence of the 
$\psi^{-1}$-theta lifting of $(\tilde{\pi} \otimes \eta) \otimes \epsilon$ in the Witt tower of the symplectic spaces. By the Rallis theta tower property \cite{Rallis}, then $\pi_0'$ is cuspidal.
Let $2p_0$ be the dimension of the symplectic space corresponding to $\pi_0'$. We want to prove
that $p_0 = p$. By the unramified correspondence we know 
$\tilde{\pi}_v$ in all but finitely many places $v$. It correspond to $\tilde{\pi}_v$ a {\it unique} 
(see e.g. \cite{M3}) Arthur packet $\Psi_v$ which contains $\pi_v$ --- the restriction of $\pi_v$ to $\SO(V)(F_v)$. And $\Psi_v$ contains a factor $\eta_v \otimes R_{a'}$ with $a' = m-2p_0 -1$. In particular $a' \geq a$. As in the proof of Lemma \ref{Lem:L} writing explicitely the partial $L$-function 
$L^S (s , \pi \times \eta )$ on a right half-plane of absolute convergence, we get a product of 
$L^S (s- (a'-1)/2 , \eta \times \eta)$ by factors which are holomorphic in $(a'+1)/2$. This forces
the partial $L$-function $L^S (s , \pi \times \eta )$ to have a pole in $s= (a'+1)/2$ and it follows from Lemma \ref{Lem:L} that $a' \leq a$. Finally
$a=a'$ and $p_0 = m-a-1$.

\subsection{} The main theorem of \cite{Moeglin97b} and \cite[Theorem 1.2]{JiangSoudry} now apply
to the representation $(\tilde{\pi} \otimes \eta) \otimes \epsilon$ to show that 
$$\Theta_{\psi, X}^V (\Theta_{\psi^{-1} , V}^X ((\tilde{\pi} \otimes \eta) \otimes \epsilon ) ) = (\tilde{\pi} \otimes \eta) \otimes \epsilon.$$
In otherwords: $\Theta_{\psi, X}^V (\pi ') = (\tilde{\pi} \otimes \eta) \otimes \epsilon$. This concludes the proof of Theorem \ref{Thm1} with $\chi = \eta \otimes \epsilon$.
\end{proof}

\part{Local computations}

\section{Cohomological unitary representations} \label{sec:4}

\subsection{Notations} Let $p$ and $q$ two non-negative integers with $p+q =m$. In this section
$G= \SO_0 (p,q)$ and $K = \SO (p ) \times \SO (q)$ is a maximal compact subgroup of $G$. 
We let $\g_0$ the real Lie algebra of $G$ and $\g_0 = \k_0 \oplus \p_0$ be the Cartan decomposition associated to the choice of the maximal
compact subgroup $K$. We denote by $\theta$ the corresponding Cartan involution. If $\mathfrak{l}_0$ is a real Lie algebra we denote by $\mathfrak{l}$ its complexification $\mathfrak{l} = \mathfrak{l}_0 \otimes \C$.   

\subsection{Cohomological $(\g , K)$-modules} Let $(\pi , V_{\pi})$ be an irreducible unitary
$(\g , K )$-module and $E$ be a finite dimensional irreducible representation of $G$. We say that $(\pi , V_{\pi})$ is {\it cohomological (w.r.t. the local system associated to $E$)} if it has nonzero $(\g , K)$-cohomology $H^{\bullet} (\g , K ; V_{\pi} \otimes E)$.

Cohomological $(\g , K)$-modules are classified by Vogan and Zuckerman in \cite{VZ}: Let $\mathfrak{t}_0$ be a Cartan subalgebra of $\k_0$. A {\it $\theta$-stable parabolic subalgebra} $\mathfrak{q} = \mathfrak{q} (X) \subset \g$ is associated to an element $X \in i \mathfrak{t}_0$. It is defined as the direct sum
$$\mathfrak{q} = \mathfrak{l} \oplus \mathfrak{u},$$
of the centralizer $\mathfrak{l}$ of $X$ and the sum $\mathfrak{u}$ of the positive eigenspaces of
$\mathrm{ad} (X)$. Since $\theta X = X$, the subspaces $\mathfrak{q}$, $\mathfrak{l}$ and $\mathfrak{u}$ are all invariant under $\theta$, so 
$$\mathfrak{q} = \mathfrak{q} \cap \k \oplus \mathfrak{q} \cap \p,$$
and so on. 

The Lie algebra $\mathfrak{l}$ is the complexification of $\mathfrak{l}_0 = \mathfrak{l} \cap \mathfrak{g}_0$. Let $L$ be the connected subgroup of $G$ with Lie algebra $\mathfrak{l}_0$. Fix positive system $\Delta^+(\mathfrak{l})$ of the roots of $\mathfrak{t}$ in $\mathfrak{l}$. Then $\Delta^+ (\g) = \Delta^+ (\mathfrak{l}) \cup \Delta (\mathfrak{u})$ is a positive system of the roots of $\mathfrak{t}$ in $\g$. 
Now extend $\mathfrak{t}$ to a Cartan subalgebra $\mathfrak{h}$ of $\g$, and choose $\Delta^+ (\g , \mathfrak{h})$ a positive system of roots of $\mathfrak{h}$ in $\g$ such that its restriction to $\mathfrak{t}$ gives $\Delta^+ (\g)$. Let 
$\rho$ be half the sum of the roots in $\Delta^+ (\mathfrak{h} , \g)$ and $\rho (\mathfrak{u} \cap \p)$ half the sum of the roots in $\mathfrak{u} \cap \p$. 
A one-dimensional representation $\lambda : \mathfrak{l} \rightarrow \C$ is admissible if it satisfies the following two conditions: 
\begin{enumerate}
\item $\lambda$ is the differential of a unitary character of $L$,
\item if $\alpha \in  \Delta (\mathfrak{u})$, then $\langle \alpha , \lambda_{\mathfrak{t}} \rangle \geq 0$.
\end{enumerate}
Given $\mathfrak{q}$ and an admissible $\lambda$, let $\mu (\mathfrak{q} , \lambda)$
be the representation of $K$ of highest weight $\lambda_{| \mathfrak{t}} + 2 \rho (\mathfrak{u} \cap \p)$.

The following proposition is due to Vogan and Zuckerman \cite[Theorem 5.3 and Proposition 6.1]{VZ}.

\begin{prop} \label{VZ}
Assume that $\lambda$ is zero on the orthogonal complement of $\mathfrak{t}$ in $\mathfrak{h}$. There
exists a unique irreducible unitary $(\g , K)$-module $A_{\mathfrak{q}} (\lambda )$ such that:
\begin{enumerate}
\item  $A_{\mathfrak{q}} (\lambda )$ contains the $K$-type $\mu (\mathfrak{q} , \lambda )$.
\item  $A_{\mathfrak{q}} (\lambda )$ has infinitesimal character $\lambda + \rho$. 
\end{enumerate}
\end{prop}

Vogan and Zuckerman (see \cite[Theorem 5.5 and 5.6]{VZ}) moreover prove:

\begin{prop} \label{VZ2}
Let $(\pi , V_{\pi})$ be an irreducible unitary
$(\g , K )$-module and $E$ be a finite dimensional irreducible representation of $G$. Suppose 
$H^{\bullet} (\g , K ; V_{\pi} \otimes E)\neq 0$. Then there is a $\theta$-stable parabolic subalgebra
$\mathfrak{q} = \mathfrak{l} \oplus \mathfrak{u}$ of $\g$, such that:
\begin{enumerate}
\item $E/\mathfrak{u} E$ is a one-dimensional unitary representation of $L$; write $-\lambda :
\mathfrak{l} \rightarrow \C$ for its differential.
\item $\pi \cong A_{\mathfrak{q}} (\lambda )$.
Moreover, letting $R = \dim (\mathfrak{u} \cap \p)$,  we have:
\begin{eqnarray*}
H^{\bullet} (\g , K ; V_{\pi} \otimes E) & \cong & H^{\bullet - R} (\mathfrak{l} , \mathfrak{l} \cap \k , \C) \\
& \cong & \mathrm{Hom}_{\mathfrak{l} \cap \mathfrak{k}} (\wedge^{\bullet -R} (\mathfrak{l} \cap \p) , \C).
\end{eqnarray*}
\item As a consequence of (2) we have
\begin{equation}\label{Rcohomologyhasdimensionone}
H^R (\g , K ; V_{\pi} \otimes E) \cong \mathrm{Hom}_{\mathfrak{l} \cap \mathfrak{k}} (\wedge^{0} (\mathfrak{l} \cap \p) , \C) \cong \C
\end{equation}
\end{enumerate}
\end{prop}
 We can now prove the following proposition.

\begin{prop}\label{occursonce}
The representation $\mu(\mathfrak{q},\lambda)$ of $K$ has the following properties
\begin{enumerate}
\item $\mu(\mathfrak{q},\lambda)$ is the only representation of $K$ common to both $\wedge^{R}(\p)\otimes E^*$ and $A_{\mathfrak{q}} (\lambda )$.
\item $\mu(\mathfrak{q},\lambda)$ occurs with multiplicity one in $\wedge^{R}(\p)\otimes E^*$.
\item $\mu(\mathfrak{q},\lambda)$  occurs with multiplicity one in $A_{\mathfrak{q}} (\lambda )$.
\end{enumerate}
\end{prop}
\begin{proof}
From \cite{VZ}, Proposition 5.4(c) we have
$$H^R (\g , K ; V_{\pi} \otimes E) \cong \mathrm{Hom}_K( \wedge^{R}(\p), A_{\mathfrak{q}} (\lambda ) \otimes E).$$
But combining this equation with  (e) of the proposition above we conclude
$$ \mathrm{dim} \  \mathrm{Hom}_K( \wedge^{R}(\p), A_{\mathfrak{q}} (\lambda ) \otimes E) = 1.$$
This last equation implies all three statements of the proposition.
\end{proof}
Let $e(\mathfrak{q})$ be a generator of the line $\wedge^R (\mathfrak{u} \cap \p)$. Then $e(\mathfrak{q})$
is the highest weight vector of an irreducible representation $V(\mathfrak{q})$ of $K$ contained in 
$\wedge^R \p$ (and whose highest weight is thus necessarily $2\rho (\mathfrak{u} \cap \p)$). It follows that $V(\mathfrak{q})$ is the unique occurrence of $\mu(\mathfrak{q})$ in $\wedge^{R}(\p)$. We will refer to the special $K$-types $\mu (\mathfrak{q})$ as {\it Vogan-Zuckerman $K$-types}. Let $V(\mathfrak{q},\lambda)$ denote the Cartan product of $V(\mathfrak{q})$ and 
$E^*$ (this means the irreducible submodule of the tensor product $V(\mathfrak{q} ) \otimes E^*$ with highest weight the sum
$2\rho(\mathfrak{u} \cap \p) + \lambda$).  By definition $V(\mathfrak{q},\lambda)$ occurs in $\wedge^R (\mathfrak{u} \cap \p) \otimes E^*$
and hence, by (2) of Proposition \ref{occursonce} above it is the unique copy of $\mu(\mathfrak{q},\lambda)$ in $\wedge^R (\mathfrak{u} \cap \p) \otimes E^*$.
From the discussion immediately above and Propostion \ref{occursonce} we obtain

\begin{cor} \label{VZKtype}
Any nonzero element  $\omega \in \mathrm{Hom}_K ( \wedge^{R} \p \otimes E^*, A_{\mathfrak{q}} (\lambda ))$ factors through the isotypic component 
$V(\mathfrak{q},\lambda)$.    
\end{cor}

\subsection{} We will make geometric use of the isomorphism of Proposition \ref{VZ2}(2). In doing so we will need the following lemmas which are essentially due to Venkataramana \cite{Venky}. 

We let $T$ be the torus of $K$ whose Lie algebra is $i \mathfrak{t}_0$. The action of $T$ on the space $\p$ is completely reducible and we have a decomposition 
$$\p = (\mathfrak{u} \cap \p) \oplus (\mathfrak{l} \cap \p) \oplus (\mathfrak{u}^- \cap \p),$$
where the element $X \in \mathfrak{t}$ acts by strictly positive (resp. negative) eigenvalues on $\mathfrak{u} \cap \p$ (resp. $\mathfrak{u}^- \cap \p$) and by zero eigenvalue on 
$\mathfrak{l} \cap \p$. Now using the Killing form, the inclusion map $\mathfrak{l} \cap \p \rightarrow \p$ induces a {\it restriction} map
$\p \rightarrow \mathfrak{l} \cap \p$ and we have the following:

\begin{lem} \label{V1.3}
Consider the restriction map $B: [\wedge^{\bullet} \p]^T \rightarrow [\wedge^{\bullet} (\mathfrak{l} \cap \p)]^T$ and the cup-product map $A : [\wedge^{\bullet} \p ]^T \rightarrow \wedge^{\bullet} \p$ given by $y \mapsto y \wedge 
e( \mathfrak{q})$. Then the kernels of $A$ and $B$ are the same.
\end{lem}
\begin{proof} This is \cite[Lemma 1.3]{Venky}. Note that although it is only stated there for Hermitian symmetric spaces, the proof goes through without any modification.
\end{proof}

Now consider the restriction map 
\begin{equation} \label{resmap10}
[\wedge^{\bullet} \p ]^K \rightarrow [ \wedge^{\bullet} (\mathfrak{l} \cap \p) ]^{K \cap L}.
\end{equation}
An element $c \in [\wedge^{\bullet} \p ]^K$ defines --- by cup-product --- a linear map in 
$$\mathrm{Hom}_K (\wedge^{\bullet} \p , \wedge^{\bullet} \p)$$
that we still denote by $c$. The following lemma is essentially the same as \cite[Lemma 1.4]{Venky}.

\begin{lem} \label{V1.4}
Let $c \in [\wedge^{\bullet} \p ]^K$. Then we have:
$$c (V( \mathfrak{q})) = 0 \Leftrightarrow c \in \mathrm{Ker} \left( [\wedge^{\bullet} \p ]^K \rightarrow [ \wedge^{\bullet} (\mathfrak{l} \cap \p) ]^{K \cap L} \right).$$
\end{lem}
\begin{proof} As a $K$-module $V(\mathfrak{q})$ is generated by $e(\mathfrak{q})$. We therefore deduce from the $K$-invariance of $c$ that:
$$c (V( \mathfrak{q})) = 0 \Leftrightarrow c  \wedge e(\mathfrak{q}) =0.$$
The second equation  is equivalent to the fact that $c$ belongs to the kernel of the map $B$ of Lemma \ref{V1.3}. But $B ( c )=0$ if and only if $c$ belongs to the kernel of the restriction map
$$[\wedge^{\bullet} \p ]^K \rightarrow [ \wedge^{\bullet} (\mathfrak{l} \cap \p) ]^{K \cap L}.$$
This concludes the proof. 
\end{proof}



\subsection{} \label{par:weight} In the notation of \cite{Bourbaki}, we may choose a Killing-orthogonal basis $\varepsilon_i$ of $\mathfrak{h}^*$ such that the positive roots are those roots $\varepsilon_i \pm \varepsilon_j$ with $1 \leq i < j \leq \ell$ as well as the roots $\varepsilon_i$ ($1 \leq i \leq \ell$) if $m$ is odd. The finite dimensional irreducible representations of $G$ are parametrized by a highest weight
$\lambda = (\lambda_1 , \ldots , \lambda_{\ell}) = \lambda_1 \varepsilon_1 + \ldots + \lambda_{\ell} \varepsilon_{\ell}$ such that $\lambda$ is dominant (i.e. $\lambda_1 \geq \ldots \geq \lambda_{\ell-1} \geq |\lambda_{\ell}|$ and $\lambda_{\ell} \geq 0$ if $m$ is odd) and integral (i.e. every $\lambda_i \in \Z)$. 

In the applications we will mainly be interested in the following examples.

\medskip
\noindent
{\it Examples.} 1. The group $G= \SO_0 (n,1)$ and $\lambda = 0$. Then for each integer $q= 0, \ldots , \ell -1$ the $K$-representation $\wedge^q \mathfrak{p}$ is just $\wedge^q \C^n$ with $K=\SO (n)$; it is
irreducible and we denote it $\tau_q$. In addition, if $n=2\ell$, $\wedge^{\ell} \mathfrak{p}$ decomposes as a sum of two irreducible representations $\tau_{\ell}^+$ and $\tau_{\ell}^-$. From this we get that for each integer $q= 0, \ldots , \ell -1$ there exists exactly one irreducible $(\g , K)$-module $(\pi_q , V_q)$
such that $H^q (\g , K ; V_{q}) \neq 0$. In addition, if $n=2\ell$, there exists two irreducible $(\g , K)$-module $(\pi_q^{\pm} , V_q^{\pm})$ such that $H^{\ell} (\g , K ; V_{\ell}^{\pm}) \neq 0$. Moreover:
$$H^k (\g , K ; V_{q}) = \left\{
\begin{split}
0 & \mbox{ if } k \neq q, n-q , \\
\C & \mbox{ if } k = q \mbox{ or } n-q
\end{split} \right.$$
and, if $n= 2\ell$, 
$$H^k (\g , K ; V_{\ell}^{\pm}) = \left\{
\begin{split}
0 & \mbox{ if } k \neq \ell , \\
\C & \mbox{ if } k =  \ell.
\end{split} \right.$$
The Levi subgroup $L \subset G$ associated to $(\pi_q , V_q)$ ($q=0 , \ldots , \ell -1$) is 
$L = C \times \SO_0 (n-2q , 1)$ where $C \subset K$.  

2. The group $G = \SO_0 (n,1)$ and $\lambda = (1 , 0 , \ldots , 0)$ is the highest weight of its standard representation in $\C^{m}$. Then there exists a unique $(\g, K)$-module $(\pi , V_{\pi})$ such that 
$H^1 (\g , K ; V_{\pi}  \otimes \C^m ) \neq 0$. The  Levi subgroup $L \subset G$ associated to $(\pi , V_{\pi})$ is $L = C \times \SO_0 (n-2, 1)$ where $C \subset K$

3. The group $G = \SO_0 (n,2)$ and $\lambda = 0$. Then $K= \SO (n) \times \SO (2)$ and $\p = \C^n \otimes (\C^2)^*$ where $\C^n$ (resp. $\C^2$) is the standard representation of $\SO (n)$ (resp. $\SO(2)$). We denote by $\C^+$ and $\C^-$ the $\C$-span of the vectors $e_1+ie_2$ and $e_1-ie_2$ in 
$\C^2$. The two lines $\C^+$ and $\C^-$ are left stable by $\SO(2)$. This yields a decomposition 
$\p = \p^+ \oplus \p^-$ which corresponds to the decomposition given by the natural complex structure on $\p_0$. For each non-negative integer $q$ the $K$-representation $\wedge^q \p =
\wedge^q (\p^+ \oplus \p^-)$ decomposes as the sum:
$$\wedge^q \p = \bigoplus_{a+b=q}  \wedge^a \p^+ \otimes \wedge^b \p^-.$$
The $K$-representations $\wedge^a \p^+ \otimes \wedge^b \p^-$ are not irreducible in general: there is 
at least a further splitting given by the Lefschetz decomposition:
$$\wedge^a \p^+ \otimes \wedge^b \p^- = \bigoplus_{k=0}^{\min (a,b)} \tau_{a-k , b-k}.$$
One can check that for $2(a+b) <n$ each $K$-representation $\tau_{a,b}$ is  irreducible. 
Moreover: only those with $a=b$ can occur as a $K$-type of a cohomological module. 
For each non-negative integer $r$ such that $4r<n$ there thus exists exactly one cohomological
module $A_{r,r}$ such that:
$$H^q (\g , K ; A_{r,r}) = \left\{ 
\begin{split}
\C & \mbox{ if } q = 2r+2k \ ( 0 \leq k \leq n-2r), \\
0 & \mbox{ otherwise}.
\end{split}\right.
$$
Moreover: $H^{2r} (\g , K ; A_{r,r}) = H^{r,r} (\g , K ; A_{r,r})$. The Levi subgroup $L \subset G$ associated to $A_{r,r}$ ($4r< n$) is $L = C \times \SO_0 (n-2r , 2)$ where $C \subset K$.  
See e.g. \cite[\S 1.5]{HarrisLi} for more details.

\medskip

\subsection{} We now consider the general case $G = \SO_0 (p,q)= \SO_0 (V)$, where $V$ is a real quadratic space of dimension $m$ and signature $(p,q)$ two non-negative integers with $p+q=m$. We denote by $(,)$ the non-degenerate quadratic form on $V$ and let $v_{\alpha}$, $\alpha = 1 ,\ldots , p$, $v_{\mu}$, $\mu = p+1 , \ldots , q$, be an orthogonal basis of $V$ such that $(v_{\alpha} , v_{\alpha})=1$ and $(v_{\mu} , v_{\mu})=-1$. We denote by $V_+$ (resp. $V_-$) is the span of 
$\{ v_{\alpha} \; : \; 1¬¨‚?†\leq \alpha \leq p \}$ (resp. $\{v_{\mu} \; : \; p+1 \leq \mu \leq m \}$).
As a representation of $\SO(p) \times \SO (q) = \SO (V_+ ) \times \SO (V_-)$, the space $\p$ is isomorphic to $V_+ \otimes (V_-)^*$. 

First recall that, as a $\GL (V_+) \times \GL (V_-)$-module, we have (see \cite[Equation (19), p. 121]{Fulton}):
\begin{equation} \label{decGL}
\wedge^R ( V_+ \otimes V_-^*) \cong \bigoplus_{\mu \vdash R} S_{\mu} (V_+) \otimes S_{\mu^*} (V_-)^*.
\end{equation}
Here $S_{\mu} (\cdot )$ denotes the Schur functor (see \cite{FultonHarris}), we sum over all partition of $R$ (equivalently Young diagram of size $|\mu|=R$) and $\mu^*$ is the conjugate partition (or transposed Young diagram). 

We will see that as far as we are concerned with special cycles, we only have to consider the decomposition of the 
submodule $\wedge^R (V_+ \otimes V_-^*)^{\SL (V_-)}$. Then each Young diagram $\mu$ which 
occurs in \eqref{decGL} is of type $\mu = (q , \ldots , q)$. 

Following \cite[p. 296]{FultonHarris}, we may define
the harmonic Schur functor $S_{[\mu]} (V_+)$ as the image of $S_{\mu} (V_+)$ by the $\SO (V_+)$-equivariant projection of  $V_+^{\otimes R}$ onto the harmonic tensors.
The representation $S_{[\mu]} (V_+)$ is irreducible with  highest weight  $\mu$.
From now on we suppose that $\mu$ has at most $\frac{p}{2}$ (positive) parts. Then Littlewood gave a formula for the decomposition of $S_{\mu} (V_+ )$ as a representation of $\SO (V_+)$ by restriction (see \cite[Eq. (25.37), p. 427]{FultonHarris}):

\begin{prop} \label{Lit}
The multiplicity of the finite dimensional $\SO (V_+)$-representation $S_{[\nu]} (V_+ )$ in $S_{\mu} (V_+ )$ equals 
$$\sum_{\xi} \dim \mathrm{Hom}_{\GL (V_+)} (S_{\mu} (V_+) , S_{\nu} (V_+) \otimes S_{\xi} (V_+)),$$
where the sum is over all nonnegative integer partitions $\xi$ with rows of {\rm even} length.
\end{prop}

\subsubsection{The Euler form} \label{eulerform} It is well known (see e.g. \cite[Theorem 5.3.3]{GoodmanWallach}) that
$[\mathrm{Sym}^q (V_+)]^{\SO (V_+)}$ is trivial if $q$ is odd and $1$-dimensional generated by  
\begin{equation} \label{euler1}
\sum_{\sigma \in \mathfrak{S}_{q}} \sigma \cdot \theta_q
\end{equation}
where $\theta = \sum_{\alpha = 1}^p v_{\alpha} \otimes v_{\alpha}$ and
$$\theta_q = \underbrace{\theta \otimes \cdots \otimes \theta}_{\ell} = \sum_{\alpha_1 , \ldots , \alpha_{\ell}} v_{\alpha_1} \otimes v_{\alpha_1} \otimes \ldots \otimes v_{\alpha_{\ell}} \otimes v_{\alpha_{\ell}},$$
if $q=2\ell$ is even.   
Note that 
$$\wedge^q (V_+ \otimes V_-^* )^{\SO (V_+) \times \SL (V_-)} \cong [\mathrm{Sym}^q (V_+)]^{\SO (V_+)} \otimes \wedge^q (V_-)^*.$$
It is therefore trivial if $q$ is odd and $1$-dimensional if $q$ is even. Using the isomorphism of \eqref{euler1} we obtain a generator of $[\wedge^q \p]^{\SO (V_+) \times \SL(V_-)}$ as the image of $\sum_{\sigma \in \mathfrak{S}_{q}} \sigma \cdot \theta_q$ under the above isomorphism.  The associated invariant $q$-form on $D$ is called  the
{\it Euler form} $e_q$. The Euler form $e_q$ is zero if $q$ is odd and  for $q=2\ell$  is expressible in terms of the curvature two-forms $\Omega_{\mu , \nu} = \sum_{\alpha=1}^p (v_{\alpha} \otimes v_{\mu}^*) \wedge (v_{\alpha} \otimes v_{\nu}^*)$ by the formula
\begin{equation} \label{euler2}
e_q = \sum_{\sigma \in \mathfrak{S}_{q}} \mathrm{sgn} (\sigma ) \Omega_{p+\sigma (1) ,p + \sigma (2)} \wedge \ldots \wedge \Omega_{p+\sigma (2\ell
-1) ,p+ \sigma (2\ell )} \in \wedge^q \p.
\end{equation}

\medskip

As in the above example it more generally follows from  \cite[Theorem 5.3.3]{GoodmanWallach} that we have:

\begin{prop} \label{invforms}
The subspace $[\wedge^{\bullet} \p]^{\SO (V_+) \times \SL(V_-)}$ is the subring of $\wedge^{\bullet} \p$ generated by 
the Euler class $e_q$.
\end{prop}

\medskip
\noindent
{\it Remark.} Proposition \ref{invforms} implies that 
$$[\wedge^{nq} \p]^{\SO (V_+) \times \SL(V_-)} = \C \cdot e_q^n.$$
It is consistent with Proposition \ref{Lit} since one obviously have:
$$\sum_{\xi} \dim \mathrm{Hom}_{\GL (V_+)} (S_{n \times q} (V_+) , S_{\xi} (V_+)) =\left\{ 
\begin{array}{ll}
0 & \mbox{ if } q \mbox{ is odd}\\
1 & \mbox{ if } q \mbox{ is even} 
\end{array} \right. .$$

\medskip

\subsection{} Let $V_r \subset \wedge^{rq} \p$ ($0 \leq r \leq p/2$) denote the realization of the irreducible $K$-type $\mu_r$ isomorphic to
$$S_{[r \times q]} (V_+ ) \otimes (\wedge^q V_-)^r \subset \wedge^{rq} (V_+ \otimes V_-^* ).$$

\medskip
\noindent
{\it Remark.} The $K$-type $\mu_r$ is the Vogan-Zuckerman $K$-type $\mu (\mathfrak{q}_r)$ where $q_r$ is any $\theta$-stable parabolic subalgebra with corresponding Levi subgroup $L= C \times \SO_0 (p-2r , q)$ with $C \subset K$. In particular we may as well the maximal one where $L=\mathrm{U} ( r ) \times \SO (p-2r,q)$, see \S \ref{VZK} where $\mu_r$ is analysed in details. 

\medskip
Wedging with the Euler class defines a linear map in
$$\mathrm{Hom}_{\SO (V_+) \times \SL (V_-)}( \wedge^{\bullet} \p , \wedge^{\bullet} \p).$$
We still denote by $e_q$ the linear map. Note that under the restriction map \eqref{resmap10} to $\mathfrak{l} \cap \p$ the Euler class $e_q$ restricts to the Euler class in
$\wedge^{q} (\mathfrak{l} \cap \p )$. It then follows from Lemma \ref{V1.4} that if $q$ is even
$e_q^k (V_r)$ is a {\it non-trivial} $K$-type in $\wedge^{(r+k)q} \p$ if and only if $k \leq p-2r$. This leads to the following:

\begin{prop} \label{PropLit}
The irreducible $K$-types $\mu_r$ are the only Vogan-Zuckerman $K$-types that occur in the subring
$$[\wedge^{\bullet} \p ]^{\SL (V_-)} = \oplus_{n=0}^p S_{n \times q} (V_+) \otimes (\wedge^q V_-)^n.$$
Moreover:
$$\mathrm{Hom}_{\SO (V_+) \times \SO (V_-)} \left(V _r ,[ \wedge^{nq} \p ]^{\SL (V_-) }\right) = \left\{ \begin{array}{ll}
\C \cdot e_{q}^{n-r} & \mbox{ if } n=r , \ldots , p-r \\
0 & \mbox{ otherwise}. 
\end{array} \right.$$
\end{prop}
\begin{proof} Cohomological representations of $G=\SO_0 (p,q)$, or equivalently their lowest $K$-types $\mu (\mathfrak{q})$, are parametrized in \cite[\S 2.2]{MSMF} by certain types of Young diagram $\nu$ so that 
$$\mu (\mathfrak{q} ) \cong S_{[\nu]} (V_+) \otimes S_{[\nu^* ]} (V_-)^* .$$ 
Here we are only concerned with $K$-types which have trivial $\SO(V_-)$-representation. 
Therefore $\nu = r \times q$ for some integer $r$. This proves the first assertion of the proposition.

We now consider the decomposition of 
$$[\wedge^{\bullet} \p ]^{\SL (V_-)} = \oplus_{n=0}^p S_{n \times q} (V_+) \otimes (\wedge^q V_-)^n$$
into irreducibles. By (Poincar\'e) duality it is enough to consider $S_{n \times q} (V_+) \otimes (\wedge^q V_-)^n$ with $n \leq p/2$. It then follows from Proposition \ref{Lit} that the multiplicity of $S_{[r \times q]} (V_+)$ in $S_{n \times q} (V_+)$ equals
$$\sum_{\xi} \dim \mathrm{Hom}_{\GL (V_+)} (S_{n \times q} (V_+) , S_{r \times q} (V_+) \otimes S_{\xi} (V_+)),$$
where the sum is over all nonnegative integer partitions $\xi$ with rows of even length. Each term of the sum above is 
a Littlewood-Richardson coefficient. The Littlewood-Richardson rule states that $\dim \mathrm{Hom}_{\GL (V_+)} (S_{n \times q} (V_+) , S_{r \times q} (V_+) \otimes S_{\xi} (V_+))$ equals
the number of Littlewood-Richardson tableaux of shape  $(n\times q) / (r\times q) = (n-r) \times q$ and of weight $\xi$, see e.g. \cite{Fulton}. 
The point is that the shape $(n\times q) / (r\times q)$ is an $n-r$ by $q$ rectangle and it is immediate that there is only one semistandard filling of an $n$ by $q$ rectangle that satifies the reverse lattice word condition, see \cite{Fulton}, \S 5.2, page 63, (the first row must be filled with ones, the second with twos etc.). 
We conclude that 
\begin{multline*}
\dim \mathrm{Hom}_{\GL (V_+)} (S_{n \times q} (V_+) , S_{r \times q} (V_+) \otimes S_{\xi} (V_+)) \\
= \dim \mathrm{Hom}_{\GL (V_+)} (S_{(n-r) \times q} (V_+) , S_{\xi} (V_+)).
\end{multline*}
The multiplicity of $S_{[r \times q]} (V_+)$ in $S_{n \times q} (V_+)$ therefore equals $0$ is $q$ is odd and $1$ if $q$ is even.
Since in the last case
$$e_q^{n-r} (V_r) \cong S_{[r \times q]} (V_+) \otimes (\wedge^q V_-)^n \subset [\wedge^{nq} V_+ \otimes V_-^* ]^{\SL (V_-) } = S_{n \times q} (V^+) \otimes (\wedge^q V_-)^n,$$
this concludes the proof.
\end{proof}

\medskip
\noindent
{\it Remark.} Proposition \ref{PropLit} implies the decomposition \eqref{subringSC} of the Introduction. 

\medskip

We conclude the section by showing that the Vogan-Zuckerman types $\mu_r$ are the only $K$-types to give
small degree cohomology.

\begin{prop} \label{P:cohomrep}
Consider a cohomological module $A_{\mathfrak{q}} (\lambda )$. Suppose that
$R= \dim (\mathfrak{u} \cap \mathfrak{p})$ is strictly less than both $p+q-3$ and $pq/4$. Then: 
either $L = C \times \SO_0 (p-2n , q)$ with $C \subset K$ and $R=nq$ or $L=C \times \SO_0 (p, q-2n)$
with $C \subset K$ and $R=np$.
\end{prop}
\begin{proof} Suppose by contradiction that $L$ contains as a direct factor the group $\SO_0 (p-2a, q-2b)$ for some positive $a$ and $b$. Then $R \geq ab + b(p-2a) + a (q-2b)$. And writing:
\begin{multline*}
ab + b(p-2a) + a (q-2b) - (p+q-3) \\  = (a-1) (b-1) + (b-1) (p-2a-1) + (a-1)(q-2b-1) 
\end{multline*}
we conclude that $R \geq p+q-3$ except perhaps if $p=2a$ or $q=2b$. But in that last case $R \geq 
pq/4$.
\end{proof}

\section{Cohomological Arthur packets} \label{sec:AP}

For archimedean $v$, local Arthur packets $\small\prod (\Psi)$ should coincide with the packets
constructed by Adams, Barbasch and Vogan \cite{ABV}. Unfortunately this is still unproved. We may nevertheless build upon their results to obtain a conjectural description of all the real Arthur packets 
which contain a cohomological representation. 

\subsection{Notations} Let $p$ and $q$ two non-negative integers with $p+q =m$. We set $\ell = [m/2]$ and $N=2\ell$. 
In this section $G= \SO (p,q)$, $K$ is a maximal compact subgroup of $G$ and $K_0 = \SO (p) \times \SO (q)$. 
We let $\g_0$ the real Lie algebra of $G$ and $\g_0 = \k_0 \oplus \p_0$ be the Cartan decomposition associated to the choice of the maximal
compact subgroup $K$. We denote by $\theta$ the corresponding Cartan involution. If $\mathfrak{l}_0$ is a real Lie algebra we denote by $\mathfrak{l}$ its complexification $\mathfrak{l} = \mathfrak{l}_0 \otimes \C$.   

We finally let $W_{\R}$ be the Weil group of $\R$ 

\subsection{The Adams-Johnson packets} Let $\mathfrak{l} \supset \mathfrak{g}$ be the 
Levi component of a $\theta$-stable parabolic subalgebra
$\mathfrak{q}$ of $\g$. Then $\mathfrak{l}$ is defined over $\R$ and we let $L$ be the 
corresponding connected subgroup of $G$. Let $T \subset G$ be a $\theta$-stable Cartan subgroup of $G$ such that $T \cap L$ is a Cartan subgroup of $L$. 
We fix $\lambda : \mathfrak{l} \rightarrow \C$ an 
admissible one-dimensional representation that is zero on the orthogonal complement of $\mathfrak{t}\cap \mathfrak{l}$ in $\mathfrak{t}$.  We will identify $\lambda$ with its highest weight in $\mathfrak{t}^*$ that is its restriction to $\mathfrak{t}$. 

Adams and Johnson have studied the particular family --- to be called {\it Adams-Johnson parameters} ---
of local Arthur parameters 
$$\Psi : W_{\R} \times \SL_2 (\C) \rightarrow {}^L \- G$$
such that 
\begin{enumerate}
\item $\Psi$ factors through ${}^L \- L$, that is  $\Psi : W_{\R} \times \SL_2 (\C) \stackrel{\Psi_L}{\rightarrow} {}^L \- L \rightarrow {}^L \- G$ where the last map is the canonical extension \cite[Proposition 1.3.5]{Shelstad} of the injection $L^{\vee} \subset G^{\vee}$, and
\item $\varphi_{\Psi_L}$ is the $L$-parameter of a unitary character of $L$ whose differential is 
$\lambda$. 
\end{enumerate}
The restriction of 
the parameter $\Psi$ to $\SL_2 (\C)$ therefore maps $\left( \begin{smallmatrix} 1 & 1 \\ 0 & 1 \end{smallmatrix} \right)$ to a principal unipotent element in $L^{\vee} \subset G^{\vee}$.

\subsection{} The packet $\small\prod_{\rm AJ} (\Psi) $ constructed by Adams and Johnson takes the following form. 
Let $W(\g , \mathfrak{t})^{\theta}$ (resp. $W  (\mathfrak{l} , \mathfrak{t})^{\theta}$) be those elements of 
the Weyl group of $\g$ which commute with $\theta$. And let $W(G,T)$ be the real Weyl group of $G$.
The representation in $\small\prod_{\rm AJ}  (\Psi )$ are parametrized by the double cosets
$$S = W  (\mathfrak{l} , \mathfrak{t})^{\theta} \backslash W(\g , \mathfrak{t})^{\theta} / W( G , T).$$
For any $w \in S$, the Lie subalgebra $\mathfrak{l}_{w} = w \mathfrak{l} w^{-1}$
is still defined over $\R$, and is the Levi subalgebra of the $\theta$-stable parabolic subalgebra
$\mathfrak{q}_w = w \mathfrak{q} w^{-1}$. Let $L_w \subset G$ be the corresponding connected group.
The representations in $\small\prod_{\rm AJ}  (\Psi)$ are the irreducible unitary representations
of $G$ whose underlying $(\g , K_0)$-modules are the Vogan-Zuckerman modules
$A_{\mathfrak{q}_w} (w \lambda)$. 

We call {\it Adams-Johnson packets} the packets associated to Adams-Johnson parameters as above. 

The following {\it folklore} conjecture does not seem to follow from any published work.

\begin{conj} \label{conj}
\begin{enumerate}
\item For any Adams-Johnson parameter $\Psi$ one has $\small\prod (\Psi) = \small\prod_{\rm AJ} (\Psi )$. 
In otherwords, the Arthur packet associated to 
$\Psi$ coincides with the Adams-Johnson packet associated to $\Psi$. 
\item The only Arthur packets that contain a cohomological representation of $G$ are Adams-Johnson packets.
\end{enumerate}
\end{conj}

\medskip
\noindent
{\it Remark.} Adams-Johnson packets are a special case of 
the packets define in \cite{ABV}. In otherwords: calling {\it ABV packets} the packets defined
in \cite{ABV}, any Adams-Johnson packet is an ABV packet defined by a parameter $\Psi$
of a particular form. According to J. Adams (private communication) part 2 of Conjecture \ref{conj}
with ``ABV packets'' in place of ``Arthur packets'' is certainly true. It would therefore be enough to prove 
the natural extension of part 1 of Conjecture \ref{conj} to any ABV packet. In chapter 26 of \cite{ABV}, Adams, Barbasch and Vogan prove that ABV packets are compatible with endoscopic lifting.\footnote{In particular, it would be enough to check Conjecture \ref{conj} only for quasi-split groups.} Part 1 of Conjecture \ref{conj} amounts to a similar statement for twisted endoscopy. This seems still open 
but is perhaps not out of reach.

\subsection{} \label{6.5} We now give a more precise description of the possible Adams-Johnson parameters:
The Levi subgroups $L$ 
corresponding to a $\theta$-stable parabolic subalgebra of $\g$ are of the forms
\begin{equation} \label{Levi}
\mathrm{U} (p_1 , q_1) \times \ldots \times \mathrm{U} (p_{r} , q_{r}) \times \SO (p_0 , q_0)
\end{equation}
with $p_0 + 2 \sum_j p_j = p$ and $q_0 + 2 \sum_j q_j = q$. We let $m_j = p_j+q_j$ ($j=0 , \ldots , r$). Then $m_0=p_0+q_0$ has the same parity as $m$; we set $\ell_0 = [m_0 /2]$ and $N_0 = 2 \ell_0$. Here and after we assume that $p_0q_0 \neq 0$. The parameter $\Psi_L$ corresponding to such an $L$ maps $\left( \begin{smallmatrix} 1 & 1 \\ 0 & 1 \end{smallmatrix} \right)$ to a principal unipotent element in each factor of $L^{\vee} \subset G^{\vee}$. The parameter $\Psi_L$ therefore contains a $\SL_2 (\C)$ factor of the maximal dimension
in each factor of $L^{\vee}$. These factors consist of $\GL (m_j)$, $j=1 , \ldots , r$, and 
the group $G^{\vee}_{m_0} = \SO (N_0 , \C)$ or $\Sp (N_0 , \C)$ (according to the parity of $m$). In any case the biggest possible $\SL_2 (\C)$ representation is $R_{m_0-1}$.\footnote{Note that $R_a$ is symplectic if $a$ is
even and orthogonal if $a$ is odd.} We conclude that the restriction of $\Psi$ to $\C^* \times \SL_2 (\C)$
decomposes as:
\begin{equation} \label{AJparam1}
(\mu_1 \otimes R_{m_1} \oplus \mu_1^{-1} \otimes R_{m_1} ) \oplus 
\ldots \oplus (\mu_r \otimes R_{m_r} \oplus \mu_r^{-1} \otimes R_{m_r} ) \oplus \Psi_0
\end{equation}
where the $\mu_j$ are unitary characters of $W_{\R}$ and 
\begin{equation} \label{AJparam2}
\Psi_0 = \left\{
\begin{array}{ll}
\chi \otimes R_{m_0 -1} & \mbox{if } m_0 \mbox{ is odd}, \\
\chi \otimes R_{m_0 -1} \oplus \chi ' \otimes R_1 &  \mbox{if } m_0 \mbox{ is even}. 
\end{array} \right.
\end{equation}
Here  $\chi$ and $\chi'$ are quadratic characters of $W_{\R}$.

\subsection{} \label{6.6} Now consider a cohomological representation $\pi$ of $G$. Its underlying $(\g , K_0)$-module is a Vogan-Zuckerman module $A_{\mathfrak{q}} (\lambda )$. It follows from their construction 
that we may choose $\mathfrak{q}$ so that the group $L$ has no compact (non-abelian) simple factors; in that case $\mathfrak{q}$ and the unitary character of $L$ whose differential is 
$\lambda$ are uniquely determined up to conjugation by $K_0$. The group $L$ is therefore as in \eqref{Levi} with either $p_jq_j \neq 0$ or $(p_j,q_j) = (1,0)$ or $(0,1)$. So that:
$$L^{\vee} \cong \GL(1)^s \times \GL(m_1) \times \ldots \times \GL (m_t ) \times G^{\vee}_{N_0}$$
where each $m_j$, $j=1 , \ldots , t$, is $>1$, $r=s+t$ and $m_{t+1} = \ldots = m_{r} = 1$. We let $\Psi$ be the Adams-Johnson parameter associated to $L$ as in the preceeding paragraph.

\begin{lem} \label{Lem:S}
\begin{enumerate}
\item The cohomological representation $\pi$ belongs to $\small\prod_{\rm AJ} (\Psi)$. 
\item If $\pi \in \small\prod_{\rm AJ} (\Psi ')$ for some Adams-Johnson parameter $\Psi '$. Then $\Psi'$
contains 
$$ (\mu_1 \otimes R_{m_1} \oplus \mu_1^{-1} \otimes R_{m_1} ) \oplus 
\ldots \oplus (\mu_t \otimes R_{m_t} \oplus \mu_t^{-1} \otimes R_{m_t} ) \oplus \Psi_0$$
as a direct factor.
\end{enumerate}
\end{lem} 
\begin{proof} 1. The representation $\pi$ belongs to $\small\prod_{\rm AJ} (\Psi)$ by definition of
the Adams-Johnson packets. Note that the Langlands parametrization of Vogan-Zuckerman modules is given
by \cite[Theorem 6.16]{VZ}. In the case of unitary groups this is detailed in \cite{BC} where the
parametrization is moreover related to Arthur parameters. Everything works similarly for orthogonal groups. 

2. If $\pi  \in \small\prod_{\rm AJ} (\Psi ')$ for some Adams-Johnson parameter $\Psi '$ we let 
$L'$ be the Levi subgroup attached to the parameter $\Psi'$. The 
underlying $(\g , K_0)$-module of $\pi$ is isomorphic to some Vogan-Zuckerman
module associated to a $\theta$-stable parabolic algebra whose associated Levi subgroup
is $L_{w}'$ for some $w \in S$. It first follows from the Vogan-Zuckerman construction that $L$ and 
$L_{w} '$ can only differ by compact factors. 
But it follows from \cite[Lemma 2.5]{AJ} that all the $L_w'$ ($w \in S$) are inner forms of each other 
so that the dual group of $L_w'$ may be identified to $L' \- {}^{\vee}$ which therefore contains 
$\GL(m_1) \times \ldots \times \GL (m_t ) \times G^{\vee}_{N_0}$ as  direct factor. And Lemma \ref{Lem:S}
follows from \S \ref{6.5}.
\end{proof}

In particular Lemma \ref{Lem:S} implies that if $\pi$ is a cohomological representation of $\SO (p,q)$ 
associated to a Levi subgroup $L=\SO(p-2r , q) \times \mathrm{U}(1)^r$ with $p>2r$ and $m-1>3r$. Then
if $\pi$ is contained in a (local) Adams-Johnson packet $\small\prod_{\rm AJ} (\Psi)$, the (local) parameter $\Psi$ contains a factor $\eta \boxtimes R_{m-2r-1}$ where $\eta$ is a quadratic character and $3(m-2r-1)> m-1$. Theorem \ref{Thm1} --- see especially Remark (2) after its proof --- therefore implies the following:

\begin{prop} \label{Prop:conj}
Assume Conjecture \ref{conj}. Let $\pi \in \mathcal{A}^c (\SO (V))$ and let $v_0$ be an infinite place of $F$
such that $\SO (V)(F_{v_0}) \cong \SO (p,q)$. Assume that $\pi_{v_0}$ is a cohomological representation of 
$\SO (p,q)$ associated to a Levi subgroup $L=\SO(p-2r , q) \times \mathrm{U}(1)^r$ with $p>2r$
and $m-1> 3r$. 
Then: there exists an automorphic character $\chi$ such that $\pi \otimes \chi$ is in the image of the cuspidal $\psi$-theta correspondence from a smaller group associated to a symplectic space of dimension $2r$.
\end{prop}

Unfortunately Conjecture \ref{conj} seems to be open. Instead we prove the
weaker (but unconditional) following result.

\begin{prop} \label{Prop:appendix}
Let $\pi \in \mathcal{A}^c (\SO (V))$ and let $v_0$ be an infinite place of $F$
such that $\SO (V)(F_{v_0}) \cong \SO (p,q)$. Assume that $\pi_{v_0}$ is a cohomological representation of 
$\SO (p,q)$ associated to a Levi subgroup $L=\SO(p-2r , q) \times \mathrm{U}(1)^r$ with $p>2r$
and $m-1> 3r$. 
Then:  the (global) Arthur parameter $\Psi$ of $\pi$ is highly non-temptered, i.e. contains a factor $\eta \boxtimes R_a$ where $\eta$ is a quadratic character and $3a > m-1$.
\end{prop} 

Since the reader may want to believe Conjecture \ref{conj} we postpone the proof of Proposition \ref{Prop:appendix} to  Appendix (see Proposition \ref{Prop:appendix2} and the remark following it). It relies on the theory of exponents. 




As we will explain later in Remark \ref{5.11}, a cohomological representation as in Proposition \ref{Prop:appendix}
does not occur in Howe's theta correspondence from a symplectic group smaller than $\Sp_{2r}$. 
In conjonction with Theorem \ref{Thm1} we obtain:

\begin{cor} \label{Cor:6.6} Let $\pi \in \mathcal{A}^c (\SO (V))$ and let $v_0$ be an infinite place of $F$
such that $\SO (V)(F_{v_0}) \cong \SO (p,q)$. Assume that $\pi_{v_0}$ is a cohomological representation of 
$\SO (p,q)$ associated to a Levi subgroup $L=\SO(p-2r , q) \times \mathrm{U}(1)^r$ with $p>2r$
and $m-1> 3r$. 
Then: there exists an automorphic character $\chi$ such that $\pi \otimes \chi$ is in the image of the cuspidal $\psi$-theta correspondence from a smaller group associated to a symplectic space of dimension $2r$.
\end{cor}

\subsection{}{\it Remark.} We believe that the hypothesis $m-1>3r$ is optimal. This is indeed the case if
$p=3$, $q=1$ and $r=1$, see Proposition \ref{prop:dim3}. 

\medskip

We now provide support for the above remark.

\subsection{} \label{6.13} In the remaining part of this section we assume that Conjecture \ref{conj} --- which relates the Adams-Johnson's pakets and the Arthur's packet containing at least one cohomological representations --- holds. In case $m$ is odd we moreover assume that the representations of the metaplectic groups are classified with Arthur's like parameters (see below). Assuming that we will prove:

\medskip
\noindent
{\it If $3r \geq m-1$ there exist cohomological cuspidal representations  which are not in the image of the $\theta$ correspondence from $\Mp(2r)$ or $\Sp(2r)$ even up to a twist by an automorphic character.}

\medskip

\subsection{} We first assume that $m=N+1$ is odd.

A particular case of Arthur's work is the classification of square integrable representation of $\SL(2,F)$ using $\GL(3,F)$; this can be also covered by the known Gelbart-Jacquet correspondance between $\GL(2)$ and $\GL(3)$. We therefore take it for granted.

We denote by $F_{2}$ the localization of the field $F$ at the places of residual characteristic $2$. 

Let $\tau_{2}$ be a cuspidal irreducible representation of $\GL(3,F_{2})$ autodual which comes from a representation  of $\SL (2,F_{2})$ or, in other terms, whose $L$-parameter factorizes through $\SO (3,{\mathbb C})$. We denote by $\tilde{\tau}_{2}$ the corresponding representation of $\SL(2,F_{2})$. We fix $\tilde{\tau}$ a cuspidal irreducible representation of $\SL(2,F)$ whose $F_{2}$ component is $\tilde{\tau}_{2}$ and which is a discrete series at the archimedean places. We go back to $\GL(3,F)$ denoting by $\tau$ the automorphic representation corresponding to $\tilde{\tau}$; because of the condition on the $F_{2}$-component $\tau$ is necessarily cuspidal.

For each $i \in \{1,\ldots ,3r-m+1\}$, we also fix a cuspidal irreducible representation $\rho_{i}$ of $\GL(2,F)$. We assume that each $\rho_i$ is of symplectic type, i.e. its local parameter is symplectic. In other words each $\rho_i$ is coming --- by the Langlands-Arthur functoriality --- from $\SO(3,F)$, equivalently $L(s, \rho_{i}, \mathrm{Sym}^2)$ has a pole at $s=1$. We moreover assume that at each archimedean place $v | \infty$, each representation $\rho_{i, v}$ belongs to the discrete series. We consider the Arthur 
parameter 
\begin{equation} \label{6*}
\Psi = \left( \boxplus_{i=1}^{3r-m+1}\rho_{i}\boxtimes R_{1} \right) \boxplus \tau\boxtimes R_{m-1-2r}.
\end{equation}
This is the Arthur parameter of a packet of representations of $G=\SO(V)$. 

\subsection{} We now look more precisely at $\Psi$ at an archimedean place $v| \infty$. As a morphism of $W_{{\mathbb R}}\times \SL_2({\mathbb C})$ in $\Sp(N , \C) $ it is the sum:
$$
\Psi_v = \left( \boxplus_{i=1}^{3r-m+1} \phi_{\rho_{i,v}} \boxtimes R_{1} \right) \boxplus \phi_{\tilde{\tau}_{v}}\boxtimes R_{m-2r-1}\boxplus 1\boxtimes R_{m-2r-1},
$$
where for $i=1,\ldots ,3r-m+1$, $\phi_{\rho_{i,v}}$ is the Langlands parameter of the discrete series $\rho_{i,v}$ and $\phi_{\tilde{\tau}_{v}}$ is the Langlands parameter of $\tilde{\tau}_{v}$. All such local parameters may be obtained by the above construction. By a suitable choice of data we can therefore assume that $\Psi_v$ coincides with an Adams-Johnson parameter. In particular the representations attached to it are all cohomological for some fixed system of coefficients. Let 
$\lambda_v$ be the corresponding highest weight.

We may moreover assume the Adams-Johnson parameter $\Psi_v$ is associated to the Levi 
subgroup $L= \SO (m)$ if $v \neq v_0$ and $L=\mathrm{U}(1)^{3r-m+1}\times \mathrm{U}(m-2r-1)\times \SO(p-2r,q)$ if $v=v_0$.\footnote{Recall that $G(F_{v_0}) = \SO(p,q)$ and that $G(F_v) = \SO (m)$ is $v\neq v_0$.} We also ask that $\lambda_v = 0$ if $v \neq v_0$, and fix $\lambda_v = \lambda$ if $v =v_0$. It follows from (the proof of)
Lemma \ref{Lem:S}(2) that the trivial representation $\pi_v$ of $G(F_v)$ is contained in $\small\prod_{\rm AJ} (\Psi_v)$ if $v \neq v_0$ and that the cohomological  representation $\pi_{v_0}$ of $G(F_{v_0} ) = \SO (p,q)$ whose underlying $(\g, K_0)$-module is the Vogan-Zuckerman module $A_{{\mathfrak{q}}}(\lambda)$, 
with ${\mathfrak {q}}$ a $\theta$-stable parabolic subalgebra with real Levi component isomorphic to 
$\mathrm{U}(1)^{r} \times \SO(p-2r,q)$, is contained in $\small\prod_{\rm AJ} (\Psi_{v_0})$. We fix these
local components $\pi_v$. 

\subsection{} The multiplicity formula to construct a global square integrable representation of $G=\SO(V)$ from the local components is still the subject of work in progress of Arthur; but we can anticipate that we have enough freedom at the finite places to construct a square integrable representation $\pi$ in the global Arthur's packet and with local component at the archimedean places, the component we have fixed.

We want to show that the representation $\pi$ is certainly not obtained via theta correspondence from a cuspidal representation of a metaplectic group $\Mp(2n)$ with $2n\leq 2r$. To do that we continue to anticipate some results: Here we anticipate that {\it the square integrable representations of the metaplectic group can also be classified as those of the symplectic group but using $\Sp(2n,{\mathbb C})$ as dual group}; after work by Adams, Adams-Barbash, Renard, Howard this is work in progress by Wen Wei Li. 

To prove that $\pi$ is not a theta lift we can now argue by contradiction: Let $\sigma$ be a cuspidal irreducible representation of $\Mp(2n)$ (with $2n\leq 2r<p$) such that $\pi$ is a theta lift of $\sigma$. Write $\boxplus_{i=1}^v \sigma_{i}\boxtimes R_{n_{i}}$ the Arthur like parameter attached to $\sigma$. To simplify matter assume that $V$ is an orthogonal form of discriminant 1 at each place (otherwise we would have to twist by the quadratic character which corresponds by class field theory to this discriminant). Consider the parameter:
\begin{equation} \label{6**}
\oplus_{i\in [1,v]}\sigma_{i}\boxtimes R_{n_{i}} \oplus 1\boxtimes R_{m-1-2n}.
\end{equation}
Here we use the fact that $m-1-2n\geq m-1-2r \geq p-2r\geq 1$.

The local theta correspondence is known at each place where $\sigma$ and $\pi$ are both unramified, see \cite{MoeglinTheta}. This implies that at every unramified place $\pi_{v}$ is necessarily the unramified representation in the local Arthur packet associated to the parameter \eqref{6**}. But by definition, $\pi_{v}$ is also in the local packet associate to \eqref{6*}; this implies that \eqref{6*} and \eqref{6**} define automorphic (isobaric) representations of $\GL(2n)$ which are isomorphic almost everywhere. These automorphic representation are therefore isomorphic and \eqref{6*} must coincide with \eqref{6**}. This is in contradiction with the fact that there is no factor $\eta \boxtimes R_{m-1-2n}$ --
for some automorphic quadratic character $\eta$ of $\GL(1)$ --- in \eqref{6*}.

\subsection{} We now assume that $m=N$ is even. We moreover assume that $p-2r>1$. 

We then do the same: First construct $\tau$ as above. For each $i\in \{1, \ldots , 3r-m+1 \}$ we fix a cuspidal irreducible representation of $\GL(2,F)$ of orthogonal type this means that $L(s, \rho_{i},\wedge^2)$ has a pole at $s=1$; we can impose any discrete series at the archimedean places we want. The Arthur parameter we look at is now:
\begin{equation} \label{6***}
\boxplus_{i=1}^{3r-m+1} \rho_{i} \boxtimes R_{1} \boxplus \tau\boxtimes R_{m-1-2r} \boxplus \eta\boxtimes R_{1},
\end{equation}
where $\eta$ is a suitable automorphic quadratic character of $\GL(1)$ is such a way that this parameter is relevant for the quasisplit form of $\SO(V)$, see \cite{ArthurBook}. 

Now we argue as above to construct a global representation $\pi$ of $\SO(V)$ in this packet which is as we want at the archimedean places. We construct \eqref{6**} as above and here we have not to anticipate more results than those announced by Arthur. But to conclude we have to make sure that $m-1-2n>1$ because there is a factor $\eta \boxtimes R_{1}$ in \eqref{6***}.
Nevertheless: as we have hypothesised that
$$
m-1-2n\geq m-1-2r\geq p-2r>1,
$$
we obtain a contradiction which proves that $\pi$ is not a $\theta$-lift of a cuspidal representation of $\Sp(2n)$ with $n\leq r$.

\section{The classes of Kudla-Millson and Funke-Millson}

In this section we will introduce the $(\mathfrak{so}(p,q) ,K)$-cohomology classes (with $K= \mathrm{SO}(p) \times \mathrm{SO}(q)$) of Kudla-Millson and Funke-Millson, find explicit values for
these classes  on the highest weight vectors of the Vogan-Zuckerman special K-types $V(n,\lambda)$, see Proposition \ref{harmonicvalue} and Proposition \ref{finalformula}, and from those values and a result of Howe deduce the key result that the translates of the above classes
under the universal enveloping algebra of the symplectic group span $\mathrm{Hom}_K(V(n,\lambda), \mathcal{P}(V^n))$ where $\mathcal{P}(V^n)$
is the polynomial Fock space, see Theorem \ref{Thm:5.8}. Our computations  will give an new derivation of some of the formulas in \cite{KashiwaraVergne}. 

\subsection{Notations} Let $V$ be a real quadratic space of dimension $m$ and signature $(p,q)$  two non-negative integers with $p+q =m$. Set $G=\SO_0 (V)$ and let $K$ be a maximal compact subgroup. We write $\mathcal{S} (V^n)$ for the space of (complex-valued) Schwartz functions on $V^n$. Let $\Mp_{2n} (\R)$ denote the metaplectic double cover of $\Sp_{2n} (\R)$ if $m$ is odd or simply $\Sp_{2n} (\R)$ if $m$ is even.  We recall that $\Mp_{2n} (\R)$ acts in the Schwartz space $\mathcal{S} (V^n)$ via the oscillator representation $\omega = \omega_{\psi}$ associated to the additive character $\psi : x \mapsto \exp (2i\pi x)$ of $\R$. Let $K' \subset \Mp_{2n} (\R)$ be the inverse image of the maximal compact subgroup 
\begin{equation} \label{embeddingofunitary}
\left\{ \left(
\begin{array}{cc}
a & b \\
-b & a 
\end{array} \right) \; : \; a+ib \in \mathrm{U}(n) \right\}
\end{equation}
of $\Sp_{2n} (\R)$. When $m$ is odd we let $(\det)^{\frac12}$ denote a character of $K'$ whose square
descends to the determinant character of $\mathrm{U}(n)$. We let $\C_{\frac{m}{2}}$ the corresponding 
one-dimensional representation of $K'$.

As in the preceeding section we let $\g$ be the complexified Lie algebra of $G$ and $\g = \p \oplus \k$ its Cartan decomposition, where $\mathrm{Lie} (K)\otimes \C= \k$. We let $\g '$ be the complexified
Lie algebra of $\Sp_{2n} (\R )$.

In this section we make the assumption that 
\begin{equation} \label{eq:assump}
2n <p.
\end{equation}

\subsection{Three consequences of \eqref{eq:assump}}
Note that $\mathrm{rank}(G) = [\frac{m}{2}]$ and $\mathrm{rank}(\mathrm{SO}(V_+)) = [\frac{p}{2}]$. It then follows from \eqref{eq:assump} that we have:
\begin{enumerate}
\item If $m= p+q$ is even then $n < \mathrm{rank}(G) -1.$
\item If $m= p+q$ is odd then $n < \mathrm{rank}(G)$.
\item $n \leq \mathrm{rank}(\mathrm{SO}(V_+))$.
\end{enumerate}

As we will explain in Subsection \ref{KVrel}, item (3) will imply we are in the first family for the two families of formulas in \cite{KashiwaraVergne} concerning the action
of $\mathrm{GL}(n,\C) \times \mathrm{O}(V_+)$  on the (pluri)harmonic polynomials on $V_+^n$. In what follows we will refer to the formulas in this first family as belonging to ``Case (1)'' of \cite{KashiwaraVergne}. 
\medskip
\noindent

\subsection{The Fock model} \label{Fock} In what follows we will use the Fock model of the oscillator representation. We let 
$W$ be a real vector space of dimension $2n$ equipped with a non-degenerate skew symmetric form 
$\langle , \rangle$ such that $\Sp_{2n} (\R ) = \Sp (W)$. Fix a positive definite complex structure $J$ on $W$. Let $\{x_1 , \ldots , x_n , y_1 , \ldots , y_n\}$ be a symplectic basis for $W$ such that
$$Jx_i = y_i \quad \mbox{ and } \quad Jy_i = -x_i \quad \mbox{ for } 1 \leq i \leq n.$$
We may decompose $W\otimes \C$ according to 
$$W\otimes \C = W' + W''$$
where $W'$ (resp. $W''$) is the $+i$ (resp. $-i$)-eigenspace of $J$. We define complex bases 
$\{ \tilde{w}_1 ' \ldots , \tilde{w}_n ' \}$ and $\{\tilde{w}_1 '' , \ldots , \tilde{w}_n '' \}$ for $W'$ and $W''$ respectively by 
$$\tilde{w}_j ' = x_j -iy_j \quad \mbox{ and } \quad \tilde{w}_j '' =x_j +i y_j \quad \mbox{ for } 1 \leq j \leq n.$$ 

The following  will be important in the rest of the paper.  We leave the proof to the reader.
\begin{lem}\label{holomorphicsubspace}
The subgroup $\mathrm{U}(n)$ of $\Sp_{2n}(\R)$ embedded according to equation \eqref{embeddingofunitary} acts on $W''$ by the standard representation
and on $W'$ by the dual of the standard representation.
\end{lem}

We finally denote by $(,)$ the non-degenerate quadratic
form on $V$. By a slight abuse of notations we will also denote by $\langle , \rangle$ the non-degenerate skew-symmetric form $(,) \otimes \langle , \rangle$ on $V\otimes W$.

Now let $v_{\alpha}$, $\alpha = 1 ,\ldots , p$, $v_{\mu}$, $\mu = p+1 , \ldots , q$, be an orthogonal basis of $V$ such that $(v_{\alpha} , v_{\alpha})=1$ and $(v_{\mu} , v_{\mu})=-1$. We denote by $V_+$ (resp. $V_-$) the span of 
$\{ v_{\alpha} \; : \; 1¬¨‚?†\leq \alpha \leq p \}$ (resp. $\{v_{\mu} \; : \; p+1 \leq \mu \leq m \}$). 

The $(\mathfrak{sp}_{2mn} , \widetilde{\mathrm{U}}_{mn)})$-module associated to the oscillator representation is made explicit in the Fock realization, see \cite{HoweT}. Using it, one sees that the space of $\widetilde{\mathrm{U}}_{mn}$-finite vectors is the symmetric algebra $\mathrm{Sym}(V_+ \otimes W'' + V_- \otimes W' )$. We introduce linear 
functionals
$$\{z_{\alpha , j} , z_{\mu ,j} \; : \; 1 \leq \alpha \leq p, \quad p+1 \leq \mu \leq m, \quad 1 \leq j \leq n \}$$
on $V_+ \otimes W' + V_- \otimes W''$ by the formulas 
$$z_{\alpha , j } (v \otimes w) = \langle v \otimes w , v_{\alpha} \otimes \tilde{w}_j '' \rangle, \quad 1 \leq \alpha \leq p, \quad 1 \leq j \leq n ,$$
$$z_{\mu , j } (v \otimes w) = \langle v \otimes w , v_{\mu} \otimes \tilde{w}_j ' \rangle, \quad p+1 \leq \mu \leq m, \quad 1 \leq j \leq n .$$
Hence we have
$$x_j = \sum_{\alpha=1}^p z_{\alpha,j} v_{\alpha} + \sum_{\mu=p+1}^{p+q} z_{\mu,j} v_{\mu}.$$

We use these coordinates to identify 
 the space $\mathrm{Sym}(  V_+ \otimes W''+ V_- \otimes W')$ (which we identify with the space of polynomial functions on $V_+ \otimes W' + V_- \otimes W''$ using the symplectic form on $V \otimes W$)  with the space $\mathcal{P} = \mathcal{P}(\C^{mn})$ of polynomials in $mn$ variables $z_{1,j} , \ldots , z_{m,j}$ ($j=1 , \ldots , n$). We will use $\mathcal{P}_+$ to denote the polynomials in the ``positive'' variables $z_{\alpha,j}, 1 \leq \alpha \leq p, 1 \leq j \leq n$.  We will use the symbol $\mathcal{P}^{(\ell)}$ to denote the subspace of polynomials of degree $\ell$ and similarly for $\mathcal{P}_+^{(\ell)}$. It will be important to record the structure of $\mathcal{P}(\C^{mn})$ as a $K' \times K$-module.
By Lemma \ref{holomorphicsubspace} we know that the action of $\mathrm{GL}(n,\C)$ on $W''$ is equivalent to the standard one.  
We will accordingly identify $\C^n$ with $W''$ (and ($\C^n)^*$ with $W'$)  as modules for $\mathrm{GL}(n,\C)$. We can then replace $W''$   with $\C^n$ and $W'$ by $(\C^n)^*$ (since we are only concerned here with the $\mathrm{GL}(n,\C)$-module structure of $W'$ and $W''$). We will  use
$e_1,\cdots,e_n$ to denote the standard basis of $\C^n$, accordingly under the above identification $e_j$ corresponds to $\tilde{w}''_j$ for $1 \leq j \leq n$. 
\begin{lem}\label{Fockspace}
We have an isomorphism  of $K' \times K$-modules
$$\mathcal{P}(\C^{pn}\times {(\C^*)}^{qn}) \cong \mathrm{Sym}( V_+ \otimes W''+ V_- \otimes W') =  \mathrm{Sym}(V_+ \otimes W'')\otimes \mathrm{Sym}(V_- \otimes W').$$
\end{lem}
 Thus the variables $z_{\alpha,j},1 \leq \alpha \leq p,1 \leq j \leq n,$ transform (in $j$) according to the standard representation of $\mathrm{GL}(n,\C)$ and the variables $z_{\mu,j},p+1 \leq \alpha \leq p+q,1 \leq j \leq n,$ transform in $j$ according to the dual of the standard representation of $\mathrm{GL}(n,\C)$. We will henceforth use the symbol $\mathcal{P}$ to denote the above algebra of polynomials {\it equipped with the graded $\mathrm{GL}(n,\C)$-module structure just described}.

Numbering the $nm$ variables as $z_1 , \ldots , z_{mn}$ we have
\begin{eqnarray*}
\omega (\mathfrak{sp}_{\C}) & = & \mathfrak{sp}^{(1,1)} \oplus \mathfrak{sp}^{(2,0)} \oplus \mathfrak{sp}^{(0,2)} , \\
\mathfrak{sp}^{(1,1)} & = & \mathrm{span} \left\{ \frac12 \left (z_i \frac{\partial}{\partial z_j} + \frac{\partial}{\partial z_j} z_i \right) \right\} , \\
\mathfrak{sp}^{(2,0)} & = &\mathrm{span}\{ z_i z_j \} , \\
\mathfrak{sp}^{(0,2)} & = & \mathrm{span} \left\{ \frac{\partial^2}{\partial z_i \partial z_j} \right\} .
\end{eqnarray*}
Here $\mathfrak{sp}_{\C}$ is the complexification of $\mathfrak{sp}_{2mn}$. The indices $i$ and $j$ vary 
from $1$ to $nm$. Note that while $\mathfrak{sp}_{2mn}$ is a real Lie algebra, the $\mathfrak{sp}^{(a,b)}$ are complex subalgebras of $\omega (\mathfrak{sp}_{\C})$. In the Cartan decomposition 
\begin{equation} \label{Cartan}
\mathfrak{sp}_{2mn} = \mathfrak{u}_{mn} \oplus \mathfrak{q}
\end{equation}
 we have 
$$\omega (\mathfrak{u}_{\C}) = \mathfrak{sp}^{(1,1)} \quad \mbox{ and } \quad \omega (\mathfrak{q}_{\C}) = \mathfrak{sp}^{(2,0)} \oplus \mathfrak{sp}^{(0,2)}.$$
\subsubsection{The positive definite Fock model $\mathcal{P}_+$}\label{matrixinterpretation}
In what follows we will be mostly concerned with the subspace $\mathcal{P}_+$ of $\mathcal{P}$, that is, polynomials in the ``positive'' variables $z_{\alpha,j}$ alone.   It will be useful (both for our computations and to compare our formulas with those of \cite{KashiwaraVergne}) to regard these variables as coordinate functions on the space $\mathrm{M}_{p,n}(\C)$ of $p$ by $n$ matrices.  The following lemma justifies this. In what follows if $U$ is a complex vector space then $\mathrm{Pol}(U)$ denotes the polynomial functions on $U$. 
\begin{lem}
There is an isomorphism of $K' \times K$-modules 
$$\mathcal{P}_+ = \mathrm{Sym}(V_+ \otimes   \C^n  ) \cong \mathrm{Pol}(\mathrm{M}_{p,n}(\C)). $$
\end{lem}
\begin{proof}
We have 
\begin{multline*}
\mathrm{Pol}(\mathrm{M}_{p,n}(\C)) \cong \mathrm{Pol}(\mathrm{Hom}(\C^n,V_+)) \cong \mathrm{Pol}( V_+ \otimes  (\C^n)^*) \\
\cong \mathrm{Sym}( ( V_+ \otimes  (\C^n)^*)^*)\cong \mathrm{Sym}(V_+ \otimes   \C^n  ).
\end{multline*}
Here we use (as we will do repeatedly) that $V_+ \cong V_+^*$ as a $K$-module.  
\end{proof}
We will sometimes regard the $z_{\alpha,j}$ as matrix coordinates in what follows.

\subsection{Some conventions} By abuse of notations we will use, in the following, the same symbols for corresponding objects and operators in both the Schwartz and Fock models. 

We use the convention that indices $\alpha, \beta , \ldots$ run from $1$ to $p$ and indices $\mu , \nu , 
\ldots$ run from $p+1$ to $m$.  In this numbering $K= \SO (p) \times \SO (q)$ acts so that for each $j$, 
the group $\SO(p)$ rotates the variables $z_{\alpha,j}$ and $\SO (q)$ rotates the variables $z_{\mu , j}$

Note that $\p \cong \C^p \otimes (\C^q)^* \cong \mathrm{M}_{p,q} (\C)$, with our convention we let $\omega_{\alpha ,\mu}$ be the linear form which maps an element of $\p$ to its $(\alpha , \mu)$-coordinate.

Finally, for multi-indices $\underline{\alpha} = (\alpha_1 , \ldots , \alpha_q)$ and $\underline{\beta} =
(\beta_1 , \ldots , \beta_{\ell})$ we will write
$$\omega_{\underline{\alpha}} = \omega_{\alpha_1, p+1} \wedge \ldots \wedge \omega_{\alpha_q ,p+q},$$
$$z_{\underline{\alpha} , j} = z_{\alpha_1, j} \ldots z_{\alpha_q,j} ,$$
$$v_{\underline{\beta}} = v_{\beta_1} \otimes \ldots \otimes v_{\beta_{\ell}} .$$

\subsection{} We are interested in the reductive dual pair $\mathrm{O} (V) \times \Mp_{2n} (\R)$ inside $\Mp_{2mn} (\R)$. We suppose that $\mathrm{O} (V)$ and $\Mp_{2n} (\R)$ are embedded in $\Mp_{2mn} (\R)$ in such a way that the Cartan decomposition \eqref{Cartan} of $\mathfrak{sp}_{2mn}$ also induces Cartan decomposition of $\mathfrak{g}$ and $\mathfrak{g}'$. Then $\mathcal{P}$ is a $(\g , K)$-module 
and a $(\g' , K')$-module. We will make use of the structure of $\mathcal{P}$ as a $(\g \oplus \g' , K \cdot K')$-module. We first recall the definition of harmonics (see \cite{HoweT}). 

The Lie algebra $\mathfrak{k} = \mathfrak{o}_p \times \mathfrak{o}_q$ of $K$ is a member of a reductive dual pair $(\mathfrak{k} , \mathfrak{l}')$ where $\mathfrak{l}' = \mathfrak{sp}_{2n} \times \mathfrak{sp}_{2n}$. We can decompose 
$$\mathfrak{l}' = \mathfrak{l}' \- {}^{(2,0)} \oplus  \mathfrak{l}' \- {}^{(1,1)} \oplus  \mathfrak{l}' \- {}^{(0,2)}, \quad \mbox{ where }  \mathfrak{l}' \- {}^{(a,b)} =  \mathfrak{l}' \cap \mathfrak{sp}^{(a,b)}.$$
Then the harmonics are defined by 
$$\mathcal{H} = \mathcal{H} (K) = \left\{ P \in \mathcal{P} \; : \; l(P) = 0 \mbox{ for all } l \in  \mathfrak{l}' \- {}^{(0,2)} \right\}.$$ 

\medskip
\noindent
{\it Remark.} The space $\mathcal{H}$ is smaller that the usual space of harmonics $\mathcal{H} (G)$ in $\mathcal{P}$ 
associated to the  ``indefinite Laplacians'', the latter being associated to the dual pair $(\g , \g')$
rather than $(\mathfrak{k} , \mathfrak{l}')$. The space $\mathcal{H}$ is easily described 
by separating variables: Let $\mathcal{P}_+$, resp. $\mathcal{P}_-$, be the space of all polynomial functions on 
$V_+ \otimes W' \cong V_+^n$, resp. $V_- \otimes W' \cong V_-^n$. It is naturally realized as a subspace of $\mathcal{P} (\C^{mn})$, namely the space of polynomials in the variables $z_{\alpha , j}$, resp. $z_{\mu , j}$. We obviously have $\mathcal{P} = \mathcal{P}_+ \otimes \mathcal{P}_+$. Now denote by 
$\mathcal{H}_+$, resp. $\mathcal{H}_-$, the harmonic polynomials in $V_+^n$, resp. $V_-^n$, (the ``pluriharmonics'' in the terminology of Kashiwara-Vergne \cite{KashiwaraVergne}).  We then have:
$$\mathcal{H} = \mathcal{H}_+ \otimes \mathcal{H}_-.$$

\subsection{Some special harmonic polynomials} \label{antiholomorphicpolynomials}
In this subsection we will introduce  subspaces $\mathcal{H}'(V_+^n)$ and $\mathcal{H}''(V_+^n)$ of $\mathcal{H}_+$ which are closed under polynomial multiplication but not closed under the action of $K =\mathrm{O}(V_+)$.  
 
We begin by introducing   coordinates (that we will call ``Witt coordinates'') that will play a key role  in what follows.  The resulting coordinates $w'_{\alpha,j}$, $w''_{\alpha,j}$, $t_j$, $1 \leq \alpha \leq [\frac{p}{2}]$, $1 \leq j \leq n$, coincide (up to an exchange of order of the indices $\alpha,j$) with  the
coordinates $x_{j \alpha}$, $y_{j \alpha}$, $t_j$, $1 \leq \alpha \leq [\frac{p}{2}]$, $1 \leq j \leq n,$  of Kashiwara and Vergne \cite{KashiwaraVergne} (Kashiwara and Vergne use $i$ instead of $j$ and $\nu$ instead of $\alpha$). 
First we define an ordered {\it Witt basis} $\mathcal{B}$ for $V_+$.  Let $p_0 = [p/2]$. We define an involution
$\alpha \to \alpha'$ of the set $\{1,2,\cdots, 2p_0\}$ by 
 $$ \alpha' = 2p_0 - \alpha +1.$$

If $p$ is even
define
$u'_{\alpha},u''_{\alpha} $ ($1 \leq \alpha \leq p_0 $) by
$$ u'_{\alpha}  = \frac{v_{\alpha} - i v_{\alpha'}}{\sqrt{2}}$$
and 
$$u''_{\alpha} =  \frac{v_{\alpha} + i v_{\alpha'}}{\sqrt{2}}.$$
Then $(u'_1,\cdots, u'_{p_0}, u''_1, \cdots, u''_{p_0})$ is the required ordered Witt basis. 
In case $p$ is odd we define $u'_{\alpha}$ and $u''_{\alpha}$ ($1 \leq \alpha \leq p_0$) as above then add $v_p$ as the last basis vector.  In both cases we will use $\mathcal{B}$ to denote the above ordered basis. 

We note that $u'_{\alpha}$ ($1 \leq \alpha \leq p_0$), and $u''_{\alpha}$ ($1 \leq \alpha \leq p_0$), are isotropic vectors which satisfy
$$(u'_{\alpha}, u'_{\beta}) = 0, \quad  (u''_{\alpha}, u''_{\beta}) = 0 \ \text{ and } \ (u'_{\alpha}, u''_{\beta}) = \delta_{\alpha \beta}\quad \text{ for all }\quad   \alpha,\beta .$$
Of course in the odd case $v_p$ is orthogonal to all the $u'_{\alpha}$'s and $u''_{\beta}$'s.

Define coordinates 
$$(w'_1,\cdots, w'_{p_0},  w''_{p_0},\cdots, w''_1) \quad \mbox{ for }\dim (V_+) \mbox{ even},$$ 
$$\mbox{resp. }  (w'_1,\cdots, w'_{p_0},  w''_{p_0},\cdots, w''_1,t) \quad \mbox{ for } \dim (V_+) \mbox{ odd},$$
by
$$ x= \sum_{\alpha=1}^{p_0} w'_{\alpha}(x) u'_{\alpha} + \sum_{\alpha=1}^{p_0} w''_{\alpha} (x) u''_{\alpha} \quad \text{ in case $p$ is even}$$
and
$$ x= \sum_{\alpha=1}^{p_0} w'_{\alpha}(x) u'_{\alpha} + \sum_{\alpha=1}^{p_0} w''_{\alpha}(x)u''_{\alpha} + t v_p \quad \text{ in case $p$ is odd}.$$
The above coordinates on $V_+$ induce coordinates $w'_{\alpha,k}$, $w''_{\alpha,k}$, $t_k$ on $V^n_+$ for $1 \leq \alpha \leq p$, $1 \leq k \leq n$, and for  $\mathbf{x} = (x_1,\cdots ,x_n) \in V_+^n$ we have 
$$ x_k= \sum_{\alpha=1}^{p_0} w'_{\alpha,k}(\mathbf{x} ) u'_{\alpha} + \sum_{\alpha=1}^{p_0} w''_{\alpha,k}(\mathbf{x} ) u''_\alpha \quad \text{ in case $p$ is even}$$
and
$$x_k= \sum_{\alpha=1}^{p_0} w'_{\alpha,k}(\mathbf{x} ) u'_{\alpha} + \sum_{\alpha=1}^{p_0} w''_{\alpha,k}(\mathbf{x} )u''_{\alpha} + t_k v_p \quad \text{ in case $p$ is odd.}$$
We note the formulas
\begin{enumerate}
\item $z_{\alpha,j}(\mathbf{x}) = (x_j,v_{\alpha})$, $1 \leq \alpha \leq p, 1 \leq j \leq n$
\item $w'_{\alpha,j}(\mathbf{x}) = (x_j,u''_{\alpha})$, $1 \leq \alpha \leq p_0, 1 \leq j \leq n$
\item $w''_{\alpha,j}(\mathbf{x}) = (x_j,u'_{\alpha})$, $1 \leq \alpha \leq p_0, 1 \leq j \leq n$.
\end{enumerate}
We find as a consequence that 
\begin{equation}\label{zandw}
w'_{\alpha,j} = z_{\alpha,j} + i z_{\alpha',j} \ \text{and} \ w''_{\alpha,j} = z_{\alpha,j} - i z_{\alpha',j} , \quad 1 \leq \alpha \leq p_0, \ 1 \leq j \leq n.
\end{equation}

It will be convenient to define $w'_{\alpha,j}$, resp. $w''_{\alpha,j}$, for $\alpha$ satisfying $p_{0} + 1  \leq \alpha \leq 2p_0$ by using equation \eqref{zandw}.
Hence we have
$$w'_{\alpha',j} = i w''_{\alpha,j}, 1 \leq \alpha \leq 2p_0.$$


For both even and odd $p$ we denote  the algebra of polynomials in $w'_{\alpha,j}$ by $\mathcal{H}'(V_+^n)$ and the algebra of polynomials
in $w''_{\alpha,j}$ by $\mathcal{H}''(V_+^n)$. 
The following lemma is critical in what follows. 

\begin{lem} \label{harmonic}
The $\C$-algebras $\mathcal{H}'(V_+^n)$ and $\mathcal{H}''(V_+^n)$ of $\mathcal{P}_+$ lie in the vector space $\mathcal{H}_+$.
\end{lem} 
\begin{proof}
The Laplacians $\Delta_{ij},1 \leq i,j \leq n,$ whose kernels define $\mathcal{H}_+$ are given by sums of {\it mixed} partials
$$\Delta_{ij} = \sum_{\alpha =1}^p \frac{\partial^2}{\partial w'_{\alpha,i} \partial w''_{\alpha,j}} +  \frac{\partial^2}{\partial w'_{\alpha,j} \partial w''_{\alpha,i}}\ \text{for $p$ even}$$
and 
$$\Delta_{ij} = \sum_{\alpha =1}^p \frac{\partial^2}{\partial w'_{\alpha,i} \partial w''_{\alpha,j}} +  \frac{\partial^2}{\partial w'_{\alpha,j} \partial w''_{\alpha,i}} + \frac{\partial^2}{\partial t_i \partial t_j} \ \text{for $p$ odd}.$$
See \cite{KashiwaraVergne}, pages 22 and 26.
\end{proof}
\medskip
\noindent
{\it Remark.} The only harmonic polynomials we will encounter in this Chapter will belong to the subring  $\mathcal{H}''(V_+^n)$.
\medskip
\subsection{The matrix $W''(\mathbf{x})$ and the harmonic polynomials $\Delta_k(\mathbf{x})$}\label{thematrixW'}
We will use $\mathbf{x}$ to denote an $n$-tuple of vectors,  $\mathbf{x} = (x_1, x_2,\cdots,x_n) \in V ^n$. 
Let  $W''(\mathbf{x})$ be the $p_0$ by $n$ matrix  with $(\alpha,j)^{\rm th}$ entry $w''_{\alpha,j},1 \leq \alpha \leq p_0, 1 \leq j \leq n $, that is the coordinates of $x_j$ relative to $u''_1,\cdots,u''_{p_0}$. Following the notation of \cite{KashiwaraVergne} we let $\Delta_k(\mathbf{x})$ be the leading  principal $k$ by $k$ minor of the matrix  $W''(\mathbf{x})$ (by this we mean the determinant of the upper left $k$ by $k$ block). The polynomials $\Delta_k(\mathbf{x}), 1 \leq k \leq n$ belong to $\mathcal{H}''(V_+)$ and hence they belong to $\mathcal{H}$. For $1 \leq k \leq p_0$ we let $W''_k(\mathbf{x})$ be the submatrix obtained by taking the first $k$ rows of $W''(\mathbf{x})$.  We then have the following equation 
\begin{equation} \label{truncatedequivariance}
W''_k(\mathbf{x} g) = W''_k(\mathbf{x})g,  g \in \mathrm{GL}(n,\C), 1 \leq k \leq p_0.
\end{equation}

\subsection{The classes of of Kudla-Millson and Funke-Millson} \label{par:5.3}Let $\lambda$ be a dominant weight for $G$ expressed as in 
\S \ref{par:weight}. Assume that $\lambda$ has at most $n$ nonzero entries. 
By supressing the last  $m_0 -n$ zeroes the dominant weight  $\lambda$ for $G$ gives rise to a dominant
weight $\lambda_1 \geq \ldots \geq \lambda_n$ of $\mathrm{U} (n)$ (also to be denoted $\lambda$)
and as such a finite dimensional irreducible representation $S_{\lambda} (\C^n)$ of $\mathrm{U}(n)$ and thus of $K'$. Here $S_{\lambda} (\C^n)$ denotes the Schur functor (see \cite{FultonHarris}); it occurs as
an irreducible subrepresentation in $(\C^n)^{\otimes \ell}$ and we 
denote by $\iota_{\lambda}$ the inclusion $S_{\lambda} (\C^n ) \rightarrow (\C^{n})^{\otimes \ell}$.

Now let $\ell = \lambda_1 + \ldots + \lambda_n$. Following \cite[p. 296]{FultonHarris}, we may define
the harmonic Schur functor $S_{[\lambda]} (V)$ as the image of the classical Schur functor  $S_{\lambda} (V)$ of $V$ under the $G$-equivariant projection of  $V^{\otimes \ell}$ onto the harmonic tensors.
We denote by $\pi_{[\lambda]}$ the $G$-equivariant projection $V^{\otimes \ell} \rightarrow S_{[\lambda]} (V)$. 
The representation $S_{[\lambda]} (V)$ is irreducible with  highest weight  $\lambda$. Note that all $S_{[\lambda]} (V)$ we will encounter are self-dual. This is obvious if $m$ is odd as all representation of $\SO (m)$ are then self-dual, when $m$ is even this follows from the fact that $\lambda$ has $n < \mathrm{rank} (G)$ nonzero entries. In what follows, if $V'$ is a representation of $K'$ and $k$ is an integer then $V'[k/2]$ will
denote the representation of $K'$ which is the tensor product of $V'$ and $\C_{\frac{k}{2}}$.

In \cite{FM} Funke and Millson construct a nonzero Schwartz form 
$$\varphi_{nq,[\lambda]}\in \mathrm{Hom}_{K'} (S_{\lambda} (\C^n)[m/2] , \mathrm{Hom}_K (\wedge^{nq} \p , \mathcal{P}_+^{(nq + \ell)} \otimes S_{[\lambda]} (V))).$$

The finite dimensional representation $S_{[\lambda]} (V)$ gives rise to a $G$-equivariant hermitian vector bundle with fiber $S_{[\lambda]} (V)$
on the symmetric space $G/K$
where the $G$-action is via $x: (v , gK) \mapsto (x\cdot v, xgK)$. 
For each $w \in S_{\lambda } (\C^n )$ and each $\mathbf{x} \in V^n$ we may therefore interpret $\varphi_{nq,[\lambda]}(w) (\mathbf{x})$ as a differential $nq$-form on $G/K$ which takes values in the vector bundle associated to $S_{[\lambda]} (V)$. 

The main result of \cite{FM}, specialized to our setting, is the following:

\begin{prop} \label{Prop:FM}
For each $w \in S_{\lambda} (\C^n)$ and $\mathbf{x} \in V^n \cong \mathrm{M}_{m,n}(\C)$, 
$\varphi_{nq , [\lambda]} (w) (\mathbf{x})$ defines a {\it closed} differential form on $G/K$ with values in $S_{[\lambda]} (V)$. Moreover  $\varphi_{nq, [\lambda]}$ defines an element in
$$\mathrm{Hom} _{K' }(S_{\lambda} (\C^n)[m/2] , (\wedge^{nq} ( \mathfrak{p}^*) \otimes \mathcal{P} \otimes S_{[\lambda]}(V))^K )$$
or equivalently a section of the homogeneous vector bundle over the Siegel upper-half space corresponding to the representation of $K'$ given
by $S_{\lambda} (\C^n)[m/2]$ with values in the $G$-invariant {\it closed} bundle-valued $nq$-forms  $ A^{nq} ( G/K; \mathcal{P} \otimes S_{[\lambda]}(V) )^G$. 
\end{prop} 

These Schwartz forms are the generalization of the ``scalar-valued'' Schwartz forms considered by Kudla and Millson \cite{KM1,KM2,KM3} to the coefficient case. We now digress to explain how to construct
the forms with trivial coefficients  as restriction of forms associated to the unitary group $\mathrm{U}(p,q)$. 

\subsection{Some special cocycles} \label{5.4} 
The natural embedding $\OO (p,q) \subset \mathrm{U}( p,q)$ yields a totally real embedding
of $G/K$ into the Hermitian symmetric space $U/L$ where $U= \mathrm{U}(p,q)$ 
and $L= \mathrm{U} (p) \times
\mathrm{U}(q)$. The tangent space $\p$ of $G/K$ identifies with the holomorphic tangent space $\p_U^{1,0}$ of $U/L$.

In \cite[Section 5]{BMM2} we have introduced special $(\mathfrak{u} (p,q) , L)$-cocycles $\psi_{aq,bq}$; here we will denote by $\psi_{nq}$ the cocycle $\psi_{nq,0}$. In the case $n=1$, we have:
$$\psi_q = \sum_{\underline{\alpha}} z_{\underline{\alpha}} \otimes \omega_{\underline{\alpha}} \in \mathrm{Hom}_{L} (\wedge^{q,0} \p_U ,  \mathcal{P}_+^{(q)}). $$
 
One important feature of the cocycle is 
that, interpreted as a differential form on $U/L$, the form $\psi_{q}$ is closed, holomorphic and square integrable 
(hence  harmonic) of degree $nq$. For general $n$ we have $\psi_{nq} = \psi_q \wedge \cdots \wedge \psi_q$.

In the Fock model the ``scalar-valued'' Schwartz form
(Kudla-Millson form) $\varphi_{nq,0}$ is --- up to a constant --- the restriction of the holomorphic  form $\psi_{nq}$. In the case $n=1$, it is given by the formula:
$$\varphi_{q,0} = \sum_{\underline{\alpha}} z_{\underline{\alpha}} \otimes \omega_{\underline{\alpha}} \in \mathrm{Hom}_K (\wedge^q \p , \mathcal{P}_+^{(q)} ).$$  
More generally $\varphi_{nq,0}$ is obtained as the (external) wedge-product:
$$\varphi_{nq,0} = \underbrace{\varphi_{q, 0} \wedge \ldots \wedge \varphi_{q,0}}_{n \ \mathrm{times}}.$$ 
We will often write $\varphi_{nq}$ instead of $\varphi_{nq,0}$.

Note that, interpreted as a differential form on $G/K$, the form $\varphi_{nq,0}$ is closed but not harmonic. 

The reader will verify that there are analogous forms $\psi_{nq,\lambda}$ and restriction formulas for the general $\varphi_{nq,\lambda}$
with values in the {\it reducible} representation $S_{\lambda}(V)$  but  there is no such form for the harmonic-valued projection  $\varphi_{nq,[\lambda]}= \pi_{[\lambda]} \circ \varphi_{nq,\lambda}$ with values in the {\it irreducible} representation
$S_{[\lambda]}(V)$. 

\medskip
\noindent
{\it Remark.} By construction the form $\psi_{nq}$ factors through the subspace $[\wedge^{nq,0} \p_U]^{\mathrm{SU} (q)}$ of $\mathrm{SU} (q)$-invariants in $\wedge^{nq,0} \p_U$. In particular the form $\varphi_{nq, [\lambda]}$ defines a element in 
$$ \mathrm{Hom}_{K'} (S_{\lambda} (\C^n)[m/2] , \mathrm{Hom}_K ([\wedge^{nq} \p]^{\mathrm{SL}(q)} , \mathcal{P} \otimes S_{[\lambda]} (V))).$$ 
As announced in the Introduction it will therefore follow from Proposition \ref{prop:8.8} that the subspace of the cohomology 
$H^{\bullet}_{\rm cusp} (X_K , S_{[\lambda]} (V))$ generated by special cycles is in fact contained in 
$H^{\bullet}_{\rm cusp} (X_K , S_{[\lambda]} (V))^{\rm SC}$.
(This of course explains the notation.)

\subsection{Cocycles with coefficients}
In order to go from the forms with trivial coefficients to those with nontrivial coefficients one multiplies $\varphi_{nq,0}$ by a remarkable $K'$-invariant 
element $\varphi_{0,[\lambda]}$ of degree zero in the $(\mathfrak{g},K)$-complex with values in $S_{[\lambda]}(\C^n)^* \otimes \mathcal{P}(V^n) \otimes S_{[\lambda]}(V)$. These elements (as $\lambda$ varies) are projections of the basic element $\varphi_{0,\ell}$ whose properties will be critical to us. Thus (in the Fock model) we have
$$\varphi_{nq,[\lambda]}= \varphi_{nq,0} \cdot \varphi_{0,[\lambda]}.$$

\subsubsection{The $K'$-equivariant family of  zero  $(\mathfrak{so}(p,q),K)$-cochains  $\varphi_{0,\ell}$ } \label{5.5} In \cite{FM} Funke and Millson define (for $\C^n$ the standard representation of $\mathrm{U}(n)$) 
$$\varphi_{0, \ell} \in \mathrm{Hom}_{K' \times S_{\ell}} ( (T^{\ell}(\C^n)) , \mathrm{Hom}_K 
(\wedge^0 \p  , \mathcal{P}_+^{(\ell)}  \otimes T^{\ell}(V_+))) $$
by 
\begin{equation}\label{FMzerocochain}
\varphi_{0, \ell} (e_{I}) = \sum_{\underline{\beta}} z_{\beta_1, i_1} \cdots z_{\beta_{\ell} , i_{\ell}}  \otimes v_{\underline{\beta}}  
\end{equation}
(up to a constant factor) where $e_{I} = e_{i_1} \otimes \ldots \otimes e_{i_{\ell}}$ and $\underline{I}= (i_1,\cdots, i_{\ell})$.
Here $K'$ acts on $T^{\ell}(\C^n)$ and $\mathcal{P}_+^{(\ell)}$ and the symmetric group $S_{\ell}$ acts on $T^{\ell}(\C^n)$ and 
$T^{\ell}(V_+)$. 
They also set: 
$$\varphi_{nq , \ell} = \varphi_{nq , 0} \cdot \varphi_{0, \ell} \in \mathrm{Hom}_{K' \times S_{\ell}} ( T^{\ell}(\C^n)[m/2],   \mathrm{Hom}_K 
(\wedge^{nq} \p , \mathcal{P}_+^{(nq+\ell)}    \otimes T^{\ell}(V)) ,$$
and 
$$\varphi_{nq, [\lambda ]} = (1 \otimes \pi_{[\lambda]}) \circ \varphi_{nq, \ell} \circ \iota_{\lambda}  \in \mathrm{Hom} _{K'}(S_{\lambda} (\C^n) , \mathrm{Hom}_K (\wedge^{nq} \p , \mathcal{P}_+^{(nq+\ell)} \otimes
S_{[\lambda]}(V))).$$

In what follows it will be important to note that $\mathrm{GL}(n,\C)$ acts on $\mathcal{P}_+^{(nq+\ell)} = S^{nq+\ell}(\C^n \otimes V_+)$ by (the action induced by) the standard action of $\mathrm{GL}(n,\C)$.
Also the map $\varphi_{0,\ell}$ has image contained in $\mathrm{Hom}_K 
(\wedge^0 \p , \mathcal{P}^{\ell}_+ \otimes T^{\ell}(V_+))$.  We will now rewrite  $\varphi_{0,\ell}$ to deduce some remarkable properties that it posseses. 
\subsubsection{Three properties of $\varphi_{0,\ell}$}\label{threeproperties}
Note first that  a map $\varphi: U \to W \otimes V$  corresponds to a map $\varphi^*: V^* \otimes U \to W$. Hence, using the isomorphism 
$V_+ \cong V_+^*$ we obtain
$$\varphi^*_{0,\ell} :  T^{\ell}(V_+) \otimes  T^{\ell}(\C^n) \to \mathcal{P}^{(\ell)}_+$$
by the formula
$$\varphi_{0,\ell}^*( v_{\underline{\beta}}\otimes e_{I} ) = (\varphi_{0,\ell}(e_I), v_{\underline{\beta}}).$$
Equation \eqref{FMzerocochain} then  becomes 
\begin{equation} \label{FMzerocochainnewversion}
\varphi_{0,\ell}^*(  v_{\underline{\beta}}\otimes e_{I}  ) = z_{\beta_1, i_1} \cdots z_{\beta_{\ell},i_{\ell}}.
\end{equation}
Rearranging the tensor factors $\C^n$ and $V_+$ we may consider the map $\varphi^*_{0,\ell}$ as a map 
$$\varphi_{0,\ell}^*: T^{\ell}( V_+  \otimes \C^n) \to \mathrm{Sym}^{\ell}( V_+  \otimes \C^n) \cong \mathcal{P}^{\ell}_+.$$

This rearrangement leads immediately to two important properties of $\varphi^*_{0,\ell}$.  

First, we have the following {\it factorization property} of $\varphi^*_{0,\ell}$
that will play a critical role in the  proof of Proposition \ref{finalformula}.
Note that as a special case of Equation \eqref{FMzerocochainnewversion}, the map $\varphi^*_{0,1} :  V_+ \otimes \C^n  \to \mathcal{P}_+$ satifies the equation
$$\varphi_{0,1}^*(  v_{\alpha}\otimes e_j ) = z_{\alpha,j}.$$
We see then that  $\varphi^*_{0,\ell}: T^{\ell}( V_+ \otimes \C^n) \to \mathcal{P}^{\ell}_+$ may be factored as follows.  Given a decomposable $\ell$-tensor
$  \mathbf{x} \otimes \mathbf{z}  =   (x_1 \otimes \cdots \otimes x_{\ell}) \otimes (z_1 \otimes \cdots \otimes z_{\ell}) \in T^{\ell}(V_+) \otimes T^{\ell}(\C^n)$,  
rearrange the tensor factors to obtain $(x_1 \otimes z_1)\otimes \cdots \otimes (x_{\ell} \otimes z_{\ell}) \in  T^{\ell}( V_+ \otimes \C^n ) $. 
Then we have
$$\varphi^*_{0,\ell} (  \mathbf{x} \otimes \mathbf{z} ) = \varphi^*_{0,1}( x_1 \otimes z_1)\varphi^*_{0,1}(x_2 \otimes z_2)\cdots \varphi^*_{0,1}(x_{\ell} \otimes z_{\ell}).$$
From this the following {\it multiplicative property} is clear
\begin{lem} \label{factorizationproperty}
Let $\mathbf{z}_1 \in T^{a}(\C^n),\mathbf{z}_2 \in T^{b}(\C^n),\mathbf{x}_1 \in T^{a}(V_+),\mathbf{x}_2 \in T^{b}(V_+)$ with $a+b = \ell$.
Then
$$\varphi^*_{0,\ell} ( (\mathbf{x}_1 \otimes \mathbf{x}_2) \otimes (\mathbf{z}_1 \otimes \mathbf{z}_2))= \varphi^*_{0,a}( (\mathbf{x}_1 \otimes \mathbf{z}_1 )\cdot \varphi^*_{0,b}(\mathbf{x}_2 \otimes \mathbf{z}_2).$$
\end{lem}
The second property we will need is that    $\varphi^*_{0,\ell}$   is (up to identifications) simply the projection from the $\ell$-th graded summand of the tensor algebra on $ V_+ \otimes \C^n $ to the corresponding summand of the symmetric  algebra.  Hence,
 the map $\varphi^*_{0,\ell}$ descends to give a map
\begin{equation}\label{rightdefinition}
\varphi^*_{0,\ell}: \mathrm{Sym}^{\ell}(V_+ \otimes \C^n ) \to \mathrm{Sym}^{\ell}(V_+ \otimes \C^n )
\end{equation} 
which is clearly the identity map.  
The following lemma is then clear. 

\begin{lem} \label{preservesharmonictensors}
$\varphi^*_{0,\ell}$ carries harmonic tensors to harmonic tensors. 
\end{lem}

From now on we will abuse notation and abbreviate $\varphi^*_{0,\ell}$ to $\varphi_{0,\ell}$ for the rest of this subsection.   

\subsection{The  Vogan-Zuckerman $K$-types associated to the special Schwartz forms $\varphi_{nq,[\lambda]}$ }\label{VZK}
For the rest of this section  the symbols $V$, $V_+$ and $V_ -$ will  mean the complexifications
of the the real vector spaces formerly denoted by these symbols. 

We recall that the Vogan-Zuckerman $K$-type $\mu(\mathfrak{q})$ is the lowest $K$-type of $A_{\mathfrak{q}}$. It may be realized by the $K$-invariant subspace  $V(\mathfrak{q})\subset \wedge^{R}(\p)$ generated by the highest weight vector $e(\mathfrak{q}) \in \wedge^{R}(\mathfrak{u} \cap \mathfrak{p}) \subset \wedge^{R}(\mathfrak{p})$.  The $K$-type $\mu(\mathfrak{q,\lambda})$ is the lowest $K$-type of $A_{\mathfrak{q}}(\lambda)$. It may be realized by the $K$-invariant subspace   $V(\mathfrak{q},\lambda) \subset \wedge^{R}(\p) \otimes E^*$ which is the Cartan product of $V(\mathfrak{q})$ and $E^*$.  

Our goal in this section is to prove the following
\begin{prop} \label{containedinharmonics}
$$\varphi_{nq,[\lambda]}(V(\mathfrak{q},\lambda))  \subset  \mathcal{H}_+.$$
\end{prop}

\medskip
\noindent
{\it Remark.} We remind the reader that we modified $\varphi_{0,\ell}$ to $\varphi^*_{0,\ell}$ in subsection \ref{threeproperties}. This results in a modification of $\varphi_{0,[\lambda]}$ and hence of $\varphi_{nq,[\lambda]}$.  Thus we should have written $\varphi^*_{nq,[\lambda]}$ instead of $\varphi_{nq,[\lambda]}$ in the above theorem and in all that follows.  Since this amounts to rearranging a tensor product we will continue to make this abuse of notation in what follows.

\medskip

The key to proving the Proposition will be to  explicitly compute $\mu(\mathfrak{q}), V(\mathfrak{q})$ and $e(\mathfrak{q})$ and the harmonic Schur functors $S_{[\lambda]}$ for the case in hand in terms of the multilinear algebra of $V_+$ and the form
$(\ ,\ )$.  
 We will define a totally isotropic subspace $E_n \subset V_+ \otimes C$ of dimension $n$ and $\mathfrak{q}$ will be the
stabilizer of a fixed flag in $E_n$.  Anticipating this we change the notation from $V(\mathfrak{q},\lambda)$ to $V(n,\lambda)$.
Also we  take $R =nq$ because for all the parabolics we construct below we will have
$$\mathrm{dim} \ (\mathfrak{u} \cap \p) = nq.$$
For  special orthogonal groups associated to an even dimensional vector space   parabolic subalgebras are not in one-to-one correspondence with  isotropic flags but rather with isotropic oriflammes.   This will not be a problem here. The reason for this comes from  the following considerations.  First  all parabolics we consider here will come from flags of isotropic subspaces in $E_n$.   Second,  there is no difference betweeen oriflammes and flags if all the isotropic subspaces considered are in dimension strictly less than the
middle dimension minus 1. Finally, note that in the even case we have $n < m_0 -1$. See \cite[Chapter 11, p. 158]{Garrett}  for details.

Furthermore, again in the even case, since our highest weight $\lambda$ of $E^*$ has at most $n$ nonzero entries and $n < m_0  - 1< m_0$  the irreducible representation
of $\mathrm{SO}(V)$ with highest weight $\lambda$ will extend to an irreducible representation of $\mathrm{O}(V)$. This extension will be unique up to tensoring with the determinant representation.  In other words for us
there will be no difference (up to tensoring by the determinant representation)  between the representation theory of $\mathrm{SO}(V)$ and $\mathrm{O}(V)$.

\medskip

\subsubsection{The proof of Proposition \ref{containedinharmonics} for the case of trivial coefficients}
We will first prove Proposition \ref{containedinharmonics} for the case of trivial coefficients.

If we are interested only in  obtaining  a representation $A_{\mathfrak{q}}$ which will give  cohomology {\it with trivial coefficients} in degree $nq$
we may take  $\mathfrak{q}$ to be a maximal parabolic, hence to be the stabilizer of a totally isotropic subspace $E' \subset V $. 
We remind the reader that throughout this section $V,V_+$ and $V_-$ are the complexifications of the corresponding real subspaces which we have denoted $V, V_+$ and $V_-$ in the rest of the paper. In order for
$\mathfrak{q}$ to be $\theta$-stable it is necessary and sufficient that $E'$ splits compatibly with the splitting $V  = V_+ \oplus V_ -$.  The simplest way to arrange this is  to choose $E' \subset E'_+ \subset V_+$. Hence, we choose $E'$ to be the $n$-dimensional totally isotropic  subspace $E' =E'_n \subset V_+$ given by
$$E'_n= \mathrm{span}\{u'_1,u'_2,\cdots,u'_n\}.$$              
Now let $E_n''$ be the dual $n$ dimensional subspace of $V_+$ given by
$$E''_n =\mathrm{span}\{u''_1,u''_2,\cdots,u''_n\}.$$ 
Then $E''_n$ is a totally isotropic subspace of $V_+$  of dimension $n$ with $E'_n \cap E_n'' =0$ such that the restriction of $(,)$ to $E'_n + E_n''$ is nondegenerate (so $E'_n$ and $E_n''$ are dually-paired by $(,)$). Let $U = (E'_n + E_n'')^{\perp}$. We will abbreviate $E'_n$ and $E''_n$ to $E'$
and $E''$ henceforth.
We obtain
\begin{equation}\label{decompofV}
V = E' \oplus U \oplus E'' \oplus V_-.
\end{equation}
In what follows we will identify the Lie algebra $\mathfrak{so}(V)$ with $\wedge^2(V)$ by $\rho:\wedge^2(V) \to \mathfrak{so}(V)$ given by
$$ \rho(u \wedge v) (w) = (u,w)v - (v,w)u.$$
The reader will verify that under this identification the Cartan splitting of $\mathfrak{so}(V)$ corresponds to
$$\mathfrak{so}(V) = \mathfrak{k} \oplus \mathfrak{p} = (\wedge^2(V_+) \oplus \wedge^2(V_-)) \oplus (V_+ \otimes V_ -).$$ 
Equation \eqref{decompofV} then induces the following splitting of $\mathfrak{so}(V) \cong \wedge^2(V)$:
\begin{multline*} 
\wedge^2(V) = (E' \otimes E'')\oplus (E' \otimes U) \oplus (E' \otimes  V_-) \oplus (E''\otimes U) \oplus (E'' \otimes V_-) \\ \oplus (U \otimes V_-)
               \oplus \wedge^2(E') \oplus \wedge^2(U) \oplus \wedge^2(E'') \oplus \wedge^2(V_-).
\end{multline*}
The reader will then verify the following lemma concerning the Levi splitting of $\mathfrak{q}$ and its relation with the above Cartan splitting
of $\mathfrak{so}(V)$.  
Recall that $\mathfrak{q}$ is the stabilizer of $E'$.  Let $\mathfrak{q} = \mathfrak{l} \oplus \mathfrak{u}$ be the Levi decomposition.
\begin{lem}
\begin{enumerate}
\item $\mathfrak{q} = [(E' \otimes E'') \oplus (U \otimes V_ -) \oplus \wedge^2 (U)] \oplus [(E'\otimes U) \oplus (E' \otimes V_-) ) \oplus \wedge^2 (E')]$
\item $\mathfrak{l} =(E' \otimes E'') \oplus (U \otimes V_-) \oplus \wedge^2(U)$
\item $\mathfrak{u}= (E'\otimes U) \oplus (E' \otimes V_-) ) \oplus \wedge^2 (E')$
\end{enumerate} 
Hence, we have 
$$\mathfrak{u} \cap \mathfrak{p} =(E' \otimes V_-) =\mathrm{span} (\{ u'_j \wedge v_{p+k}: 1 \leq j \leq n, 1 \leq k \leq q\})$$
whence
$$e(\mathfrak{q}) = [(u'_1 \wedge v_{p+1}) \wedge \cdots \wedge (u'_1 \wedge v_{p+q})] \wedge \cdots \wedge [ (u'_n \wedge v_{p+1}) \wedge \cdots \wedge (u'_n \wedge v_{p+q})].$$
\end{lem}
Next we  describe the Vogan-Zuckerman subspace  $V(n) \subset \wedge^{nq} \p \cong \wedge^{nq}(V_+ \otimes V_ -)$ (underlying the realization of the Vogan-Zuckerman special $K$-type in $\wedge^{nq}\mathfrak{p}$) using the standard formula for the 
decomposition of the exterior power of a tensor product, see equation \eqref{GLdec} or Equation (19), or the formula  on page 80 of  \cite{FultonHarris}. Here $n \times q$ denotes the
partition of $nq$ given by $q$ repeated $n$ times.  

\begin{lem}
The subspace   $V(n)$ of   $\wedge^{nq}(V_+ \otimes V_ -)$ is given by 
$$V(n) = S_{[n \times q]}(V_+) \boxtimes S_{[q \times n]}(V_ -) = S_{[n \times q]}(V_+) \boxtimes \C.$$
\end{lem}
\subsection{} \label{VZK2}  We denote by $V(n , \lambda)$ the {\it Cartan product} of $V(n) \otimes S_{[\lambda]} (V)^*$ i.e. the highest $K$-type of the tensor product $V(n) \otimes S_{[\lambda]} (V)^* \cong V(n) \otimes S_{[\lambda]} (V)$.

We can now prove Proposition \ref{containedinharmonics} for the case of trivial coefficients. We refer the reader to subsection \ref{thematrixW'}
for the definition of the $p_0$ by $n$ matrix $W''(\mathbf{x})$. We recall that $\Delta_n(\mathbf{x})$ is the determinant of the leading principal $n$ by 
$n$ minor of $W''(\mathbf{x})$.

\begin{prop}\label{harmonicvalue}
We have 
\begin{equation*}\label{valueofthecocycle}
\varphi_{nq}(e(\mathfrak{q}))(\mathbf{x}) = \Delta_n(\mathbf{x})^q \in \mathcal{H}''(V_+^n).
\end{equation*}
and consequently
$$\varphi_{nq}(V(n)) \subset  \mathcal{H}_+.$$ 
\end{prop}
\begin{proof} It follows from \cite[Lemma 3.16]{BMM2} that seen as an element of $\wedge^{nq,0} \mathfrak{p}_U$ the vector $e(\mathfrak{q})$ is a Vogan-Zuckerman vector for the theta stable parabolic $\mathfrak{q}_{n,0}$ of $\mathfrak{u}(p,q)$. The cocycle $\varphi_{nq}$ being the restriction of $\psi_{nq}$ the Proposition follows from \cite[Proposition 5.24]{BMM2}.
\end{proof}

\subsection{A  derivation of the formulas  for the simultaneous highest weight harmonic polynomials in Case (1)} \label{simultaneousharmonics}
In this section we will prove the general case of Proposition \ref{containedinharmonics}, in other words, we will prove
that $\varphi_{nq,[\lambda]}(V(n,\lambda))$ takes values in $\mathcal{H}_+$. We will do this by giving an {\it explicit}
formula for $\varphi_{0,[\lambda]}(e_{\lambda} \otimes v^*_{[\lambda]})$ which will obviously be in $\mathcal{H}''(V_+^n)$.
Our computation will give a new  derivation of the formulas of Kashiwara and Vergne for the simultaneous highest weight vectors
in $\mathcal{H}_+$ for their Case (1),  Propositions (6.6) and (6.11) of \cite{KashiwaraVergne}. This derivation will be an immediate consequence of the multiplicative
property, Lemma \ref{factorizationproperty}, of $\varphi_{0,\ell}$ and standard facts in representation theory.

First we need to make an observation.  Note that $V(n , \lambda)$ is the lowest $K$-type of $A_{\mathfrak{q}} (\lambda)$, with $\mathfrak{u} \cap \p$ as above. But if not all the entries of $\lambda$ are equal then we can no longer take $\mathfrak{q}$ to be a {\it maximal} parabolic and we will have to
replace the totally isotropic subspace $E_n$ by a flag in $E_n$.   For example,  if all
the entries of $\lambda$ are different then we must take a full flag in $E_n$.  However for all parabolics $\mathfrak{q}$ obtained we obtain
the same formula for $\mathfrak{u} \cap \p$ and the formula above for $e(\mathfrak{q})$ remains valid.  In fact (because induction in stages is satisfied for
derived functor induction) one can always take $\mathfrak{q}$ to be
the stabilizer of a full flag in $E_n$ and hence will have the property that the Levi subgroup $L_0$ of $Q$ intersected with $G$
will be $\mathrm{U}(1)^n$.  This we will do for the rest of the paper.

Next, we recall that
since $V(n,\lambda)$ is embedded in $\wedge^{nq}(\mathfrak{p}) \otimes S_{[\lambda]}(V)^*$ as the Cartan product, the highest vector of 
$e(\mathfrak{q},\lambda)$ of $V(n,\lambda)$  is given by
\begin{equation} \label{productofweightvectors}
e(\mathfrak{q},\lambda) = e(\mathfrak{q}) \otimes v_{[\lambda]}^*
\end{equation} 
where $v_{[\lambda]}^*$ is the highest weight vector of $S_{[\lambda]}(V)^*$. In what follows let $e_{\lambda}$ be a highest weight vector
of $S_{\lambda}(\C^n)$.
Also since $\varphi_{nq,[\lambda]} = \varphi_{nq,0} \cdot \varphi_{0,[\lambda]}$ we have
\begin{equation}\label{productforFM}
\varphi_{nq,[\lambda]}(e_{\lambda} \otimes e(\mathfrak{q},\lambda)) = \varphi_{nq,0}(e(\mathfrak{q})) \cdot  \varphi_{0,[\lambda]}(e_{\lambda} \otimes v_{[\lambda]}^*).
\end{equation}
Here and in the formula just above we think of $\varphi_{nq,0}$ as an element of $\mathrm{Hom}_K(\wedge^{nq}(\p), \mathcal{P})$
and $\varphi_{0,[\lambda]}$ as an element of 
$\mathrm{Hom}_{K' \times K}(S_{\lambda}(\C^n)[m/2]\otimes S_{[\lambda]}(V)^*, \mathcal{P})$.
Since they both take values in the $\C$-algebra $\mathcal{P}$ we can multiply those values.  The resulting product is what is used in the
equation \eqref{productforFM} and induces the product in the previous equation. Since $\mathcal{H}''(V_+^n)$ is closed under multiplication if we can prove 
that both factors of the right-hand side of equation \eqref{productforFM} are contained in the ring $\mathcal{H}''(V_+^n)$ we can conclude  that 
$$\varphi_{nq,[\lambda]}(e_{\lambda} \otimes e(\mathfrak{q},\lambda)) \in \mathcal{H}''(V_+^n) \subset \mathcal{H}_+.$$
Accordingly, since $V(n,\lambda)$ is an irreducible $K' \times K$-module, the action of $K' \times K$ preserves $\mathcal{H}_+$, and $\varphi_{nq,[\lambda]}$ is a $K' \times K$ homomorphism,   Proposition \ref{containedinharmonics} will follow.

In fact we now prove a  stronger statement than the required $\varphi_{0,[\lambda]}(e_{\lambda} \otimes v_{[\lambda]}^*) \in \mathcal{H}''(V_+^n)$.  Namely we give a new proof of the
formulas of \cite{KashiwaraVergne} in Case (1). This proof makes clear that their formula follows  with very little  computation.  Namely
we first realize the representation $V(\lambda)$ with highest weight $\lambda$ of $\mathrm{O}(V_+)$ in a tensor power of symmetric powers of fundamental representations (exterior powers of the standard representation $V_+$) corresponding
to representing the highest weight $\lambda$ in terms of the fundamental weights $\varpi_j, 1 \leq j \leq p_0$
\begin{equation}\label{fundamentalweights}
\lambda = \sum_{i=1}^{p_0} a_i \varpi_i .
\end{equation}
We then write down the standard realization  of the highest weight vector in this tensor product, see equation \eqref{highestweightvector} below.
The point is that this realization is represented as a product, it is ``factored''. 
We then apply the $K$-homomorphism $\varphi_{0,\ell}$ to this vector using the {\it  multiplicative property}, Lemma \ref{factorizationproperty}  and obtain the desired realization of it in $\mathcal{P}^{\ell}_+$ as a product of powers of leading principal minors of the matrix $W''(\mathbf{x})$. Thus
the new feature of the proof is the existence and factorization property of the map $\varphi_{0,\ell}$.  We now give the details.

The reader will verify that (after changing from the basis of the dual of the Cartan given by the fundamental weights to
the standard basis) that the following formula is the same as those of \cite{KashiwaraVergne}, Proposition (6.6), Case (1) and Proposition (6.11), Case (1). 
\begin{prop} \label{finalformula} Write the highest weight $\lambda$ in terms of the fundamental weights according to
$\lambda = \sum_{i=1}^{n} a_i \varpi_i$. 
Then we have (up to a constant multiple)

 \begin{equation}\label{fformula}
\varphi_{0,[\lambda]}(e_{\lambda}\otimes v_{[\lambda]}^*)(\mathbf{x}) = \Delta_1(\mathbf{x})^{a_1}\Delta_2(\mathbf{x})^{a_2} \cdots \Delta_{n}(\mathbf{x})^{a_n}. 
\end{equation}
and consequently
\begin{equation}\label{goestoharmonics}
\varphi_{0,[\lambda]}(S_{\lambda}(\C^n) \otimes S_{[\lambda]}(V_+)) \subset \mathcal{H}_+
\end{equation}
Combining  equation \eqref{fformula} with the equation  of Proposition \ref{harmonicvalue} we obtain
\begin{equation} \label{vfformula}
\varphi_{nq,[\lambda]}(e_{\lambda}\otimes e(\mathfrak{q}) \otimes  v_{[\lambda]}^*)(\mathbf{x}) = \Delta_1(\mathbf{x})^{a_1}\Delta_2(\mathbf{x})^{a_2} \cdots \Delta_{n}(\mathbf{x})^{a_n +q}.
\end{equation} 
\end{prop}
\begin{cor} \label{whatwewanted}
$$\varphi_{nq,[\lambda]}(e_{\lambda}\otimes v_{[\lambda]}^*) \in \mathcal{H}''(V_+^n) \subset \mathcal{H}_+ .$$
\end{cor}

\begin{proof}
First, as stated, we give the standard realization  of  the  of highest weight vectors of   in the tensor product $T^{\ell}(\C^n) \otimes T^{\ell}(V^*)$, namely we have  the formula
\begin{multline}\label{highestweightvector}
 e_{\lambda} \otimes v^*_{\lambda} = [e_1^{\otimes a_1} \otimes  (e_1 \wedge e_2)^{\otimes a_2} \otimes \cdots \otimes (e_1 \wedge e_2 \wedge \cdots \wedge e_n)^{\otimes a_n}] \\ \otimes  [u_1^{\otimes a_1} \otimes (u_1 \wedge u_2)^{\otimes a_2} \otimes  \cdots \otimes (u_1 \wedge u_2 \wedge \cdots \wedge u_n)^{\otimes a_n}].
\end{multline}
Indeed, note that the above tensor is annihilated by the nilradicals of both Borels.  It is obviously annihilated by $\wedge^2(E)$. Since the rest of  the nilradical of the  Borel subalgebra $\mathfrak{b}$ for $\mathfrak{so}(V_+)$ map  is spanned by the root vectors $u_j \wedge u'_i, 1 \leq j < i \leq p_0,$ that map $u_i$ to $u_j$ with $j<i$ (and $u_1$ to $0$), the claim follows for $\mathfrak{so}(V_+)$.  Similarly the Borel subalgebra  for $\mathfrak{gl}(n,\C)$ is spanned by the elements $E_{ij}, 1 \leq j < i \leq n,$ that map $e_i$ to $e_j$ with $j<i$ (and $e_1$ to $0$) and are zero on all basis vectors other than $e_i$ . Note also that the $u_i$'s are orthogonal isotropic vectors hence the above tensor in the $u_i$'s is a harmonic
tensor in $T^{\ell}(V_+)$. Lastly 
the above vector has weight $\lambda$ by construction. 
Note that with the above  realization of $e_{\lambda} \otimes v^*_{[\lambda]}$ in $T^{\ell}(\C^n \otimes V_+^*)$ we have
$$\varphi_{0,[\lambda]}( e_{\lambda} \otimes v^*_{[\lambda]}) = \varphi_{0,\ell}( e_{\lambda} \otimes v^*_{[\lambda]})$$
where on the right-hand side we consider $e_{\lambda} \otimes v^*_{[\lambda]} \in T^{\ell}(\C^n \otimes V_+^*)$.
We now apply the factorization property of $\varphi_{0,\ell}$, see Lemma \ref{factorizationproperty}. 

Indeed, factoring the right-hand side of equation \eqref{highestweightvector} into $n$ factors (not counting the powers $a_i$) we obtain
\begin{multline*}
\varphi_{0,\ell }(e_{\lambda}\otimes v_{[\lambda]}^*)  = \varphi_{0,1}( e_1 \otimes u_1)^{a_1} \varphi_{0,2}([e_1 \wedge e_2] \otimes  [u_1 \wedge u_2])^{a_2} \cdots \\ \cdots \varphi_{0,n}( [e_1 \wedge \cdots e_n] \otimes [ u_1\wedge \cdots \wedge u_n])^{a_n}.
\end{multline*}
But then observe that
$$\varphi_{0,k}([e_1 \wedge e_2 \wedge \cdots \wedge e_k]\otimes [u_1 \wedge u_2 \wedge \cdots \wedge _k]) = \Delta_k(\mathbf{x}).$$
The Proposition follows.
\end{proof}

\subsection{The derivation of the correspondence of representations on the harmonics}\label{repcorrespondence}
In this subsection we will see how the map 
$\varphi_{0,\ell}: \mathrm{Sym}^{\ell}(\C^n \otimes V_+) \to \mathrm{Sym}^{\ell}(\C^n \otimes V_+) \cong \mathrm{Pol}^{\ell}(\mathrm{M}_{p\times n}(\C))$ induces the decomposition formula for the dual pair $\GL(n,\C) \times \mathrm{O}(V_+)$
acting on $\mathcal{H}^{\ell}(\C^n \otimes V_+)$. In what follows let $P(\ell,n)$ be the set of ordered partitions of $\ell$ into less than or equal to $n$ parts (counting repetitions).  We will assume the known result  that {\it  $\mathrm{GL}(n,\C) \times \mathrm{O}(V_+)$ 
acting on $\mathcal{H}^{\ell}(\C^n \otimes V_+)$ forms a dual pair}, see \cite{Howe1}, and compute what the resulting correspondence is using $\varphi_{0,\ell}$. 
\begin{prop}\label{KVcorrespondence}
Under the assumption $n \leq [p/2]$ the map $\varphi_{0,\ell}$ induces an isomorphism of $\GL(n,\C) \times \mathrm{O}(V_+)$-modules 
$$\varphi_{0,\ell}: \sum_{\lambda \in P(\ell,n)}  S_{\lambda}(\C^n) \boxtimes S_{[\lambda]}(V_+) \to \mathcal{H}^{\ell}(\C^n \otimes V_+).$$
As a consequence the correspondence between $\GL(n,\C)$-modules and $\mathrm{O}(V_+)$ modules induced by the action of the dual pair on the harmonics is $S_{\lambda}(\C^n) \leftrightarrow S_{[\lambda]}(V_+)$ (so loosely put ``take the same partition''). 
\end{prop}
\begin{proof} 
We first recall the decomposition of the $\ell$-th symmetric power of a tensor product, see \cite{FultonHarris}, page 80,
$$\mathrm{Sym}^{\ell}(\C^n \otimes V_+) = \bigoplus_{\lambda \in P(\ell,n)} S_{\lambda}(C^n) \otimes S_{\lambda}(V_+).$$
Hence we obtain an isomorphism of $\GL(n,\C) \times \GL(V_+)$-modules
$$\varphi_{0,\ell} :\bigoplus_{\lambda \in P(\ell,n)} S_{\lambda}(\C^n) \otimes S_{\lambda}(V_+) \to \mathrm{Sym}^{\ell}(\C^n \otimes V_+) \cong \mathrm{Pol}^{\ell}((\C^n)^* \otimes V_+^*).$$
But by equation \eqref{goestoharmonics} of Proposition \ref{finalformula} we have (under the assumption $n \leq [p/2]$), 
$$\varphi_{0,\ell}(S_{\lambda}(\C^n) \otimes S_{[\lambda]}(V)) \subset \mathcal{H}^{\ell}(\C^n \otimes V_+).$$
The map is obviously an injection.  

To prove the map is a surjection let $\lambda \in P(\ell,n)$.  Then, from the assumption preceding the statement of the theorem, the $\mathrm{O}(V_+)$-isotypic subspace 
$$\mathrm{Hom}_{\mathrm{GL}(n,\C)}( S_{\lambda}(\C^n),
\mathcal{H}^{\ell}(\C^n \otimes V_+))$$ 
is an irreducible representation for $\mathrm{O}(V_+)$.  But we have just seen it contains the subspace $S_{[\lambda]}(V_+)$.  Hence it coincides with $S_{[\lambda]}(V_+)$.
From this the proposition follows.
\end{proof}

\subsection{The relation with the work of Kashiwara and Vergne}\label{KVrel}
The previous results are closely related to the work of Kashiwara and Vergne \cite{KashiwaraVergne} studying the action of $\mathrm{GL}(n,\C) \times \mathrm{O}(k)$
on the harmonic polynomials on $M(n,k)$. We first note that  Propositions \ref{harmonicvalue} and
\ref{finalformula} do not follow from the results of \cite{KashiwaraVergne} since we do not know a priori that the cocycles $\varphi_{nq}$ resp.  $\varphi_{nq, [\lambda]}$ take harmonic
values on the highest weight vectors $e(\mathfrak{q})$ resp. $e(\mathfrak{q},\lambda)$ (and this and the results of Howe are the key tools underlying the proof of Theorem \ref{Thm:5.8}).
 
For the benefit of the reader  in comparing the results of  Section \ref{simultaneousharmonics}, Proposition \ref{finalformula}  and Section \ref{repcorrespondence}, Proposition \ref{KVcorrespondence}  with the corresponding results of \cite{KashiwaraVergne} we  provide a dictionary between the notations
of our paper and theirs.  We are studying   the action of the dual pair $\mathrm{GL}(n,\C) \times \mathrm{O}(V_+)$ on the harmonic polynomials $\mathcal{H}_+$. Thus our $p$ corresponds to their $k$ and their $n$ coincides with our $n$.  Their $\ell$ is the rank of $\mathrm{O}(V_+)$ which we have denoted 
$p_0$.   Kashiwara and Vergne take for the Fock module the polynomials on $\mathrm{M}_{n\times p}(\C) \cong V_+^* \otimes \C^n 
 $, that is the $\mathrm{GL}(n,\C) \times \mathrm{O}(V_+)$-module $\mathrm{Sym}( V_+ \otimes (\C^n)^* )$ whereas we take  the polynomials
on $    V_+^* \otimes (\C^n)^*  \cong  V_+ \otimes (\C^n)^* \cong \mathrm{M}_{p\times n}(\C)$ that is 
 $\mathrm{Sym}( V_+ \otimes \C^n )$ for our Fock model.  In the two correspondences between $\mathrm{GL}(n,\C)$ modules and $\mathrm{O}(V_+)$ modules    the  $*$ on the second factor causes their $\mathrm{GL}(n,\C)$-modules to be the contragredients of ours. 

 There are two results in \cite{KashiwaraVergne} that are reproved here using $\varphi_{0,\ell}$  (in ``Case 1'', $n = \mathrm{rank}(\mathrm{GL}(n,\C) \leq p_0 = \mathrm{rank}(\mathrm{O}(V_+))$.
First, we give a new derivation of their formula for the simultaneous $\mathrm{GL}(n,\C) \times \mathrm{O}(V_+)$-highest weight vectors in the
 space of harmonic polynomials $\mathcal{H}_+$ in Proposition \ref{finalformula}. Second, we give a new proof of the correspondence between $\mathrm{GL}(n,\C)$ modules and $\mathrm{O}(V_+)$ modules in Proposition \ref{KVcorrespondence}.  Both of our proofs here are based on the 
properties of the  element $\varphi_{0,\ell}$.  As noted above, the correspondence between irreducible representations of $\mathrm{GL}(n,\C)$ and $\mathrm{O}(V_+)$ is different from that of Kashiwara and Vergne (the representations of $\mathrm{GL}(n,\C)$ they obtain are the contragredients of ours).

\subsection{The computation of $\mathrm{Hom}_K(V(n,\lambda),\mathcal{P})$ as a $U(\mathfrak{g}')$ module}

In this subsection we will prove the following theorem by combining Proposition \ref{containedinharmonics} with a result of Howe.

\begin{thm} \label{Thm:5.8}
As a $(\g' , K')$-module, 
$ \mathrm{Hom}_K (V(n , \lambda ) , \mathcal{P})$
is generated by the restriction of $\varphi_{nq , [\lambda]}$ to $V(n , \lambda)$, i.e.
$$ \mathrm{Hom}_K (V(n , \lambda ) , \mathcal{P}) = U(\g') \varphi_{nq,[\lambda]}.$$
Moreover: there exists a $(\g \times \g ' , K \times K')$-quotient $\mathcal{P} / \mathcal{N}$ of $\mathcal{P}$ such that the $(\g' , K')$-module $\mathrm{Hom}_K (V(n, \lambda) , \mathcal{P} / \mathcal{N})$ is irreducible, generated by the image of $\varphi_{nq , [\lambda]} |_{V(n , \lambda)}$ and isomorphic  
to the underlying $(\g' , K')$-module of the (holomorphic) unitary discrete series representation with lowest $K'$-type  (having highest weight) $S_{\lambda} (\C^n) \otimes \C_{\frac{m}{2}}$.
\end{thm}
Here $U(\g')$ denotes the universal enveloping algebra of $\g'$.

The theorem will be a consequence of general results of Howe and the results obtained above  combined with the following lemmas. 
\subsection{} We consider the decomposition of 
$\mathcal{P}$ into $K$-isotypical components:
$$\mathcal{P}= \bigoplus_{\sigma \in \mathcal{R} (K , \omega)} \mathcal{J}_{\sigma},$$
see \cite[\S 3]{HoweT}.

\medskip
\noindent
{\it Remark.}
However, we have to be careful about two group-theoretic points concerning our maximal compact subgroups $K$ and $K'$.  First we address $K'$. The action of $\GL(n,\C)$ on the Fock model 
$\mathcal{P}$ is the standard action on polynomial functions
{\it twisted by a character}. Recall we identify $\mathcal{P} = \mathcal{P} (\C^{mn})$ with the space  $\mathcal{P} (\mathrm{M}_{m,n} (\C) )$ of 
polynomials on $m$ by $n$ complex matrices, $\mathrm{M}_{m,n} (\C)$. Then the action of the
restriction of the Weil representation $\omega$ to $\GL(n,\C)$ on $\mathcal{P} (\mathrm{M}_{m,n} (\C) )$ is given by the following formula.  Let $Z \in \mathrm{M}_{m,n} (\C)$ and $P \in \mathcal{P} (\mathrm{M}_{m,n} (\C) )$. Then we have for $g \in \GL(n,\C)$
\begin{equation}\label{actioninWeil}
\omega (g)P(Z) = \det(g)^{\frac{p-q}{2}}P(Z g).
\end{equation}
To be precise, we do not get an action of the general linear group but rather of its connected two-fold cover $\widetilde{\GL}(n,\C)$, the metalinear group.  We will ignore this point in what follows as we have done with the difference between the symplectic and metaplectic groups.
Second we address $K$. In what follows the theory of dual pairs requires us to use $\mathrm{O}(V_+)$ and $\mathrm{O}(V_ -)$ below.  The reader will verify
that in fact we may replace $\mathrm{O}(V_+)$ and $\mathrm{O}(V_ -)$ by $\mathrm{SO}(V_+)$ and $\mathrm{SO}(V_ -)$. However we  note  $\varphi_{nq,[\lambda]}$ transforms by a power of the determinant representation of $\mathrm{O}(V_ -)$ which will consequently be ignored.  
Thus in what follows $K$ will denote the product $\mathrm{SO}(V_+) \times \mathrm{SO}(V_ -)$.
The main point that
allows us to make this restriction to the connected group $K = \mathrm{SO}(V_+) \times \mathrm{SO}(V_ -)$  is that 
the restriction of the representation $V(n,\lambda)$ of $\mathrm{O}(V_+)$ to $\mathrm{SO}(V_+)$ is irreducible. 

\medskip

In what follows the key point will be to compute the $V(n,\lambda)$-isotypic component in $\mathcal{H} =\mathcal{H}_+ \otimes \mathcal{H}_-$
as a $\GL(n,\C)$ module under the action induced by the Weil representation.  We will temporarily ignore the twist by $\det^{ (p-q)/2}$. 
Then denoting  the above isotypic subspace by $\mathcal{H}_{V(n,\lambda)}$ we have
\begin{equation}
\mathcal{H}_{V(n,\lambda)} =   \mathrm{Hom}_K(V(n,\lambda), \mathcal{H})\otimes V(n,\lambda).
\end{equation}
Here the first  factor is a $\GL(n,\C)$ module where $\GL(n,\C)$ acts by post-composition.

In what follows it will be very important that the representation $V(n,\lambda)$ of $\mathrm{SO}(V_+) \times \mathrm{SO}(V_ -)$ has trivial
restriction to the second factor.  To keep track of this, up until the end of the proof of Lemma \ref{Ktypegenerates},   we will denote 
 the restriction of the representation $V(n,\lambda)$ to the first factor $\mathrm{SO}(V_+)$ of $K$ by $V(n,\lambda)_+$. 
Thus as a representation of the product $\mathrm{SO}(V_+) \times \mathrm{SO}(V_ -)$ we have
$$V(n,\lambda) = V(n,\lambda)_+ \boxtimes 1$$
and
\begin{multline}\label{bothfactors} 
\mathrm{Hom}_K(V(n,\lambda),\mathcal{H}_+ \otimes \mathcal{H}_ -) \\
\begin{split}
& = \mathrm{Hom}_{\mathrm{SO}(V_+)}(V(n,\lambda)_+,\mathcal{H}_+) \otimes
\mathrm{Hom}_{\mathrm{SO}(V_ -)}(1,\mathcal{H}_-) \\
&= \mathrm{Hom}_{\mathrm{SO}(V_+)}(V(n,\lambda)_+,\mathcal{H}_+) \otimes \C.
\end{split}
\end{multline}
Thus it remains to compute the first tensor factor.

\begin{lem}
The $\GL(n,\C)$-module $\mathrm{Hom}_{\mathrm{SO}(V_+)}(V(n,\lambda)_+,\mathcal{H}_+)$ is the ireducible module with highest weight $(q +\lambda_1 , \ldots , q + \lambda_n)$. Hence we have an isomorphism of $\GL(n,\C)$-modules
$$\mathrm{Hom}_{\mathrm{SO}(V_+)}(V(n,\lambda)_+,\mathcal{H}_+) \cong S_{\lambda}(\C^n) \otimes \C_{q}.$$
\end{lem}
\begin{proof}
The lemma is an immediate consequence of Proposition \ref{KVcorrespondence}. Since the irreducible $\mathrm{O} (V_+)$-module $V(n, \lambda )_+$
has dominant weight $\mu = (q+\lambda_1 , \ldots , q + \lambda_n , 0 \ldots ,0,\cdots, 0)$ it  follows 
from Proposition \ref{KVcorrespondence} that the corresponding  irreducible module 
for $\GL (n, \C)$ is  isomorphic to $S_{\lambda} (\C^n )\otimes \C_{q}$ and consequently has highest weight $(q +\lambda_1 , \ldots , q + \lambda_n)$. .
\end{proof}

Taking into account the  twist of the standard $\GL(n,\C)$-action by $\det^{ (p-q)/2}$ in  the action of the Weil representation on $\mathcal{P}$ we find that the final determinant twist is $q +(p-q)/2 = m/2$.
We conclude:

\begin{lem} \label{Linv}
Under the action coming from the
restriction of the Weil representation the $V(n,\lambda)$ isotypic subspace of $\mathcal{H} = \mathcal{H}_+ \otimes \mathcal{H}_ -$   decomposes as a $\GL(n,\C) \times [\mathrm{SO}(V_+) \times  \mathrm{SO}(V_ -)]$ -module according to :
\begin{equation} \label{DP}
\mathcal{H}_{V(n,\lambda)} = (S_{\lambda}(\C^n) \otimes \C_{\frac{m}{2}} ) \otimes V(n, \lambda).
\end{equation}
\end{lem}

 Recall by Proposition \ref{containedinharmonics}
we have
\begin{equation*} \label{phiTransf}
\varphi_{nq, [\lambda]} \in \mathrm{Hom}_K (V(n , \lambda) , \mathcal{H}_+ \otimes \mathcal{H}_ -)
\end{equation*}
But $\varphi_{nq, [\lambda]}$ takes values in $\mathcal{H}_+$ thus induces an element $\varphi^+_{nq, [\lambda]} \in \mathrm{Hom}_K (V(n , \lambda) , \mathcal{H}_+ )$ such that 
$$\varphi_{nq, [\lambda]} = \varphi^+_{nq, [\lambda]} \otimes 1.$$
 
Note that  $\mathrm{GL}(n) \times \mathrm{O}(V_ +)$ acting on $\mathcal{H}_+$ forms a dual pair and $\mathrm{O}(V_ -)$ acting on $\mathcal{H}_-$
forms a dual pair, $V(n,\lambda) = V(n,\lambda)_+ \boxtimes \C$  and $\varphi_{nq, [\lambda]} = \varphi^+_{nq, [\lambda]} \otimes 1$. We have

\begin{lem}\label{almostdone}
We have:
\begin{enumerate}
\item $ \mathrm{Hom}_{\SO(V_+)} (V(n , \lambda)_+ , \mathcal{H}_+ ) = U(\mathfrak{gl}(n, \C))\cdot \varphi^+_{nq, [\lambda]}$.
\item  $ \mathrm{Hom}_{\SO(V_-)} (\C , \mathcal{H}_ - ) = U(\mathfrak{gl}(n, \C))\cdot 1 = \C$.  
\item $\mathrm{Hom}_K(V(n,\lambda), \mathcal{H}_+  \otimes \mathcal{H}_ -) = [U(\mathfrak{gl}(n, \C))\otimes U(\mathfrak{gl}(n, \C))]\cdot  \varphi_{nq,[\lambda]}$.
\end{enumerate}
\end{lem}

Note that we have a product of dual pairs $(\mathrm{GL}(n,\C) \times  \mathrm{O}(V_+)) \times (\mathrm{GL}(n,\C) \times  \mathrm{O}(V_ -))$. 
Since $\varphi_{nq, [\lambda]}$ takes values in $\mathcal{H}_+ \otimes \C$  the action of  $U(\mathfrak{gl}(n, \C))\otimes U(\mathfrak{gl}(n, \C))$ of (3) of Lemma \ref{almostdone} on $\varphi_{nq, [\lambda]}$ coincides with the ``{\it diagonal} '' action of $U(\mathfrak{gl}(n, \C))$ (i.e. coming from the diagonal inclusion of $\mathrm{GL}(n,\C)$ into the product   of the two first factors of the two dual pairs above). It is critical in what follows that this diagonal  $\mathfrak{gl}(n,\C)$ is the complexification of the Lie algebra  of the maximal compact subgroup $K'$  of the metaplectic group $\mathrm{Sp}_{2n}(\R)$ in our basic dual pair $\mathrm{Sp}_{2n}(\R) \times \mathrm{O}(V)$. 
Hence we obtain the improved  version of (3) of the previous lemma that we will need to prove the first statement of  Theorem \ref{Thm:5.8}.

\begin{lem}
$$ \mathrm{Hom}_K (V(n , \lambda) , \mathcal{H}_+ \otimes \mathcal{H}_ -) = U(\mathfrak{gl}(n, \C))\cdot  \varphi_{nq, [\lambda]}.$$
\end{lem}

\medskip

The next lemma proves the first assertion of Theorem \ref{Thm:5.8}. 
\begin{lem} \label{Ktypegenerates}
$$\mathrm{Hom}_K (V(n , \lambda) , \mathcal{P})= U(\mathfrak{\g}')\cdot  \varphi_{nq,[\lambda]}.$$
\end{lem}
\begin{proof}
By Howe \cite[Proposition 3.1]{HoweT} the space $\mathcal{H} = \mathcal{H}_+ \otimes \mathcal{H}_ -$ generates $\mathcal{P}$ as a 
$U(\g')$-module, i.e.
$$\mathcal{P} = U (\g') ( \mathcal{H}_+ \otimes \mathcal{H}_ -).$$
We obtain:
\begin{eqnarray*}
\mathrm{Hom}_K (V(n , \lambda) , \mathcal{P})) & = & 
 \mathrm{Hom}_K (V(n , \lambda) , U (\g' )( \mathcal{H}_+ \otimes \mathcal{H}_ -)) \\
& = &U(\g') (\mathrm{Hom}_K (V(n , \lambda) , \mathcal{H}_+ \otimes \mathcal{H}_ -)) \\
& = & U(\g' ) U(\mathfrak{gl}(n, \C))\cdot \varphi_{nq, [\lambda]} = U(\g' )\cdot \varphi_{nq, [\lambda]} .
\end{eqnarray*}
\end{proof}

\subsection{} We now prove the second assertion of Theorem \ref{Thm:5.8}. 

It follows from Li \cite{Li} that the $(\g , K)$-module $A_{\mathfrak{q}} (\lambda)$ occurs in Howe's theta correspondence (see \cite[Theorem 2.1]{HoweT}).
In particular: there exists a $(\g \times \g' , K \times K')$ quotient of $\mathcal{P}$
which has the form
$$\mathcal{P} / \mathcal{N} \cong A_{\mathfrak{q}} (\lambda ) \otimes \pi ',$$
where $\pi '$ is a finitely generated, admissible, and quasisimple $(\g' , K')$-module.
This yields a projection
\begin{multline*}
\mathrm{Hom}_K (V(n,\lambda) , \mathcal{P} )  \rightarrow  \mathrm{Hom}_K (V(n,\lambda) , 
\mathcal{P} / \mathcal{N}  )   \\ = \pi' \otimes \mathrm{Hom}_K (V(n,\lambda) , 
A_{\mathfrak{q}} (\lambda )) = \pi' \otimes \C. 
\end{multline*}
But since $\mathrm{Hom}_K (V(n,\lambda) , \mathcal{P} )\cong U(\g') \varphi_{nq,[\lambda]}$, the projection is nonzero and $\pi'$ is irreducible the projection must map the generator 
$\varphi_{nq,[\lambda]} |_{V(n, \lambda )}$ of the $U(\g')$-module $\mathrm{Hom}_K (V(n,\lambda) , \mathcal{P} )$ to a generator of $\pi'$.  
Finally: Li makes Howe's correspondence explicit. In our case $\pi'$ is the underlying 
$(\g', K')$-module of the holomorphic unitary discrete series representation with lowest $K'$-type  $S_{\lambda} (\C^n) \otimes \C_{\frac{m}{2}}$ which is the $K'$ type generated by 
$\varphi_{nq,[\lambda]}$.

This concludes the proof of Theorem \ref{Thm:5.8}.

\subsection{} \label{5.11} {\it Remark.} The $(\g , K)$-module $A_{\mathfrak{q}} (\lambda)$ does not occur in Howe's theta correspondence from a symplectic group smaller than $\Sp_{2n} (\R)$.

\medskip
\noindent
\begin{proof} Let $k<n$. It follows e.g. from \cite[Cor. 3 (a)]{Fulton} that as a $K$-module
$\mathcal{P} (\C^{km})_+ \cong \mathrm{Sym} ((\C^p)^{\oplus k})$ does not contain $V(n)$. As $V(n)$ occurs as a $K$-type in
$A_{\mathfrak{q}} (\lambda) \otimes S_{[\lambda]} (V)$ it follows from the proof of Theorem \ref{Thm:5.8} that $A_{\mathfrak{q}} (\lambda)$ does not occur in Howe's theta correspondence from $\Sp_{2k} (\R)$.
\end{proof}

\part{Geometry of arithmetic manifolds}

\section{Cohomology of arithmetic manifolds} \label{7}

\subsection{Notations} \label{par:7.1}
Let $F$ be a totally real field of degree $d$ and $\A$ the ring of adeles of $F$. Let $V$
be a nondegenerate quadratic space over $F$ with $\dim_F V =m$. We assume that $G=\SO(V)$
is compact at all but one infinite place. We denote by $v_0$ the infinite place where $\SO(V)$ is non
compact and assume that $G(F_{v_0}) = \SO (p,q)$. Let $\widetilde{G} = \GSpin (V)$ be the set of all invertible elements in the even Clifford algebra such that $gVg^{-1} =V$.
There is an exact sequence
\begin{equation}
1 \rightarrow F^* \rightarrow \widetilde{G} \rightarrow G \rightarrow 1,
\end{equation}
where $F^*$ is the subgroup of the center which acts trivially on $V$. We denote by 
$\mathrm{Nspin} : \widetilde{G} \rightarrow F^*$ the spinor norm map and let $\widetilde{G}^{\mathrm{der}}$ be its kernel. 
We finally let
$$D = \SO_0 (p,q) / (\SO (p) \times \SO (q)).$$

\subsection{Arithmetic manifolds} \label{par:7.2} In this paragraph we mainly follow \cite[\S 1]{Kudla}.
For any compact open subgroup 
$K\subset G (\A_f )$, we denote by $\widetilde{K}$ its preimage in $\widetilde{G} (\A_f )$ and let
$$X_K = \widetilde{G}(F ) \backslash (\SO (p,q) \times \widetilde{G}(\A_f)) / (\SO (p) \times \SO (q))\widetilde{K}.$$
The connected components of $X_K$ can be described as follows. Write 
$$\widetilde{G}(\A_f ) = \sqcup_j \widetilde{G}(F)_+ g_j \widetilde{K}.$$
Here $\widetilde{G}(F )_+$ consists of those elements whose spinor norm --- viewed as an element of $F^*$ --- is totally positive, i.e. lies in $F_{\infty +}^* = (\R_+^*)^d$ where $d$ is the degree of $F/\Q$. Then 
$$X_K = \sqcup_j \Gamma_{g_j} \backslash D,$$
where $\Gamma_{g_j}$ is the image in $\SO (p,q)_0$ of 
$$\Gamma_{g_j} ' = \widetilde{G}(F )_+ \cap g_j \widetilde{K} g_j^{-1}.$$
Since the group $\widetilde{G}^{\mathrm{der}}$ is connected, simply connected\footnote{Here we work in the algebraic category: connected means Zariski-connected and a semisimple group $G$ is simply connected if any isogeny 
$G' \rightarrow G$ with $G'$ connected is an isomorphism.} and semisimple, the strong approximation theorem implies (see e.g. \cite[Thm. 5.17]{Milne}) that
$$\pi_0 (X_K ) \cong \widetilde{G}(F )_+ \backslash \widetilde{G}(\A_f ) /\widetilde{K} \cong \A^* / F_c^* \mathrm{Nspin} (\widetilde{K})$$
where $F_c^*$ denote the closure of $F^* F^*_{\infty +}$ in $\A^*$. 

We let $\Gamma_K = \Gamma_1$ and $Y_K = \Gamma_K \backslash D$ be the associated connected
component of $X_K$. These are the arithmetic manifolds we are interested in. Note that the manifolds considered in the introduction are particular cases of these.

\subsection{} \label{par:7.3} The group $\pi_0 (X_K)$ acts on $X_K$ by permutation of the connected 
components:
$$X_K = \sqcup_{\sigma \in \A^* / F_c^* \mathrm{Nspin} (\widetilde{K})} Y_K^{\sigma}.$$
For any $x \in \widetilde{G}(\A_f )$, let $L$ be the projection of $x\widetilde{K}x^{-1}$ in 
$G (\A_f )$. There is a natural map $X_K \rightarrow X_L$ given by right
multiplication by $x^{-1}$. Its restriction to the connected component $Y_K$ gives a map
$$Y_K \rightarrow Y_L^{\mathrm{Nspin} (x)}.$$

\subsection{Differential forms} \label{7.4} Let $(\rho, E)$ be a finite dimensional irreducible representation 
of $G_{\infty} = \SO_0 (p,q)$. Let $K_{\infty} = \SO (p) \times \SO (q)$. The representation $\rho_{|K_{\infty}}$ on $E$ gives rise to a $G_{\infty}$-equivariant Hermitian bundle on $D$, namely, $(E \times G_{\infty}) / K_{\infty}$, where the $K_{\infty}$-action (resp. $G_{\infty}$-action) is given by  
$(v, g) \stackrel{k}{\mapsto} (\rho^* (k^{-1}) v , gk)$ (resp. $(v,g) \stackrel{x}{\mapsto} (v , xg)$). There
is, up to scaling, one $G_{\infty}$-invariant Hermitian metric on $E$. We fix an inner 
product $(,)_{E}$ in this class. We denote this Hermitian vector bundle also by $E$. 
Note that this bundle is $G_{\infty}$-equivariantly 
isomorphic to the trivial vector bundle $E \times D$, where the $G_{\infty}$-action is via 
$x: (v, gK_{\infty}) \mapsto (\rho^* (x) v  , xgK_{\infty})$. 
Smooth sections of $E$ are identified 
with maps $C^{\infty} (G ,E)$ with the property that $f(gk) = \rho^* (k^{-1}) f(g)$.

\subsection{} Now let $\Gamma_K = \Gamma_1$ as above and keep notations as in section \ref{sec:4} with $G$ (resp. $K$) replaced by $G_{\infty}$ (resp. $K_{\infty}$). 
The bundle of $E$-valued differential $k$-forms on 
$Y_K = \Gamma_K \backslash D$ can be identified with the vector bundle associated with the 
$K_{\infty}$-representation $\wedge^k \p^* \otimes E$. Note that $\wedge^k \p^* \otimes E$ is naturally endowed with a $K_{\infty}$-invariant scalar product: the tensor product of $(,)_{E}$ with the scalar
product on $\wedge^k \p^*$ defined by the Riemannian metric on $D$. The space of differentiable 
$E$-valued $k$-forms on $Y_K$, denoted $\Omega^k (Y_K, E)$, is therefore identified with 
$$\left( C^{\infty} (\Gamma_K \backslash G_{\infty}) \otimes E \otimes \wedge^k \p^* \right)^{K_{\infty}} \cong \mathrm{Hom}_{K_{\infty}} \left( \wedge^k \p , C^{\infty} (\Gamma_K \backslash G_{\infty}) \otimes E \right).$$
A compactly supported element $\varphi \in \Omega^k (Y_K, E)$ defines a smooth map 
$\Gamma_K \backslash G_{\infty} \rightarrow \wedge^k \p^* \otimes E$ which satisfies:
$$\varphi (gk)  = \wedge^k \mathrm{ad}_{\p}^* (k^{-1}) \otimes \rho^* (k^{-1}) (\varphi (g)) \quad (g \in G_{\infty} , \quad k \in K_{\infty})$$
so that the norm
$$\varphi \mapsto \int_{Y_K} || \varphi (xK_{\infty}) ||_{\wedge^k \p^* \otimes E}^2 dx$$
is well defined. The space of square integrable $k$-forms $\Omega^k_{(2)} (Y_K , E)$ is the
completion of the space of compactly supported differentiable 
$E$-valued $k$-forms on $Y_K$ with respect to this latter norm.

\subsection{The de Rham complex} The de Rham differential 
$$d: \Omega^k (Y_K , E) \rightarrow
\Omega^{k+1} (Y_K, E)$$ 
turns $\Omega^{\bullet} (Y_K , E)$ into a complex. Let $d^*$ be the formal
adjoint. We refer to $dd^* + d^*d$ as the {\it Laplacian}. It extends to a self-adjoint non-negative
densely defined elliptic operator $\Delta_k^{(2)}$ on $\Omega_{(2)}^k (Y_K , E)$, the form Laplacian.
We let 
$$\mathcal{H}^k (Y_K , E) = \{ \omega \in \Omega_{(2)}^k (Y_K , E) \; : \; \Delta_k^{(2)} \omega = 0 \}$$
be the space of harmonic $k$-forms. Hodge theory shows that $\mathcal{H}^k (Y_K , E)$ is isomorphic
to $\overline{H}_{(2)}^k (Y_K , E)$ --- the {\it reduced} $L^2$-cohomology group --- when $Y_K$ is compact the latter group is just $H^{k} (Y_K , E)$ the $k$-th cohomology group of the de Rham complex $\Omega^{\bullet} (Y_K , E)$. We will mainly work with $\mathcal{H}^k (Y_K , E)$.

\subsection{} \label{7.7} Let $(\pi , V_{\pi})$ be an irreducible $(\g , K_{\infty})$-module and consider the linear map:
$$T_{\pi} : \mathrm{Hom}_{K_{\infty}} (\wedge^* \p , V_{\pi}) \otimes \mathrm{Hom}_{\g , K_{\infty}} (V_{\pi}, L^2 (\Gamma_K \backslash G_{\infty}) \otimes E) \rightarrow \Omega_{(2)}^k (Y_K , E)$$
which maps $\psi \otimes \varphi$ to $\varphi \circ \psi$. The image of $T_{\pi}$ is either orthogonal to  
$\mathcal{H}^k (Y_K , E)$ or $H^k ( \g , K_{\infty} ; V_{\pi} \otimes E) \neq 0$ i.e. $\pi$ is cohomological. In the latter case we denote by $\mathcal{H}^k (Y_K , E)_{\pi}$ (resp. 
$\overline{H}_{(2)}^k (Y_K , E)_{\pi}$) the subspace of $\mathcal{H}^k (Y_K , E)$ (resp. 
$\overline{H}_{(2)}^k (Y_K , E)$) corresponding to the image of $T_{\pi}$. 

A global representation $\sigma \in \mathcal{A}^c (\SO(V))$ with $K$-invariant vectors and such that the restriction of $\sigma_{v_0}$
to $\SO_0 (p,q)$ is isomorphic to $\pi$, and $\sigma_v$ is trivial for every infinite place $v \neq v_0$, contributes to 
$$\mathrm{Hom}_{\g , K_{\infty}} (V_{\pi}, L^2 (\Gamma_K \backslash G_{\infty}) \otimes E).$$ 
We denote by $H_{\mathrm{cusp}}^k (Y_K , E)_{\pi}$ the corresponding subspace of $\mathcal{H}^k (Y_K , E)_{\pi}$ (obtained using the map $T_{\pi}$).

Let $m_{K} (\pi)$ be the multiplicity with which $\pi$ occurs as an irreducible cuspidal summand in $L^2 (\Gamma_K \backslash G_{\infty})$. It follows from Matsushima's formula, see e.g.  \cite{BorelWallach},
that
$$H_{\rm cusp}^k (Y_K , E)_{\pi} \cong m_K (\pi) H^k (\g , K_{\infty} ; V_{\pi} \otimes E).$$

\subsection{} Since $X_K$ is a finite disjoint union of connected manifolds $Y_L$ 
we may easily translate the above definitions into $\mathcal{H}^k (X_K , E)$, $\mathcal{H}^k (X_K , E)_{\pi}$, $H_{\mathrm{cusp}}^k (X_K , E)_{\pi}$, etc$\ldots$ 

We set\footnote{These are only notations. We won't consider such spaces as $\mathrm{Sh} (G)$ or $\mathrm{Sh}^0 (G)$.} 
$$H^{k}_{\rm cusp} (\mathrm{Sh} (G) , E )_{\pi} = \lim_{\substack{\rightarrow \\ K}}  H^{k}_{\rm cusp} (X_K , E )_{\pi}$$
and
$$H^{k}_{\rm cusp} (\mathrm{Sh}^0 (G) , E )_{\pi} = \lim_{\substack{\rightarrow \\ K}} H^{k}_{\rm cusp} (Y_K , E )_{\pi}.$$

Now the inclusion map $Y_K \rightarrow X_K$ yields a surjective map
$$H_{\mathrm{cusp}}^k (X_K , E)_{\pi} \rightarrow H_{\mathrm{cusp}}^k (Y_K , E)_{\pi}.$$
As these inclusions have been chosen in a compatible way we get a surjective map:
$$H^{k}_{\rm cusp} (\mathrm{Sh} (G) , E )_{\pi} \rightarrow H^{k}_{\rm cusp} (\mathrm{Sh}^0 (G) , E )_{\pi}.$$

\subsection{Cohomology classes arising from the $\theta$-correspondence} 
Fix $(\pi , V_{\pi})$ a cohomological irreducible $(\g , K_{\infty})$-module such that
$$H^k (\g , K_{\infty} ; V_{\pi} \otimes E) \neq 0.$$
Note that $H_{\rm cusp}^k (\mathrm{Sh} (G) , E)$ is generated by the images of 
$H^k ( \g , K_{\infty} ; \sigma \otimes E)$ where $\sigma$ varies among all irreducible cuspidal
automorphic representations of $G(\A)$ which occur as irreducible subspaces in the space
of cuspidal automorphic functions in $L^2 (G(F) \backslash G(\A))$ and such that $\sigma_v$ is the trivial representation for each infinite place $v \neq v_0$. We let 
$$H_{\theta}^k (\mathrm{Sh} (G) , E) \subset H_{\rm cusp}^k (\mathrm{Sh} (G) , E)$$
be the subspace generated by those $\sigma \in \mathcal{A}^c (\SO(V))$ that are in the image of the cuspidal $\psi$-theta correspondence from a smaller group, see \S \ref{3.1}. And we
define 
$$H_{\theta}^k (\mathrm{Sh} (G) , E)_{\pi} = H_{\theta}^k (\mathrm{Sh} (G) , E) \cap H_{\rm cusp}^k (\mathrm{Sh} (G) , E)_{\pi}.$$

\medskip
\noindent
{\it Remark.} We can fix one choice of nontrivial additive character $\psi$: Every other nontrivial
additive character is of the form $\psi_t : x \mapsto \psi (tx)$ for some $t \in F^*$. Notations being 
as in \S \ref{par:1.4}, one then easily checks that
$$\Theta_{\psi_t , X}^V (\pi ') = \Theta_{\psi , X}^V (\pi_t ')$$
where the automorphic representation $\pi_t '$ is obtained from twisting $\pi '$ by an automorphism of 
$\Mp (X)$.
We may thus drop explicit reference to $\psi$.

\medskip

We can now state and prove the main theorem of this section.

\begin{thm} \label{thm:main7}
Assume $\pi$ is associated to a Levi subgroup $L= \SO (p-2r,q) \times \mathrm{U}(1)^r$ with $p>2r$ and $m-1 > 3r$. Then: the natural map 
$$H^{k}_{\theta} (\mathrm{Sh} (G) , E )_{\pi} \rightarrow H^{k}_{\rm cusp} (\mathrm{Sh}^0 (G) , E )_{\pi}$$
is surjective.
\end{thm}
\begin{proof} It follows from Corollary \ref{Cor:6.6} and Theorem \ref{Thm1} that 
$H_{\rm cusp}^k (\mathrm{Sh} (G) , E)_{\pi}$ is generated by the images of 
$H^k ( \g , K_{\infty} ; (\sigma \otimes \eta) \otimes E)$ where the representations
$\sigma \in \mathcal{A}^c (\SO(V))$ are in the image of the cuspidal $\psi$-theta correspondence from a smaller group and such that the underlying $(\g , K_{\infty})$-module of $\sigma_{v_0}$ is isomorphic to $\pi$, each $\sigma_v$ (for $v\neq v_0$ infinite) is the trivial representation, and 
$\eta$ varies among all automorphic characters of $G(\A)$. 

Now let $\omega$ be an element of the image of $H^k ( \g , K_{\infty} ; (\sigma \otimes \eta) \otimes E)$ in $H_{\rm cusp}^k (\mathrm{Sh} (G) , E)_{\pi}$. Choose $K \subset G (\A_f )$ a compact open subgroup 
such that $\omega \in H_{\rm cusp}^k (X_K , E)_{\pi}$ and $\eta$ is $\widetilde{K}$-invariant. Seeing 
$\omega$ as an element of 
$$\lim_{\substack{\rightarrow \\ K}} \Omega^k (X_K , E) = \mathrm{Hom}_{K_{\infty}} (\wedge^k \p , C^{\infty} (\widetilde{G} (F) \backslash (\SO (p,q) \times \widetilde{G} (\A_f)) \otimes E)$$
we may form the tensor product $\omega \otimes \eta^{-1}$. It defines an element of the image of $H^k ( \g , K_{\infty} ; \sigma \otimes E)$ in $H_{\rm cusp}^k (X_K , E)_{\pi}$ whose restriction to $Y_K$
is equal to $\omega_{|Y_K}$. We conclude that $\omega$ and $\omega \otimes \eta$ have the same image in $H^{k}_{\rm cusp} (\mathrm{Sh}^0 (G) , E )_{\pi}$ and the theorem follows.
\end{proof}

It is not true in general that the automorphic representations of $\SO(V)$ that are in the image of the 
cuspidal theta corespondence from a smaller group are cuspidal. This is the reason why in the next theorem we assume that $V$ is anisotropic.

One can therefore deduce from Theorem \ref{thm:main7} the following:

\begin{thm} \label{thm:7.13} Assume that $V$ is anisotropic.
Let $r$ be a positive integer such that $p>2r$ and $m-1 > 3r$ and 
let $\pi= A_{\mathfrak{q}} (\lambda)$ be a cohomological $(\g, K_{\infty})$ whose associated Levi subgroup $L$ is isomorphic to $\SO (p-2r, q) \times \mathrm{U}(1)^r$ and such that $\lambda$ has at most $r$ nonzero entries. Then: the global theta correspondence induces an isomorphism between the space of cuspidal holomorphic Siegel modular forms, of weight $S_{\lambda} (\C^r)^* \otimes \C_{- \frac{m}{2}}$ at $v_0$ and of weight $\C_{-\frac{m}{2}}$ at all the others infinite places, on the 
{\rm connected} Shimura variety associated to the symplectic group $\Sp_{2r} |_F$ and the space $H^{rq}_{\rm cusp} (\mathrm{Sh}^0 (G) , S_{[\lambda]} (V))_{\pi}$.
\end{thm}
\begin{proof} The surjectivity follows from Theorem \ref{thm:main7} and Theorem \ref{Thm:5.8}. The injectivity follows from Rallis inner product formula \cite{Rallis}. In our case it is due to Li, see the proof \cite[Theorem 1.1]{Li}. More precisely: let $f_1$, $f_2$  two cuspidal holomorphic Siegel modular forms
of weight $(S_{\lambda} (\C^r )^* \otimes \C_{- \frac{m}{2}}) \otimes \C_{-\frac{m}{2}} \otimes \ldots \otimes \C_{-\frac{m}{2}}$ on the 
{\rm connected} Shimura variety associated to the symplectic group $\Sp_{2r} |_F$. These are functions in $L^2 ( \Mp (X) \backslash \Mp_{2r} (\A))$ which respectively belong to the spaces of two cuspidal 
automorphic representations $\sigma_1'$, $\sigma_2' \in \mathcal{A}^c (\Mp (X))$. And Rallis' inner product formula --- as recalled in \cite[\S 2]{Li} --- implies that if $\phi_1$ and $\phi_2$ are functions in 
$\mathcal{S} (V(\A)^r)$ then:
$$\langle \theta_{\psi , \phi_1}^{f_1} , \theta_{\psi , \phi_2}^{f_2} \rangle = \left\{ 
\begin{array}{ll}
\int_{\Mp (\A)} \langle \omega_{\psi} (h) \phi_1 , \phi_2 \rangle \langle \sigma ' (h) f_1 , f_2 \rangle dh & \mbox{ if } \sigma_1'=\sigma_2' = \sigma ', \\
0 & \mbox{ if } \sigma_1' \neq \sigma_2'.
\end{array} \right.$$
We are thus reduced to the case where $\sigma_1'= \sigma_2'$. But the integral on the right-hand side
then decomposes as a product of local integrals. At each unramified finite place these are special (non-vanishing) values of local $L$-functions, see \cite[\S 5]{Li}. It therefore remains to evaluate the remaining
local factors. And Li proves that these are non-zero in our special case where $\sigma_{v_0} '$ is a holomorphic 
discrete series of weight $S_{\lambda} (\C^r) \otimes \C_{\frac{m}{2}}$ and $\sigma_v$ ($v$ infinite, $v\neq v_0$) is the trivial
representation.  
\end{proof}

\section{Special cycles}\label{specialcyclessection}

\subsection{Notations} \label{8.1} We keep notations as in \S \ref{par:7.1} and keep following the adelization \cite{Kudla}
of the work of  Kudla-Millson. We denote by $( , )$ the quadratic form on $V$ and let $n$ be an integer $0\leq n \leq p$. Given an $n$-tuple $\mathbf{x}=(x_1 ,  \ldots , x_n ) \in V^n$ we let 
$U=U(\mathbf{x})$ be the $F$-subspace of $V$
spanned by the components of $\mathbf{x}$. We write $(\mathbf{x},\mathbf{x})$ for the $n\times n$ symmetric matrix  with $ij$th
entry equal to $(x_i , x_j)$.  
Assume $(\mathbf{x},\mathbf{x})$ is totally positive semidefinite of rank $t$. Equivalently:
as a sub-quadratic space $U \subset V$ is totally positive definite of dimension $t$. 
In particular: $0 \leq t \leq p$ (and $t\leq n$). 
The constructions of the preceeding section can therefore be made with the space $U^{\perp}$ in place
of $V$. Set $H= \SO (U^{\perp})$. There is a natural morphism $H \rightarrow G$. Recall that we can realize $D$ as the set of negative $q$-planes in $V_{v_0}$. We then let $D_H$ be the
subset of $D$ consisting of those $q$-planes which lie in $U^{\perp}_{v_0}$.

\subsection{Special cycles with trivial coefficients} \label{8.2} Let $U=U(\mathbf{x})$ as above. 
Fix $K \subset G(\A_f)$ a compact open subgroup.
As in \S \ref{par:7.2} we write
$$\widetilde{G}(\A_f ) = \sqcup_j \widetilde{G}(F)_+ g_j \widetilde{K}.$$
Recall that 
$$\Gamma_{g_j}' = \widetilde{G}(F)_+ \cap g_j \widetilde{K} g_j^{-1}.$$
We set 
$$\Gamma_{g_j , U}' = \widetilde{H} (F)_+ \cap g_j \widetilde{K} g_j^{-1} = \widetilde{H} (F) \cap \Gamma_{g_j} '.$$ 
Let $\Gamma_{g_j , U}$ be the image of $\Gamma_{g_j , U} '$ in $\SO (p-t , q)_0$. We denote by 
$c(U,g_j,K)$ the image of the natural map
\begin{equation} \label{cycle}
\Gamma_{g_j , U} \backslash D_H \rightarrow  \Gamma_{g_j} \backslash D, \quad \Gamma_{g_j , U} z \mapsto \Gamma_{g_j} z.
\end{equation}

\medskip
\noindent
{\it Remark.} For $K$ small enough the map \eqref{cycle} is an embedding. The cycles $C_U:=c(U,1,K)$ are therefore connected totally geodesic codimension $t$ submanifolds in $Y_K$. These are the totally
geodesic cycles of the introduction for the particular $Y_K$ considered there.

\subsection{} We now introduce composite cycles as follows. For $\beta \in \mathrm{Sym}_n (F)$
totally positive semidefinite, we set
$$\Omega_{\beta} = \left\{ \mathbf{x} \in V^n \; : \; \frac12 (\mathbf{x},\mathbf{x}) = \beta \mbox{ and } \dim U(\mathbf{x}) = \mbox{rank} \beta \right\}.$$
Then $\Gamma_{g_j}'$ acts on $\Omega_{\beta}(F)$ with 
finitely many orbits. Given a $K$-invariant Schwartz function $\varphi \in \mathcal{S} (V(\A_f)^n)$ we define
\begin{equation} \label{cycle2}
Z(\beta , \varphi , K) = \sum_j \sum_{\substack{\mathbf{x} \in \Omega_{\beta} (F) \\ {\rm mod} \ \Gamma_{g_j}'}} \varphi (g_j^{-1} \mathbf{x})  c(U(\mathbf{x}) , g_j , K).
\end{equation}

Let $t = \mathrm{rank} (\beta )$. Suppose first $\beta$ is nonsingular.  Then in Subsection  \ref{cuspidalprojection} we have associated an element of $H_{\rm cusp}^{qt} (X_K)$ to the class of the  cycle $Z(\beta , \varphi , K)$  which we called the cuspidal projection of the class.  We now give a new construction of this projection.

Let $t = \mathrm{rank} (\beta )$. Because it is  rapidly decreasing any cuspidal $q(p-t)$-form can be integrated along $Z(\beta , \varphi , K)$. 
We claim the  canonical pairing between the $qt$-forms with cuspidal coefficients and the $p(q-t)$-forms with cuspidal coefficients is a perfect pairing.  Indeed the forms with cuspidal coefficients are $L^2$ and  they are stable under the restriction of the Hodge star operator so the claim follows.  Hence the induced pairing between
$H_{\rm cusp}^{qt} (X_K)$ and $H_{\rm cusp}^{q(p-t)} (X_K)$ is a perfect pairing.  We can thefore associate to  $Z(\beta , \varphi , K)$ a class
$[\beta , \varphi]^0 \in H^{qt}_{\rm cusp} (X_K)$. We let
$$[\beta , \varphi] : =  [\beta , \varphi]^0 \wedge e_q^{n-t} \in H^{qn}_{\rm cusp} (X_K)$$
where we abusively denote by $e_q$ the Euler form (an invariant $q$-form) dual to the Euler class of Section \ref{sec:4}.

\subsection{Special cycles with nontrivial coefficients} \label{8.4} Following \cite{FM} we now promote
the cycles \eqref{cycle2} to cycles with coefficients.  

Let $\lambda$ be a dominant weight for $G$ expressed as in 
\S \ref{par:weight}. Assume that $\lambda$ has at most $n$ nonzero entries and that $n\leq p$. 
Then $\lambda$ defines a dominant
weight $\lambda_1 \geq \ldots \geq \lambda_n$ of $\mathrm{U} (n)$ 
and as such a finite dimensional irreducible representation $S_{\lambda} (\C^n)$ of $\mathrm{U}(n)$ and 
thus of $K'$.  As above we denote by $S_{[\lambda]} (V)$ the finite dimensional irreducible representation of $G$ 
with highest weight $\lambda$.
 
Fix a neat level $K$ so that each $\Gamma_{g_j , U(\mathbf{x})}$ in \eqref{cycle2} acts trivially on $U(\mathbf{x})$.
The components $x_1 , \ldots , x_n$ of each $\mathbf{x}$ are therefore all fixed by $\Gamma_{g_j , U(\mathbf{x})}$. 
Hence any tensor word in these components will be fixed by $\Gamma_{g_j , U(\mathbf{x})}$. Given 
a tableau $T$ on $\lambda$, see \cite{Fulton},\footnote{Note that a tableau is called a semistandard filling in \cite{FM}.} we set
$$c(U(\mathbf{x}) , g_j , K)_T = c(U(\mathbf{x}) , g_j , K) \otimes \mathbf{x}_T.$$
Here $\mathbf{x}_T \in S_{[\lambda]} (V)$ is the harmonic tensor corresponding to $T$. We can similarly
define $Z(\beta , \varphi , K)_T$ as a cycle with coefficient in $S_{[\lambda]}(V)$. We let 
$[\beta , \varphi]_T^0$, resp. $[\beta , \varphi]_T$, be the corresponding element in 
$ H^{qt}_{\rm cusp} (X_K , S_{[\lambda]} (V))$, resp.  $H^{qn}_{\rm cusp} (X_K , S_{[\lambda]}(V))$. 

We finally define 
$$[\beta , \varphi ]_{\lambda}  \in \mathrm{Hom} (S_{\lambda} (\C^n) , H^{qn}_{\rm cusp} (X_K , S_{[\lambda]} (V))$$ as
the linear map defined by 
$$[\beta , \varphi ]_{\lambda}  (\epsilon_T ) =  [\beta , \varphi]_T$$
where $(\epsilon_1 , \ldots , \epsilon_n)$ is the canonical basis of $\C^n$ and $\epsilon_T$ is the standard basis of $S_{\lambda} (\C^n) $ parametrized by the tableaux on $\lambda$.  

\subsection{Back to the forms of Kudla-Millson and Funke-Millson}
Recall from \S \ref{par:5.3} that we have defined an element 
$$\varphi_{nq,[\lambda]} \in \mathrm{Hom} (S_{\lambda} (\C^n) , \mathrm{Hom}_{K_{\infty}} (\wedge^{nq} \p, \mathbf{S} (V(F_{v_0})^n ) \otimes 
S_{[\lambda]} (V) ))  .$$ 

The quotient space 
$$\widehat{D} := G(F_{v_0}) / K_{\infty} \cong \SO (p,q) / (\SO(p) \times \SO (q))$$
is disconnected and is the disjoint union of two copies of $D$. We let $\Omega^{\bullet} (\widehat{D} , S_{[\lambda]}(V))$
denote the complex of smooth $S_{[\lambda]} (V)$-valued differentiable forms on $\widehat{D}$:
$$\Omega^{\bullet} (\widehat{D} ) = \left[ C^{\infty} (G(F_{v_0}))\otimes S_{[\lambda]} (V) \otimes \wedge^{\bullet} \p^* \right]^{K_{\infty}}.$$
Fixing the base point $z_0 = eK$ in $\widehat{D}$, we have an isomorphism
\begin{multline} \label{eq:isom}
\left[\mathbf{S} (V (F_{v_0})^n) \otimes  \Omega^{nq} (\widehat{D} , S_{[\lambda]} (V)) \right]^{G(F_{v_0})} 
\\ \stackrel{\sim}{\rightarrow}  \left[\mathbf{S} (V(F_{v_0})^n ) \otimes S_{[\lambda]}(V) \otimes \wedge^{nq}  
\p^* \right]^{K_{\infty}},
\end{multline}
given by evaluating at $z_0$. We therefore see $\varphi_{nq, [\lambda]}$ as an element 
$$\varphi_{nq,[\lambda]}\in \mathrm{Hom} \left(S_{\lambda} (\C^n) ,  \left[\mathbf{S} (V (F_{v_0})^n) \otimes  \Omega^{nq} (\widehat{D} , 
S_{[\lambda]}(V)) \right]^{G(F_{v_0})} \right).$$

\subsection{} Consider now a positive definite inner product space $V_+$ of dimension $m$ over $\R$.
We may still consider the Schwartz form $\varphi_{0}  \in \mathcal{S} (V_+^n)$. Recall that through the identification of the Fock space as a subspace of the Schwartz space we get:
$$\varphi_0 (\mathbf{x}) = \exp (-\pi \mathrm{tr} (\mathbf{x},\mathbf{x})).$$
Then, under the Weil representation $\omega_+$ of $\Mp (n, \R)$ associated to $V_+$, we have
$$\omega_+ (k') \varphi_0 = \det (k')^{\frac{m}{2}} \varphi_0.$$ 
If $\mathbf{x} \in V_+^n$ with $\frac12 (\mathbf{x},\mathbf{x}) = \beta \in \mathrm{Sym}_n (\R)$, then for $g' \in \Mp_{2n} (\R)$ 
we set
\begin{equation} \label{eq:whittaker}
W_{\beta} (g' ) = \omega_+ (g') \varphi_0 (\mathbf{x}).
\end{equation}

\subsection{Dual forms} Now we return to the global situation. Let $n$ be an integer with $1 \leq n \leq p$. We fix a level $K$ and a $K$-invariant Schwartz function $\varphi \in \mathcal{S} (V(\A_f)^n)$. 
Define
\begin{multline} \label{eq:phitilde}
\phi = \varphi_{nq,[\lambda]}\otimes \big( \bigotimes_{\substack{v | \infty \\ v \neq v_0}} \varphi_0 \big)
\otimes \varphi \\ \in \mathrm{Hom} \left( S_{\lambda} (\C^n ) , \left[ \mathbf{S} (V(\A)^n) \otimes \Omega^{nq} (\widehat{D}, S_{[\lambda]}(V)) 
\right]^{G (F_{v_0})} \right).
\end{multline}
We may then form the theta function $\theta_{\psi, \phi}(g,g')$ as in \S \ref{par:1.3}. As a function
of $g$ it defines a $S_{[\lambda]} (V)$-valued closed $nq$-form on $X_K$ which we abusively denote by 
$\theta_{nq, \lambda} (g', \varphi)$. Let $[\theta_{nq , \lambda} (g' , \varphi)]$ be the (projection of the) corresponding class 
in $\mathrm{Hom} (S_{\lambda} (\C^n) , H^{nq}_{\rm cusp} (X_K , S_{[\lambda]} (V)^*))$. 

For $g' \in \Mp_{2n}( \R)^d \subset \Mp_{2n} (\A)$ and for $\beta \in \mathrm{Sym}_n (F)$ with $\beta \geq 0$, set 
$$W_{\beta} (g') = \prod_{v | \infty} W_{\beta_v} (g_v ') .$$

The following result  is proved by Funke and Millson \cite[Theorems 7.6 and 7.7]{FM}; it is a generalization to twisted
coefficients of the main theorem of \cite{KM3}. The way we rephrase it here in the adelic language is due to Kudla \cite{Kudla}. Recall that rapidly decreasing $q(n-p)$-forms can be paired with degree $nq$
cohomology classes. We denote by $\langle , \rangle$ this pairing. 

\begin{prop}\label{prop:8.8}
As a function of $g' \in \Mp_{2n} (\A)$ the cohomology class $[\theta_{nq, \lambda} (g', \varphi )]$ is a holomorphic
Siegel modular form of weight $S_{\lambda} (\C^n)^* \otimes \C_{- \frac{m}{2}}$ with coefficients in $H^{qn}_{\rm cusp} (X_K , E_{\lambda}^*)$. Moreover: for any rapidly decreasing closed $q(p-n)$-form $\eta$ on $X_K$ with values in $S_{[\lambda]}(V)$, for  any element $g' \in \Mp_{2n}( \R)^d \subset \Mp_{2n} (\A)$ and for any $K$-invariant map $\varphi \in \mathcal{S} (V(\A_f)^n)$ the Fourier
expansion of $\langle [\theta_{nq, \lambda} (g', \varphi )] , \eta \rangle$ is given by
$$\langle [\theta_{nq, \lambda} (g', \varphi )] , \eta \rangle = \sum_{\beta \geq 0} \langle [\beta , \varphi ]_{\lambda} , \eta \rangle W_{\beta} (g') .$$
\end{prop}

\begin{defn} We let 
$$H^{nq}_{\theta_{nq,\lambda}}(\mathrm{Sh}(G),S_{[\lambda]}(V))_{A_{\mathfrak{q}}(\lambda)}$$ 
be the subspace of 
$H^{nq}_{\theta}(\mathrm{Sh}(G),S_{[\lambda]}(V))_{A_{\mathfrak{q}}(\lambda)}$
generated  by those $\sigma \in \mathcal{A}^c(\SO(V))$ that are in the image of the cuspidal $\psi$-theta correspondence from
$\Mp_{2n}( \mathbb{A})$ where the infinite components $\varphi_v$ of the global Schwartz function $\phi$ satisfies
$$\varphi_{v_0} = \varphi_{nq,[\lambda]} \quad \text{ and } \quad \varphi_v = \varphi_0, \ v|\infty, \ v \neq v_0. $$

We will call the corresponding map from the space of Siegel modular forms tensored with the Schwartz space of the finite adeles to automorphic forms for $\mathrm{SO}(V)$ the {\it special} theta lift and the correspondence
between Siegel modular forms and automorphic forms for $\mathrm{SO}(V)$ the special theta correspondence.  We will denote the special
theta lift evaluated on $f' \otimes \varphi$ by $\theta_{nq,\lambda}(f' \otimes \varphi)$.  
\end{defn}

\section{Main theorem} \label{9}

\subsection{Notations} We keep notations as in the preceeding paragraphs. In particular we let 
$\lambda$ be a dominant weight for $G$ expressed as in 
\S \ref{par:weight}. We assume that $\lambda$ has at most $n$ nonzero entries and that $n\leq p$. 
We let $\mathfrak{q}$ be the $\theta$-stable parabolic subalgebra of $\g$ described in \S
\ref{VZK}. 

We let $SC^{nq} (\mathrm{Sh} (G) , S_{[\lambda]} (V))$ be the subspace of $H^{nq}_{\rm cusp} (\mathrm{Sh} (G) , S_{[\lambda]} (V))_{A_{\mathfrak{q}} (\lambda)}$ spanned by the projection on the
relevant $K_{\infty}$-type of the images of the classes $[\beta , \varphi]_{\lambda}$.

\subsection{} Let $K$ be a compact-open subgroup of $G(\mathbb{A}_f )$. Any  $K$-invariant classes in $SC^{nq} (\mathrm{Sh} (G) , S_{[\lambda]} (V))$ defines a class in $H^{nq}_{\rm cusp} (X_K ,S_{[\lambda]} (V))_{A_{\mathfrak{q}} (\lambda)}$. And it follows from \cite[Corollary 5.11]{Kudla} that 
$$SC^{nq} (X_K , S_{[\lambda]} (V)) := SC^{nq} (\mathrm{Sh} (G) , S_{[\lambda]} (V))^K$$ 
is precisely the subset of $H^{nq}_{\rm cusp} (X_K , S_{[\lambda]} (V))_{A_{\mathfrak{q}} (\lambda)}$
spanned by the projections  of the images of the classes $[\beta , \varphi]_{\lambda}$ for $K$-invariant
functions $\varphi$. 

Note that if $\mathbf{x} \in V^n$ and $(\mathbf{x},\mathbf{x})$ is totally positive semidefinite of rank $t\leq n$, 
the wedge product with $e_q^{n-t}$ of the 
special cycle with coefficients $c(U(\mathbf{x}) , 1 , K)_T$ defines a class in $H^{nq}_{\rm cusp} (Y_K , S_{[\lambda]} (V))$. We let $Z^{nq} (Y_K , S_{[\lambda]} (V))$ be the subspace of $H^{nq}_{\rm cusp} (Y_K , S_{[\lambda]} (V))_{A_{\mathfrak{q}} (\lambda)}$ spanned by the projections of these classes. The restriction map $X_K \rightarrow Y_K$ (restriction to a connected component) obviously yields a map 
\begin{equation} \label{SC-Z}
SC^{nq} (X_K , S_{[\lambda]} (V))  \rightarrow Z^{nq} (Y_K , S_{[\lambda]} (V) ).
\end{equation}
We don't know in general if this map is surjective or not. It will nevertheless follow from Theorem \ref{Thm:main9} that in small degree the map \eqref{SC-Z} is indeed surjective. 

\subsection{The   special theta lift is onto the  $A_{\mathfrak{q}}(\lambda)$-isotypic component of the image of the general theta lift}
In this subsection we carry out what was called step 2 in the introduction. This subsection is the analogue for general $\SO(p,q)$ of subsections 6.8 - 6.11 of \cite{HoffmanHe}.  In particular we now  recall their Lemma 6.9.  We need the following definition of the complex linear  antiautomorphism 
$Z \to Z^*$ of $U(\mathfrak{g}')$. For $Z \in U(\mathfrak{g}')$ with $Z = X_1X_2\cdots X_n$ we define 
$$Z^* = (-1)^n X_n X_{n-1}\cdots X_1.$$
Now for general Schwartz functions $\varphi$ and Siegel automorphic forms $f$ we have the following 
\begin{lem} \label{moveover} For $Z \in U(\mathfrak{g}')$ we have
$$\theta(Zf, \varphi) = \theta(f, Z^* \varphi).$$
\end{lem}
\begin{proof}
See  \cite[Lemma 6.9]{HoffmanHe}.
\end{proof}

We  now show that the  projected {\it special} theta lift $\theta_{nq,\lambda}$ is onto
$H^{nq}_{\theta} (\mathrm{Sh} (G) , S_{[\lambda]} (V))_{A_{\mathfrak{q}} (\lambda)}$ that is we have: 

\begin{thm} \label{StepTwo}
$$H^{nq}_{\theta_{nq, \lambda}} (\mathrm{Sh} (G) , S_{[\lambda]} (V))_{A_{\mathfrak{q}} (\lambda)}= H^{nq}_{\theta} (\mathrm{Sh} (G) , S_{[\lambda]} (V))_{A_{\mathfrak{q}} (\lambda)}  .$$
\end{thm}

\begin{proof} In what follows we will use the following simple observation to produce the commutative diagram \eqref{diagramofthetalift} below.   Suppose $H$ is a group and we have $H$-modules $A,B,U,V$. Suppose further that we have $H$-module homomorphisms $\Phi:U \to V$ and $\Psi: B \to A$.  Then we have a commutative diagram
\begin{equation}\label{genprin}
\begin{CD}
\mathrm{Hom}_H(A,U) @>\Phi_{*}>> \mathrm{Hom}(A,V) \\
@V\Psi^*VV                          @VV\Psi^*V \\
\mathrm{Hom}(B,U)  @>\Phi_{*}>> \mathrm{Hom}(B,V)
\end{CD}
\end{equation} 
Here $\Phi_*$ is postcomposition with $\Phi$ and $\Psi^*$ is precomposition with $\Psi$.

In what follows  $H$ will be the group $K_{\infty}$. We now define $K_{\infty}$-module  homomorphisms $\Phi$ and $\Psi$ that we will concern us here.
We begin with the $(\g,K_{\infty})$-module homomorphism $\Phi$.  
Let $H_{A_{\mathfrak{q}}(\lambda)}$ be the 
subspace of smooth vectors in $L^2(G(\Q) \backslash \ G(\A))$ which is the sum of the spaces $H_{\sigma}$ of those representations  $\sigma \in \mathcal{A}^c(\SO(V))$ such that 
\begin{itemize}
\item $\sigma_{v_0} |_{\SO_0 (p,q)} = A_{\q}(\lambda)$.
\item $\sigma_v$ is the trivial representation for all the infinite places $v \neq v_0$ (note that at such places $v$ we have  $G(F_{v}) \cong \SO(p+q)$). 
\item 
$\sigma$ is in the image of the cuspidal $\psi$-theta correspondence from $\Mp (X)$.
\end{itemize} 
As explained just above Corollary \ref{Cor:6.6}, Remark~\ref{5.11} forces the dimension of the 
symplectic space $X$ to be exactly $2n$. 
We now realize $H_{A_{\q}(\lambda)}$ as a subspace in 
$$L^2 (\widetilde{G} (F) \backslash \SO_0 (p,q) \times \widetilde{G} (\A_f)).$$ 
As explained in \S \ref{7.7} we have:
$$H^{nq}(\mathrm{Sh} (G) , S_{[\lambda]} (V))_{A_{\mathfrak{q}} (\lambda)} \cong H^{nq}(\g,K_{\infty}, H_{A_{\q}(\lambda)}\otimes S_{[\lambda]} (V)).$$
But by Proposition 5.4 of \cite{VZ} we have 
\begin{equation}\label{cohomologyequalsKtype}
H^{nq}(\g,K_{\infty} , H_{A_{\q}(\lambda)}\otimes S_{[\lambda]}(V)) \cong  \mathrm{Hom}_{K_{\infty}}( V(n,\lambda),  H_{A_{\q}(\lambda)}).
\end{equation}

Now let $\pi'_0$ (resp. $\pi'$) 
be the holomorphic discrete series representation of $\mathrm{Mp}(2n,\R)$ with lowest
$K$-type (having highest weight) $S_{\lambda} (\C^n) \otimes \C_{\frac{m}{2}}$ (resp. $\C_{\frac{m}{2}}$). As recalled in \S \ref{3.1} the lowest $K_{\infty}$-type $V(n, \lambda)$ has a canonical lift 
$\widetilde{V(n, \lambda)}$ to $\mathrm{O} (p) \times \mathrm{O} (q)$. We let $\widetilde{A_{\q}(\lambda)}$
be the unique irreducible unitary representation of $\mathrm{O} (p,q)$ with lowest $K$-type $\widetilde{V(n, \lambda)}$ and the same infinitesimal character as $A_{\q}(\lambda)$, see \cite[\S 6.1]{HarrisLi} for more details. It follows from \cite{Li} that $\pi_0 '$ corresponds to
$\widetilde{A_{\q}(\lambda)}$ under the local theta correspondence $\Mp(2n,\R) \times \mathrm{O}(p,q)$ 
and that $\pi'$ corresponds to 
the trivial representation of $\mathrm{O}(p+q)$ under the local theta correspondence $\Mp(2n,\R) \times \mathrm{O}(p+q)$.  Let  $ H'_{\pi_0'}$ be the 
subspace of $L^2(\Mp(2n,\Q) \backslash \Mp(2n,\A))$ which is the sum of the subspaces of smooth vectors $H'_{\sigma'}$  of those representations $\sigma' \in \mathcal{A}^c(\Mp(X))$ such that 
\begin{itemize}
\item $\sigma'_{v_0} = \pi_0'$.
\item $\sigma'_v = \pi'$  for all the infinite places $v \neq v_0$ (note that at such places $v$ we have  $G(F_{v}) \cong \SO(p+q)$). 
\end{itemize} 

We realize the oscillator representation as a $(\g, K_{\infty}) \times G(\A_f)$-module 
in the subspace 
$$\mathbf{S} (V(F_{v_0})^n) \times \mathcal{S} (V(\A_f )^n) \subset \mathbf{S} (V(\A)^n).$$
Here the inclusion maps an element $(\varphi_{\infty} , \varphi)$ of the right-hand side to 
$$\phi = \varphi_{\infty} \otimes \big( \bigotimes_{\substack{v | \infty \\ v \neq v_0}} \varphi_0 \big)
\otimes \varphi$$
where the factors $\varphi_0$ at the infinite places $v$ not equal to $v_0$
are Gaussians, the unique element (up to scalar multiples) of the Fock space $\mathcal{P}(V(F_v)^n)$ which is fixed by the compact group $\SO(V(F_v))$. We abusively write elements of $\mathbf{S} (V(F_{v_0})^n) \times \mathcal{S} (V(\A_f )^n)$ as $\varphi_{\infty} \otimes \varphi$.

From now on we abreviate $H=H_{A_{\q}(\lambda)}$, $H' = H'_{\pi_0'}$ and 
$$\mathbf{S} =\mathbf{S} (V(F_{v_0})^n) \times \mathcal{S} (V(\A_f )^n).$$
It follows from the definition of the global theta lift (see \S \ref{par:1.4}) that for any 
$f \in H'$  and   $\phi \in \mathbf{S}$  
the map $f \otimes \phi \mapsto \theta_{\psi , \phi}^f$ is a 
$(\g , K_{\infty})$-module 
homomorphism from $H'\otimes \mathbf{S}$ to the space $H$. We will drop the dependence of 
$\psi$ henceforth and abbreviate this map to $\theta$ whence $f\otimes \phi \mapsto \theta (f \otimes \phi)$.  Then in the diagram 
\eqref{genprin}
we take $U= H'\otimes \mathbf{S}$ and $V= H$ and $\Phi =\theta$.

We now define the map $\Psi$.  We will take $A$ as above to be the vector space  $\wedge^{nq}(\p)\otimes S_{[\lambda]}(V)^*$, the vector space 
$B$ to be the submodule $V(n,\lambda)$ and $\Psi$ to be the inclusion $i_{V(n,\lambda)} : V(n,\lambda) \to \wedge^{nq}(\p)\otimes S_{[\lambda]}(V)^*$ (note that there is a unique embedding up to scalars and the scalars are not important here). 

From the general diagram \eqref{genprin} we obtain the desired  commutative diagram
\begin{equation} \label{diagramofthetalift}
\begin{CD}
H' \otimes \mathrm{Hom}_{K_{\infty}}(\wedge^{nq}(\p)\otimes S_{[\lambda]}(V)^*, \mathbf{S}) @>\theta_*>> \mathrm{Hom}_{K_{\infty}}(\wedge^{nq}(\p)\otimes S_{[\lambda]}(V)^*, H )  \\
@Vi_{V(n,\lambda)}^*VV       @VVi_{V(n,\lambda)}^*V\\
H'\otimes \mathrm{Hom}_{K_{\infty}}(V(n,\lambda), \mathbf{S}) @>\theta_* >> \mathrm{Hom}_{K_{\infty}}(V(n,\lambda), H )
=H^{nq}(\mathrm{Sh}(G), S_{[\lambda]}(V))
\end{CD}
\end{equation}
Then
$$H^{nq}_{\theta}(\mathrm{Sh}(G), S_{[\lambda]}(V)) = \mathrm{Image}(i_{V(n,\lambda)}^* \circ \theta_*).$$

Now we examine the diagram.  Since $V(n,\lambda)$ is a summand, the   map on the left is {\it onto}.  Also by Corollary \ref{VZKtype} the map on the right
is an isomorphism. 

We now define $U_{\varphi_{nq,[\lambda]}}$ to be 
the affine subspace of $\mathrm{Hom}_{K_{\infty}}(\wedge^{nq}(\p)\otimes S_{[\lambda]}(V)^*, \mathbf{S} )$ defined by $\varphi_{\infty} = \varphi_{nq,[\lambda]}$. 

The theorem is then equivalent to the equation
\begin{equation}\label{firstequivversionforthm95}
i_{V_{nq,\lambda}}^* \circ \theta_* (H' \otimes  U_{\varphi_{nq,[\lambda]}}) = 
\mathrm{Image} (i_{V_{nq,\lambda}}^* \circ \theta_*).
\end{equation}
Since the above diagram is commutative,  equation \eqref{firstequivversionforthm95} holds if and only if we have
\begin{equation}\label{secondequivversionforthm95}
\theta_* \circ   i_{V_{nq,\lambda}}^*(H' \otimes  U_{\varphi_{nq,[\lambda]}}) = \mathrm{Image} (i_{V_{nq,\lambda}}^* \circ \theta_*) =H^{nq}_{\theta} (\mathrm{Sh}(G), S_{[\lambda]}(V))_{A_{\q}(\lambda)}.
\end{equation}
Put $\overline{U}_{\varphi_{nq,[\lambda]}}    =i_{V_{nq,\lambda}}^*(U_{\varphi_{nq,[\lambda]}})$.  
Since the left-hand vertical arrow $i_{V(n,\lambda)}^*$ is onto, equation \eqref{secondequivversionforthm95} holds if and only if 
\begin{equation}\label{thirdequivversionforthm95}
\theta_*(H' \otimes \overline{U}_{\varphi_{nq,[\lambda]}}) = \theta_*(H' \otimes \mathrm{Hom}_{K_{\infty}}(V(n,\lambda), \mathbf{S})).
\end{equation} 

 We now prove equation \eqref{thirdequivversionforthm95}.

To this end let $\xi \in  \theta_*(H' \otimes \mathrm{Hom}_{K_{\infty}}(V(n,\lambda), \mathbf{S}))$.
Hence, by definition, there exists 
 $\phi = \varphi_{\infty} \otimes \varphi \in \mathrm{Hom}_{K_{\infty}} (V(n, \lambda ) , \mathbf{S})$
and $f \in H'$ such that 
\begin{equation}\label{hitxi}
 \theta_*( f \otimes \phi) = \xi.
\end{equation}

We claim that in equation \eqref{hitxi} (up to replacing the component  $f_{v_0}$) we may replace the  factor  $\varphi_{\infty}$ of $\phi$ by $\varphi_{nq,[\lambda]}$
 without changing the right-hand side $\xi$ of equation \eqref{hitxi}.  
Indeed by Theorem \ref{Thm:5.8} 
there exists $Z \in U(\mathfrak{sp}_{2n})$ such that
\begin{equation}\label{everythingcomesfromKM}
 \varphi_{v_0} = Z \varphi_{nq,[\lambda]}.
\end{equation}
Now by Lemma \ref{moveover} (with slightly changed notation)  we have
\begin{equation}
\theta_*( f \otimes Z\phi) = \theta_*( Z^*f \otimes  \phi).
\end{equation}

Hence setting $f' = Z^*f$ we obtain, for all $f \in H'$,
\begin{equation}\label{movedZover}
\begin{split}
\xi  & = \theta_* (f \otimes (\varphi_{\infty} \otimes  \varphi)) \\ 
& =  \theta_* (f \otimes (Z\varphi_{nq,[\lambda]} \otimes  \varphi)) \\
& = \theta_* (Z^*f \otimes  (\varphi_{nq,[\lambda]} \otimes  \varphi)) \\
& = \theta_*( f' \otimes (\varphi_{nq,[\lambda]} \otimes  \varphi)) .
\end{split}
\end{equation}
We conclude that 
the image of the space  $H' \otimes U_{\varphi_{nq,[\lambda]}}$ under $\theta_*$ coincides with the image 
of $H' \otimes \mathrm{Hom}_{K_{\infty}}  (V(n,\lambda), \mathbf{S})$ as required. 

\end{proof}

\medskip
\noindent
{\it Remark.}
The reader will observe that equation \eqref{cohomologyequalsKtype} plays a key role in the paper. Roughly speaking, it converts problems
concerning  the functor 
$$H^{nq}(\g,K_{\infty} ,\bullet \otimes  S_{[\lambda]} (V))$$ 
on $(\g,K_{\infty})$-modules  to problems concerning the functor  
$\mathrm{Hom}_{K_{\infty}}(V(n,\lambda), \bullet)$.
The first functor is not exact whereas the second is. 
\medskip

\subsection{Special cycles span}
In this section we will prove that the special cycles span at any fixed level.  First as a consequence of  Theorem \ref{StepTwo} we have:

\begin{prop} \label{SC=T} We have an inclusion 
$$H^{nq}_{\theta} (\mathrm{Sh} (G) , S_{[\lambda]} (V))_{A_{\mathfrak{q}} (\lambda)}\subset SC^{nq}(\mathrm{Sh}(G), S_{[\lambda]}(V)).$$
\end{prop}

\begin{proof}
By Theorem \ref{StepTwo}, in the statement of Proposition \ref{SC=T} we may replace 
$$H^{nq}_{\theta} (\mathrm{Sh} (G) , S_{[\lambda]} (V))_{A_{\mathfrak{q}} (\lambda)}$$ 
by 
$$H^{nq}_{\theta_{nq,\lambda}} (\mathrm{Sh} (G) , S_{[\lambda]} (V))_{A_{\mathfrak{q}} (\lambda)}.$$
We then use Proposition \ref{prop:8.8}. Since cusp forms are rapidly decreasing we can pair classes in $H^{nq}_{\rm cusp} (\mathrm{Sh} (G) , 
S_{[\lambda]} (V) )$ with classes in $H^{(n-p)q}_{\rm cusp} (\mathrm{Sh} (G) , S_{[\lambda]} (V)^* )$. We denote by 
$\langle , \rangle$ this pairing. It is a perfect
pairing. Hence letting $SC^{nq}(\mathrm{Sh}(G), S_{[\lambda]}(V))^{\perp}$ and $H^{nq}_{\theta} (\mathrm{Sh} (G) , S_{[\lambda]} (V))_{A_{\mathfrak{q}} (\lambda)}^{\perp}$ denote the respective
annihilators in $H^{(n-p)q}_{\rm cusp} (\mathrm{Sh} (G) , S_{[\lambda]} (V)^* )$ it suffices to prove
\begin{equation} \label{inclusofperps}
SC^{nq}(\mathrm{Sh}(G), S_{[\lambda]}(V))^{\perp} \subset H^{nq}_{\theta} (\mathrm{Sh} (G) , S_{[\lambda]} (V))_{A_{\mathfrak{q}} (\lambda)}^{\perp}.
\end{equation}
To this end let  
$\eta \in SC^{nq}(\mathrm{Sh}(G), S_{[\lambda]}(V))^{\perp} \subset  H^{(p-n)q}_{\rm cusp} (\mathrm{Sh} (G) , 
S_{[\lambda]} (V)^* )$. Assume $\eta$ is $K$-invariant for some level $K$. 
It then follows from proposition \ref{prop:8.8} (see \cite[Theorem 7.7]{FM} in the classical setting) that for $g' \in \Mp_{2n} (\R)^d \subset \Mp_{2n} (\A)$ the Fourier
expansion of the Siegel modular form 
$$\theta_{\varphi} (\eta) := \langle [\theta_{nq , \lambda} (g' , \varphi )] , \eta \rangle (=  \int_{X_K} \theta_{nq , \lambda} (g' , \varphi ) \wedge \eta )$$ 
is given by
$$\theta_{\varphi} ( \eta) = \sum_{\beta \geq 0} \langle [\beta , \varphi]_{\lambda} , \eta \rangle W_{\beta} (g').$$
In particular: since $\eta$ is orthogonal to the subspace 
$$SC^{nq} (\mathrm{Sh} (G) , S_{[\lambda]} (V)) \subset H^{nq}_{\rm cusp} (\mathrm{Sh} (G) , 
S_{[\lambda]} (V) )$$ 
generated by the projections of the forms in the image of $[\beta , \varphi]_{\lambda}$ then all the Fourier 
coefficients of the Siegel modular form $\theta_{\varphi} ( \eta)$ vanish and therefore $\theta_{\varphi}
(\eta )=0$. But then for any $f\in H_{\sigma'}$ ($\sigma ' \in \mathcal{A}^c (\Mp (X))$), we have:
\begin{equation*}
\int_{X_K} \eta \wedge \theta (f , \varphi)[\varphi_{nq,[\lambda]}]  = \int_{\Mp (X) \backslash \Mp_{2n} (\A)}
 \theta_{\varphi } (\eta) f(g') dg' = 0. 
\end{equation*}
This forces $\eta$ to be orthogonal to all forms $\theta (f , \varphi)[\varphi_{nq,[\lambda]}]$ hence we have verified equation \eqref{inclusofperps} and hence proved the proposition. 
\end{proof}

\subsection{}{\it Remark.} We do not know if the reverse inclusion holds in proposition \ref{SC=T}. It obviously follows from theorem \ref{thm:main7} that it holds when $p>2n$ and $m-1 > 3n$ but the relation between 
cycles and $\theta$-lifts is less transparent near the middle degree, we address this problem in 
conjecture \ref{Conj2:periods}.

\subsection{} Now Proposition \ref{SC=T} and Theorem \ref{thm:main7} imply\footnote{Note that the $(\g , K_{\infty})$-module $A_{\mathfrak{q}} (\lambda)$ is associated to the Levi subgroup 
$L= \SO (p-2n , q)\times  \mathrm{U} (1)^n$.} that if $p>2n$ and $m-1 > 3n$
the natural map 
$$SC^{nq} (\mathrm{Sh}(G) , S_{[\lambda]} (V)) \rightarrow H^{nq}_{\rm cusp} (\mathrm{Sh}^0 (G) , S_{[\lambda]} (V))_{A_{\mathfrak{q}} (\lambda)}$$
is surjective. 

Taking invariants under a compact open subgroup $K \subset G (\A_f)$ we get the following:

\begin{thm} \label{Thm:main9}
Suppose $p>2n$ and $m-1 > 3n$. Then, for any compact open subgroup $K \subset G(\A_f)$, the natural 
map
$$SC^{nq} (X_K , S_{[\lambda]} (V)) \to H^{nq}_{\rm cusp} (Y_K , S_{[\lambda]} (V))_{A_{\mathfrak{q}} (\lambda)}$$
is surjective. In particular: 
$$Z^{nq} (Y_K , S_{[\lambda]} (V)) = H^{nq}_{\rm cusp} (Y_K , S_{[\lambda]} (V))_{A_{\mathfrak{q}} (\lambda)}.$$
\end{thm}

\subsection{}{\it Remark.} As motivated in \S \ref{6.13} we believe that the conditions $p>2n$ and $m-1 > 3n$ are necessary. 
The decomposition of $L^2 (G(F) \backslash G(\A))$ into irreducible automorphic representations yields
a decomposition of $H^{\bullet} (\mathrm{Sh}^0 (G) , S_{[\lambda]} (V))$ and, when $m-1\leq 3n$, one may try to classify which automorphic representations contribute to the part of the cohomology generated 
by special cycles. This is adressed in Theorem \ref{Conj:periods} and Conjecture \ref{Conj2:periods}. 

\medskip

\part{Applications}

\section{Hyperbolic manifolds} \label{10}

\subsection{Notations}
Let notations be as in the preceeding section. Assume moreover that $G$ 
is anisotropic over $F$ and that $q=1$. For any compact open subgroup $K \subset G(\A_f)$,  the 
connected component $Y_K$ is therefore a closed congruence hyperbolic $p$-manifold. These are
called ``of simple type'' in the introduction. We point out that in that case the Euler form is trivial ($q$ is odd) so that special cycles are totally geodesic cycles.

\subsection{Proof of Theorem \ref{thm:intro1}} 
First note that if $n < \frac{p}{3}$ then $2n <p$ and $3n < p=m -1$ so that Theorem \ref{Thm:main9} applies. We consider the case where $\lambda = 0$.

Let $K \subset G(\A_f)$ be a compact open subgroup. Since $Y_K$ is closed the cohomology of $Y_K$ is the same as its cuspidal cohomology. And it follows from the first example of \S \ref{par:weight} that $A_{\mathfrak{q}} (0)$ --- with $\mathfrak{q}$ as in \S \ref{VZK} --- is the 
unique cohomological module with trivial coefficients which occurs in degree $n$. The natural projection
map
$$H^{n} (Y_K , \C) \to H^{n}_{\rm cusp} (Y_K , \C)_{A_{\mathfrak{q}} (0)}$$
is therefore an isomorphism and we conclude from Theorem \ref{Thm:main9} that $H^{n} (Y_K , \C)$
is spanned by the classes of the special (totally geodesic) cycles. As these define {\it rational} cohomology classes Theorem \ref{thm:intro1} follows.

\medskip
\noindent
{\it Remark.} As motivated in \S \ref{6.13}, the bound $n < \frac{p}{3}$ is certainly optimal in general as cohomology classes
of $3$-dimensional hyperbolic manifolds are not generated by classes of special classes in general, see Proposition \ref{prop:dim3} for an explicit counterexample when 
$p=3$ and $n=1$.

\medskip

\subsection{Proof of Theorem \ref{thm:intro2}}
First note that if $p \geq 4$ then $m-1 >3$ so that  Theorem \ref{Thm:main9} applies. We consider the case where $\lambda = (1, 0 , \ldots , 0)$ or $(2,0 , \ldots , 0)$. In the first case $S_{[\lambda]} (V) \cong \C^{p+1}$ and in the second case $S_{[\lambda]} (V)$ is isomorphic to the complexification of $\mathcal{H}^2 (\R^{p+1})$ --- the space of harmonic (for the Minkowski metric) degree two polynomials on $\R^{p+1}$.

Let $K \subset G(\A_f)$ be a compact open subgroup. Since $Y_K$ is closed the cohomology of $Y_K$ is the same as its cuspidal cohomology. And it follows from the second example of \S \ref{par:weight} that $A_{\mathfrak{q}} (\lambda )$ --- with $\mathfrak{q}$ as in \S \ref{VZK} --- is the 
unique cohomological module with coefficients in $S_{[\lambda]} (V)$ which occurs in degree $1$. The natural projection map
$$H^{1} (Y_K , S_{[\lambda]} (V) ) \to H^{1}_{\rm cusp} (Y_K , S_{[\lambda]} (V) )_{A_{\mathfrak{q}} (\lambda )}$$
is therefore an isomorphism and we conclude from Theorem \ref{Thm:main9} that $H^{1} (Y_K , S_{[\lambda]} (V))$
is spanned by the classes of the special cycles. As these define {\it real} cohomology classes Theorem \ref{thm:intro2} follows.

\section{Shimura varieties associated to $\OO (p,2)$} \label{11}

\subsection{Notations}
Let notations be as in section \ref{9}. Assume moreover that $G$ 
is anisotropic over $F$ and that $q=2$. For any compact open subgroup $K \subset G(\A_f)$,  the 
connected component $Y_K$ is therefore a closed projective complex manifold of (complex) dimension $p$. These are the connected Shimura varieties associated to $\OO (p,2)$ in the introduction. 
We consider the case where $\lambda = 0$. 

\subsection{The Lefschetz class} In that setting $D$ is a bounded symmetric domain in $\C^p$. Let 
$k (z_1 , z_2)$ be its Bergmann kernel function and $\Omega (z) = \partial \overline{\partial} \log k (z,z)$
the associated K\"ahler form. For each compact open subgroup $K \subset G(\A_f)$, 
$\frac{1}{2i \pi} \Omega$ induces a $(1,1)$-form on $Y_K$ which is the Chern form of the canonical bundle $\mathcal{K}_{Y_K}$ of $Y_K$. As such it defines a {\it rational} cohomology class dual to  
(a possible rational multiple of) a complex subvariety: given a projective embedding this is the class of (a possible rational multiple of) a hyperplane section. 
Moreover: the cup product with $\frac{1}{2i \pi} \Omega$ induces the Lefschetz operator 
$$L : H^k (Y_K , \Q) \rightarrow H^{k+2} (Y_K , \Q)$$
on cohomology. This is the same operator as the multiplication by the Euler form $e_2$. 

\subsection{} According to the Hard Lefschetz theorem the map
$$L^k : H^{p-k} (Y_K , \C) \to H^{p+k} (Y_K , \C)$$
is an isomorphism; and if we define the {\it primitive} cohomology 
$$H^{p-k}_{\rm prim} (Y_K , \C) = \mathrm{ker} ( L^{k+1} : H^{p-k} (Y_K , \C) \to H^{p+k+2} (Y_K , \C))$$
then we have the {\it Lefschetz decomposition} 
$$H^m (Y_K , \C) = \oplus_k L^k H^{m-2k}_{\rm prim} (Y_K , \C).$$

The Lefschetz decomposition is compatible with the Hodge decomposition. In particular we set: 
$$H^{n,n}_{\rm prim} (Y_K , \C) = H^{n,n} (Y_K , \C) \cap H^{2n}_{\rm prim} (Y_K , \C).$$
And it follows from the third example of \S \ref{par:weight} that 
$$H^{n,n}_{\rm prim} (Y_K , \C) = H^{2n} (Y_K , \C)_{A_{n,n}}.$$

\subsection{Proof of Theorem \ref{thm:intro3}}
First note that if $n < \frac{p+1}{3}$ then $2n < p$ and $3n < p+1= m  -1$ so that Theorem \ref{Thm:main9} applies. Note moreover that $A_{\mathfrak{q}} (0) = A_{n,n}$.

Let $K \subset G(\A_f)$ be a compact open subgroup. Since $Y_K$ is closed the cohomology of $Y_K$ is the same as its cuspidal cohomology. The natural projection
map
$$H^{2n} (Y_K , \C) \to H^{2n}_{\rm cusp} (Y_K , \C)_{A_{\mathfrak{q}} (0)}$$
is nothing else but the projection onto $H^{n,n}_{\rm prim} (Y_K , \C)$ and we conclude 
from Theorem \ref{Thm:main9} that $H^{n,n}_{\rm prim} (Y_K , \C)$
is spanned by the projection of the classes of the special cycles. Those are complex subvarieties of $Y_K$.

To conclude the proof we recall that 
$$H^{n,n} (Y_K , \C) = \oplus_k L^k H^{n-k , n-k}_{\rm prim} (Y_K , \C).$$
And since $n-k\leq n$ the preceeding paragraph applies to show that $H^{n-k , n-k}_{\rm prim} (Y_K , \C)$ is spanned by the projection of classes of complex subvarieties of $Y_K$. As wedging with $L^k$ amounts to take the intersection with another complex subvariety we conclude that the whole 
$H^{n,n} (Y_K , \C)$ is spanned by the classes of complex subvarieties of $Y_K$. These are rational
classes we can therefore conclude that every rational cohomology class of type $(n,n)$ on $Y_K$ 
is a linear combination with rational coefficients of the cohomology classes of complex subvarieties
of $Y_K$.

\section{Arithmetic manifolds associated to $\SO(p,q)$}

\subsection{} In the general case where $G=\SO (p,q)$ we can proceed as in the preceeding section 
to deduce Corollary \ref{cor:intro} from Theorem \ref{thm:intro4} and Proposition \ref{P:cohomrep}. Note however that in this general case we have to use both totally geodesic cycles associated to $\SO (p-n,q)$ (related to cohomological representations associated to a Levi of the form $L=C \times \SO_0 (p-2n,q)$) and totally geodesic cycles associated to $\SO(p,q-n)$ (related
to cohomological representations associated to a Levi of the form $L=C \times \SO_0 (p,q-2n)$). 

\subsection{Proof of Theorem \ref{thm:intro6}} Recall that the arithmetic manifold $Y=Y_K$ is associated to a non-degenerate quadratic space $(V, (,))$ over a totally real field $F$. Consider an totally imaginary quadratic extension $L/F$ and let  $\tau$ be the corresponding Galois involution. Then
$$h(x,y) = (x , \tau (y)) \quad (x, y \in V \otimes_F L)$$
defines a non-degenerate hermitian form. We let $U$ be the the corresponding special unitary group
$\mathrm{SU}(h)$. Then $U$ is defined over $F$ and $U(F \otimes_{\Q} \R) \cong \mathrm{SU}(p,q)
\times \mathrm{SU} (m)^{d-1}$ where $d$ is the degree of $F$. 

We define $\mathrm{Sh}(U)$ as in the orthogonal case and consider $S_{\lambda} (V)$ as finite $U$-module (here $V$ is complexified). The group $\mathrm{U}(h)$ is a member of a reductive dual pair 
$$\mathrm{U}(W) \times  \mathrm{U}(h) \subset \Sp (V \otimes W)$$ 
where $\mathrm{U}(W)$ is the centralizer of $J$ (our choice of positive definite complex structure on $W$). In this way we can define the ($\psi$-)theta correspondence from $\mathrm{U}(W)$ to $\mathrm{U}(h)$ where the test functions $\phi$ still vary in $\mathbf{S} (V(\A)^n)$.

At the place $v_0$ we have introduced the cocycle $\psi_q$, see \S \ref{5.4}. We set
$$\psi_{nq, \ell} = (\underbrace{\psi_q \wedge \ldots \wedge \psi_q}_{n \ \mathrm{times}}) \cdot \varphi_{0, \ell}$$
and 
$$\psi_{nq, \lambda} = (1 \otimes \pi_{\lambda}) \circ \psi_{nq, \ell}(\bullet ) \circ \iota_{\lambda},$$
where $\pi_{\lambda} : V^{\otimes \ell} \rightarrow S_{\lambda} (V)$ is the natural projection. We may therefore define
$$H_{\theta_{nq, \lambda}}^{nq,0} (\mathrm{Sh}(U) , S_{\lambda} (V))$$
as the subspace of $H^{nq,0} (\mathrm{Sh}(U) , S_{\lambda} (V))$ generated by those $\sigma \in 
\mathcal{A}^c (\mathrm{SU}(h))$ that are in the image of the cuspidal $\psi$-theta correspondence from
$\mathrm{U}(W)$ where the infinite components $\varphi_v$ of the global Schwartz function $\phi$
satisfies 
$$\varphi_{v_0} = \psi_{nq, \lambda} \quad \mbox{ and } \quad \varphi_v = \varphi_0, \ v|\infty, \ v \neq v_0. $$
We consider the natural map 
\begin{equation} \label{resmap}
H_{\theta_{nq, \lambda}}^{nq,0} (\mathrm{Sh}(U) , S_{\lambda} (V)) \rightarrow H^{nq}_{\theta_{nq, \lambda}} (\mathrm{Sh} (G) , S_{[\lambda]} (V))_{A_{\mathfrak{q}}(\lambda)}.
\end{equation}
obtained by composing the restriction map with 
the projection $\pi_{[\lambda]}$ onto the harmonic tensors.

\begin{lem}
The map \eqref{resmap} is onto.
\end{lem}
\begin{proof} This follows immediately from two facts:
First the reductive dual pairs $(\mathrm{U}(W), 
\mathrm{U}(h))$ and $(\Mp(W) , \SO (V))$ form a see-saw pair, in the terminology of Kudla \cite{Kudla3}.
Secondly, the restriction of the holomorphic  form $\psi_{nq, \lambda}$ composed with 
the projection $\pi_{[\lambda]}$ onto the harmonic tensors is equal to $\varphi_{nq, [\lambda]}$, see \S \ref{5.4}. 
\end{proof}

It therefore follows from Theorem \ref{StepTwo} and Theorem \ref{thm:main7} that if 
$2n < p$ and $3n<m -1$, the natural map:
$$H^{nq,0} (\mathrm{Sh}(U) , S_{\lambda} (V)) \rightarrow H^{nq}_{\rm cusp} (\mathrm{Sh}^0 (G) , S_{[\lambda]} (V))_{A_{\mathfrak{q}}(\lambda)}$$
is surjective. Taking invariants under a compact open subgroup $L \subset U(\A_f)$ such that $L \cap G(F)=K$ we get Theorem
\ref{thm:intro6} with 
$$Y^{\C} = \Lambda_L \backslash D^{\C},$$ 
where 
$$D^{\C} = \mathrm{SU}(p,q) / \mathrm{S}(\mathrm{U}(p) \times \mathrm{U}(q))$$
and 
$\Lambda_L$ is the image of $U(F) \cap L$ inside $\mathrm{SU}(p,q)$.

\section{Growth of Betti numbers}

\subsection{Notations} Let $F$ be a totally real field and $\mathcal{O}$ its ring of integer. Let $V$
be a nondegenerate quadratic space over $F$ with $\dim_F V =m$. We assume that $G=\SO(V)$
is anisotropic over $F$ and compact at all but one infinite place. We denote by $v_0$ the infinite place where $\SO(V)$ is non compact and assume that $G(F_{v_0}) = \SO (p,1)$.
We fix $L$ an integral lattice in $V$ and let $\Phi= G (\mathcal{O})$ be the subgroup of $G(F)$ consisting of those elements that take $L$ into itself. We let $\mathfrak{b}$ be an ideal in $\mathcal{O}$
and let $\Gamma = \Gamma (\mathfrak{b}$) be the congruence subgroup of $\Phi$ of level 
$\mathfrak{b}$ (that is, the elements of $\Phi$ that are congruent to the identity modulo $\mathfrak{b}$). 
We let $\p$ be a prime ideal of $\mathcal{O}$ which we assume to be prime to $\mathfrak{b}$. 
Set $\Gamma (\p^k) = \Gamma (\mathfrak{b} \p^k)$. The quotients $Y_{\Gamma (\p^k)}$ are real hyperbolic $p$-manifolds and, as first explained in \cite{CossuttaMarshall}, the combination of Theorem \ref{thm:main7} and works of Cossutta \cite[Theorem 2.16]{Cossutta} and Cossutta-Marshall \cite[Theorem 1]{CossuttaMarshall} implies the following strengthening of a conjecture of Sarnak
and Xue \cite{SarnakXue} in this case.

\begin{thm} \label{thm:growth}
Suppose $i <\frac{p}{3}$. Then:
$$b_i (\Gamma (\p^k)) \ll \mathrm{vol} (Y_{\Gamma (\p^k)})^{\frac{2i}{p}}.$$
\end{thm}
\begin{proof} Cossutta and Marshall prove this inequality with 
$b_i (\Gamma (\p^k))$ replaced by the dimension of the part of $H^i (\Gamma (\p^k))$ which come from
$\theta$-lift. But it follows from the first example of  \S \ref{par:weight} 
that only one cohomological module can contribute to $H^i (\Gamma (\p^k))$ and since  $i <\frac{p}{3}$ it follows from Theorem \ref{thm:main7} that the classes which come from
$\theta$-lifts generate $H^i (\Gamma (\p^k))$. Theorem \ref{thm:growth} follows.
\end{proof}

\subsection{}{\it Remark.} 1. Raising the level $\mathfrak{b}$ one can prove that the upper bound given by Theorem \ref{thm:growth} is sharp when $p$ is even, see \cite[Theorem 1]{CossuttaMarshall}.

2. The main theorem of Cossutta and Marshall is not limited to hyperbolic manifolds. In conjunction with
Theorem \ref{thm:main7} one may get similar asymptotic results for the multiplicity of cohomological automorphic forms in orthogonal groups. In fact Theorem \ref{thm:7.13} relates the multiplicities of certain cohomological automorphic forms to the multiplicities of certain Siegel modular forms. The latter are much easier to deal with using limit formulas as in \cite{Clozel,Savin}. The main issue then is to control
the level; this is exactly what Cossutta and Marshall manage to do. 

\section{Periods of automorphic forms}

\subsection{Notations} We keep notations as in \S \ref{par:7.1} and \ref{8.1}. In particular we let $F$ be a totally real field of degree $d$ and denote by $\A$ its ring of adeles of $F$. Let $V$
be a nondegenerate quadratic space over $F$ with $\dim_F V =m$. We assume that $G=\SO(V)$
is compact at all but one infinite place. We denote by $v_0$ the infinite place where $\SO(V)$ is non
compact and assume that $G(F_{v_0}) = \SO (p,q)$. Hereafter $U$ will always denote a totally positive definite subquadratic space of dimension $n \leq p$ in $V$. And we denote by $H$ the group $\SO (U^{\perp})$.

We furthermore let $(\pi, V_{\pi})$ be an irreducible $(\g, K_{\infty})$-module such that 
$$H^k (\g , K_{\infty} ; V_{\pi} \otimes E) \neq 0$$ 
for some finite dimensional irreducible representation $(\rho , E)$ of $\SO_0 (p,q)$ of dominant weight $\lambda$ with at most $n$ nonzero entries. 

For any cusp form $f$ in $L^2 (\widetilde{G} (F) \backslash (\SO (p,q) \times \widetilde{G}(\A_f))$ and any character $\chi$ of $\widehat{\pi}_0 = \A^* / F_c^*$, we define the period integral
$$P(f , U , \chi) = \int_{\widetilde{H}(F) \backslash (\SO(p-n, q) \times \widetilde{H}(\A_f))} f (h) \chi (\mathrm{Nspin}_{U^{\perp}} (h)) dh$$
where $dh = \otimes_v dh_v$ is a fixed Haar measure on $\widetilde{H}(F) \backslash (\SO(p-n, q) \times \widetilde{H}(\A_f))$. 

\subsection{Distinguished representations} Let $(\sigma , V_{\sigma})$ be an irreducible cuspidal 
automorphic representation of $G(\A)$ which occurs as an irreducible subspace $V_{\sigma}$ in the space of cuspidal automorphic functions in $L^2 (G(F) \backslash G(\A))$ and such that $\sigma_v$ 
is trivial for any infinite place $v \neq v_0$. 

Note that a function $f \in V_{\sigma}$ lifts, in a canonical way, as a function in 
$$L^2 (\widetilde{G} (F) \backslash (\SO (p,q) \times \widetilde{G}(\A_f))).$$ 
By a slight abuse of notation we still denote by $f$ this function.
We call $P(f , U , \chi)$ the $(\chi , U)$-period of $\sigma$. We write $P(\sigma , U , \chi) \neq 0$ if there
exists $f \in V_{\sigma}$ such that $P(f,U,\chi)$ is nonzero and $P(\sigma , U , \chi)=0$ otherwise. 
Let us say that $\sigma$ is {\it $\chi$-distinguished} if $P (\sigma , U , \chi) \neq 0$ for some $U$. 

\subsection{}  Given a place $v$ of $F$ we denote by $V(\sigma_v)$ the space of the representation $\sigma_v$. We furthermore fix 
an isomorphism 
$$\tau : \otimes_v ' V(\sigma_v) \rightarrow V_{\sigma} \subset L^2 (G(F) \backslash G(\A)).$$
We finally write $\sigma_f = \otimes_{v\not \; | \infty }'  \sigma_{v}$ and $V(\sigma_f) = \otimes_{v\not \; | \infty }' V( \sigma_{v})$. 

We are only be concerned with the very special automorphic representation $(\sigma , V_{\sigma})$ such that the restriction of $\sigma_{v_0}$ to $\SO_0 (p,q)$ is isomorphic to $\pi$. We abusively write ``$\sigma_{v_0} \cong \pi$'' for ``the restriction of $\sigma_{v_0}$ to $\SO_0 (p,q)$ is isomorphic to $\pi$''.  
 
Classes in the cuspidal cohomology $H_{\rm cusp}^{\bullet} (\mathrm{Sh}(G) , E)$ are represented
by cuspidal automorphic forms. We may therefore decompose $H_{\rm cusp}^{nq} (\mathrm{Sh}(G) , E)_{\pi}$ as a sum
$$H_{\rm cusp}^{\bullet} (\mathrm{Sh}(G) , E)_{\pi} = \bigoplus_{(\sigma, V_{\sigma}) \; : \; \sigma_{v_0} \cong \pi} 
H_{\rm cusp}^{\bullet} (\mathrm{Sh}(G) , E)(\sigma) $$
where we sum over irreducible cuspidal 
automorphic representations $(\sigma, V_{\sigma})$ of $G(\A)$ which occurs as an irreducible subspace $V_{\sigma}$ in the space of cuspidal automorphic functions in $L^2 (G(F) \backslash G(\A))$. 

\subsection{} Recall from \S \ref{8.2} and \ref{8.4} that for $K$ small enough we have associated to a subspace $U$ and a tableau $T$ some connected cycles-with-coefficient $c(U,g_j,K)_T=  c(U,g_j,K)\otimes \epsilon_T$ where $\epsilon_T \in E$. We restrict to the case where the $g_j$'s belong to $H(\A_f )$ so that
the $c(U,g_j,K)$ are the images the connected components of 
\begin{equation} \label{eq:SC}
\widetilde{H} (F) \backslash (\SO (p-n, q) \times \widetilde{H} (\A_f)) / (SO(p-n) \times SO(q)) \widetilde{K \cap H}.
\end{equation}
Any character $\chi$ of finite order of $\A^* / F_c^*$ defines a locally constant function on \eqref{eq:SC}.
Let $Z_{U,T}^{\chi} = \sum_j \chi (\mathrm{Nspin}_U (g_j)) c(U,g_j, K)_T$ and $[Z_{U,T}^{\chi}]$ the corresponding 
cohomology class in $H^{qn}_{\rm cusp} (X_K , E)$.
We let $Z^{nq} (\mathrm{Sh}(G) , E)(\sigma )$ be the subspace of $H^{nq}_{\rm cusp} (\mathrm{Sh} (G) , E)(\sigma )$ spanned by the projection of the classes $[Z_{U,T}^{\chi}]$. 

\begin{prop} \label{prop:13.5}
Suppose that $\pi$ is associated to a Levi subgroup $L= \SO(p-2n, q) \times \mathrm{U} (1)^n$ with $p-2n \geq 0$. If the space $Z^{nq} (\mathrm{Sh}(G) , E)(\sigma )$ is non-trivial then $\sigma$ is $\chi$-distinguished for some finite character $\chi$.
\end{prop}
\begin{proof} 
This proposition is part of \cite[Theorem 6]{BlasiusRogawski} and we recall the proof for the reader's convenience.

 Recall that  $\sigma$ is an irreducible cuspidal automorphic representation of 
$G(\A)$ which occurs as an irreducible subspace $V_{\sigma} \subset L^2 (G(F) \backslash G(\A))$.  
We identify $V_{\sigma}$ as a subspace 
$$V_{\sigma} \subset L^2 (\widetilde{G} (F) \backslash (\SO (p,q) \times \widetilde{G}(\A_f))).$$

We first relate the period integral to the pairing of a cohomology class with a special cycle, as done in \cite[Theorem 6]{BlasiusRogawski}.

\subsection{} Using the above identification we have:
\begin{equation} \label{eq:13.5}
H^{nq}_{\rm cusp} (\mathrm{Sh} (G) , E ) (\sigma) = \mathrm{Hom}_{K_{\infty}} (\wedge^{nq} \p , V_{\sigma} \otimes E).
\end{equation}
Recall from \S \ref{VZ2} that the $K_{\infty}$-module $V(n)$  occurs with multiplicity one in 
$\wedge^{nq} \p$ and that any non-zero element in the right-hand side of \eqref{eq:13.5} factorizes
through the isotypical component $V(n)$. We furthermore note that $V(n)$ occurs with multiplicity 
one in $V(\pi)\otimes E$. Since $V(\pi) = V(\sigma_{v_0})$ this leads to a canonical (up to multiples)
non-zero element  $\omega \in \mathrm{Hom}_{K_{\infty}} (V(n) , V(\sigma_{v_0}) \otimes E)$.

Let $\{v_j \}$ and $\{l_j \}$ be dual bases of $V(n)$ and $V(n)^*$, respectively. Fix a neat level 
$K \subset G( \A_f )$ such that $\sigma$ has $K$-invariant vectors. For $x_f$ in the space $V(\sigma_f)$
of $\sigma_f$, the element 
$$\omega \otimes x_f \in \mathrm{Hom}_{K_{\infty}} (V(n) , V(\sigma_{v_0})\otimes E ) \otimes V(\sigma_f)^K$$
corresponds to a $E$-valued harmonic form $\Omega_{x_f}$ on 
$$X_K =\widetilde{G} (F) \backslash (\SO (p,q) \times \widetilde{G}(\A_f)) / K_{\infty} \widetilde{K}$$
which decomposes as:
$$\Omega_{x_f} = \sum_j f_{j, x_f} \otimes e_{j} l_j$$
where $f_{j, x_f} \in V_{\sigma}  \subset L^2 (\widetilde{G} (F) \backslash (\SO (p,q) \times \widetilde{G}(\A_f)))$ is a cusp form, $e_{j} \in E$ and the tensor product $f_{j, x_f} \otimes e_{j}$ is the image of the vector $\omega (v_j) \otimes x_f$ under the map $\tau$.
The form $\Omega_{x_f}$ depends only on $x_f$ and $\omega$ but not on the choice of bases 
 $\{v_j \}$ and $\{l_j \}$. 
 
The Hodge $*$-operator establishes a one-to-one $K_{\infty}$-equivariant correspondence between 
$\wedge^{nq} \p$ and $\wedge^{(p-n)q} \p$ as well as on their dual spaces. It maps $\Omega_{x_f}$
onto the $E^*$-valued $(p-n)q$-harmonic form $*\Omega_{x_f} = \sum_j f_{j, x_f} \otimes e_{j} *l_j$.

Recall from \S \ref{7.4} that we have fixed a $\SO_0 (p,q)$-invariant inner product $(,)_{E}$ on $E$ and that this induces a natural inner product $\langle , \rangle$ on $E$-valued differential forms.
 
\subsection{} \label{Q} Let $\omega_H$ be a fixed non-zero element in the dual space of $\wedge^{(p-n)q} \p_H$ where 
$\p_H = \p \cap \mathfrak{h}$ and $\mathfrak{h}$ is the complexified Lie algebra of $\SO(p-n,q)$. Since
$\dim \wedge^{(p-n)q} \p_H =1$, $\omega_H$ is unique up to multiples. Each $*l_j$ may be represented as an element in the dual space of $\wedge^{(p-n)q} \p$. The restriction of 
$*\Omega_{x_f}$ to a connected cycle $c(U,g,K)_T$ is therefore of the form $(\sum_j c_j f_{j,x_f} \otimes e_{j} ) \omega_H$
where the $c_j$ are complex constants. Let 
$$f_{x_f,T} = \sum_j c_j (e_{j} , e_T)_{E^*} f_{j,x_f}.$$
Then for a suitable normalization of Haar measure $dh$, we have:
\begin{equation} \label{eq:period}
P(f_{x_f , T} , U , \chi ) = \langle [Z_{U,T}^{\chi} ] , \Omega_{x_f} \rangle
\end{equation}
where $\langle, \rangle$ denotes the inner product on differential forms. This remains true even if $G$
is not anisotropic: the projection of $[Z_{U,T}^{\chi}]$ in the cuspidal part of the cohomology belongs
to the $L^2$-cohomology.

If $Z^{nq} (\mathrm{Sh}(G) , E)(\sigma ) \neq \{ 0 \}$ equation \eqref{eq:period} therefore implies 
that 
$$P(f_{x_f , T} , U , \chi )  \neq 0$$ 
and $\sigma$ is $\chi$-distinguished. 
\end{proof}

\begin{ques}
Does the converse to Proposition \ref{prop:13.5} also hold ?
\end{ques}

Answering this question seems to lie beyond the tools of this paper. As a corollary of Proposition \ref{prop:13.5}  we nevertheless get the following:

\begin{thm} \label{thm:period}
Let $\sigma$ be an irreducible cuspidal automorphic representation of $G(\A)$ which occurs as an irreducible subspace in $L^2 (G(F) \backslash G(\A))$.
Suppose that the restriction of $\sigma_{v_0}$ to $\SO_0 (p,q)$ is a cohomological representation 
$\pi=A_{\mathfrak{q}} (\lambda)$ whose associated Levi subgroup $L$ is isomorphic to $\SO(p-2r, q) \times \mathrm{U} (1)^r$ with $p>2r$ and $m-1 >3r$ and such that $\lambda$ has at most $r$ nonzero entries. And suppose that for all infinite places $v \neq v_0$, the representation $\sigma_v$ is trivial. Then: $\sigma$ is $\chi$-distinguished for some finite character $\chi$.
\end{thm}
\begin{proof} Applying Theorem \ref{Thm:main9} --- with $r=n$ --- to each components of $X_K$ we conclude that the 
projections of the classes $[Z_{U,T}^{\chi}]$ for varying $U$, $T$ and $\chi$ generate $H^{qn}_{\rm cusp} (X_K , S_{[\lambda]} (V))_{\pi}$. (Here we should note that classes obtained by wedging a cycle
class with a power of the Euler form are not primitive and therefore project trivially into $H^{qn}_{\rm cusp} (X_K , S_{[\lambda]} (V))_{\pi}$.) Therefore: if $H^{nq}_{\rm cusp} (X_K , S_{[\lambda]} (V)) (\sigma)$
is non trivial, then $Z^{nq} (X_K , S_{[\lambda]} (V)) (\sigma) \neq \{ 0 \}$. And Theorem \ref{thm:period}
follows from Proposition \ref{prop:13.5}.
\end{proof}

\medskip
\noindent
{\it Remark.} We believe that the conditions $p>2r$ and $m-1 >3r$ are necessary in general but even assuming the results of \S \ref{6.13} we would still have to answer positively to question \ref{Q}.  

\medskip

Note that \cite[Theorem 1.1]{GJS} and the proof of Theorem \ref{thm:period} imply the following:

\begin{thm} \label{Conj:periods}
Let $\sigma$ be an irreducible cuspidal automorphic representation of $G(\A)$ which occurs as an irreducible subspace in $L^2 (G(F) \backslash G(\A))$.
Suppose that the restriction of $\sigma_{v_0}$ to $\SO_0 (p,q)$ is a cohomological representation 
$\pi=A_{\mathfrak{q}} (\lambda)$ whose associated Levi subgroup $L$ is isomorphic to $\SO(p-2r, q) \times \mathrm{U} (1)^r$ with $p\geq2r$, $q\geq 1$ and such that $\lambda$ has at most $r$ nonzero entries. And suppose that for all infinite places $v \neq v_0$, the representation $\sigma_v$ is trivial. 
Assume moreover that there exists a quadratic character $\eta$ of $F^* \backslash \A^*$ such that
the partial $L$-function $L^S (s, \sigma \times \eta )$ has a pole at $s_0=\frac{m}{2} - r$ and is holomorphic for $\mathrm{Re} (s) >s_0$.
Then: $\sigma$ is $\chi$-distinguished for some finite character $\chi$.
\end{thm}

We don't believe that belonging in the image of the theta correspondence --- which is a subtle properties away from stable range --- can be characterized
only by the existence of a pole for a partial $L$-function. We propose the following:

\begin{conj} \label{Conj2:periods}
Let $\sigma$ be an irreducible cuspidal automorphic representation of $G(\A)$ which occurs as an irreducible subspace in $L^2 (G(F) \backslash G(\A))$.
Suppose that the restriction of $\sigma_{v_0}$ to $\SO_0 (p,q)$ is a cohomological representation 
$\pi=A_{\mathfrak{q}} (\lambda)$ whose associated Levi subgroup $L$ is isomorphic to $\SO(p-2r, q) \times \mathrm{U} (1)^r$ with $p\geq2r$, $q\geq 1$ and such that $\lambda$ has at most $r$ nonzero entries. And suppose that for all infinite places $v \neq v_0$, the representation $\sigma_v$ is trivial. Then there exists an automorphic character $\chi$ such that $\sigma \otimes \chi$ is in the image of the 
cuspidal $\psi$-theta correspondence from a smaller group associated to a symplectic space of dimension $2r$ if and only if $\sigma$ is $\eta$-distinguished for some finite quadratic character $\eta$ 
and the global (Arthur) $L$-function $L (s, \sigma^{\rm GL} \times \eta)$ has a pole at $s_0=\frac{m}{2} - r$ and is holomorphic for $\mathrm{Re} (s) >s_0$.
\end{conj}

\subsection{A $3$-dimensional example} 
Let $F=\Q$ and let $V$ be the $4$-dimensional $\Q$-vector space with basis $e_1,e_2,e_3,e_4$ and quadratic form
$$q(x_1e_1 + x_2e_2+ x_3e_3+x_4e_4) = x_1^2 - x_2^2 - x_3^2 - d x_4^2.$$
Assume that $d$ is positive and not a square in $\Q$ and let 
$$D  = \left\{
\begin{array}{ll}
4d & \mbox{ if } d \equiv 1 , 2 \ (\mathrm{mod} \ 4) ,\\
d & \mbox{ if } d \equiv 3 \  \ (\mathrm{mod} \ 4) .
\end{array}
\right.$$ 
Consider the group $G= \mathrm{Res}_{\Q(\sqrt{-d}) / \Q} \GL(2)$. Recall --- e.g. from \cite[Proposition 3.14]{EGM} --- that $G$ is isomorphic to 
$\mathrm{GSpin} (V)$ over $\Q$. The corresponding arithmetic manifolds are therefore the same; these are the $3$-dimensional hyperbolic manifolds obtained as quotient of the hyperbolic $3$-space by Bianchi groups.

Distinguished representations of $G$ have already attracted a lot of attention: Consider periods associated to totally positive vectors in $V$. The corresponding groups $H$ are inner forms of 
$\GL(2)_{|\Q}$ and the corresponding cycles are totally geodesic surfaces. 
 
\begin{prop} \label{prop:bc}
Suppose that $\pi$ is an irreducible cuspidal automorphic representation of $G(\A)$. Then if $\pi$ is $\chi$-distinguished for some finite character $\chi$ then some twist of
$\pi$ is the base-change lift of a cuspidal automorphic representation of $\GL (2, \A)$.
\end{prop}
\begin{proof} This follows from the proof of \cite[\S 5]{Flicker}.
\end{proof}

Finis, Grunewald and Tirao \cite{FGT} compare two classical ways of constructing cuspidal
cohomology classes in $H^1 (\SL(2 , \mathcal{O}_{-d})$ where $\mathcal{O}_{-d}$ is the ring of integers 
of $\Q(\sqrt{-d})$. The first is the base-change construction where the corresponding cuspidal automorphic representations are obtained as twists of base-change lifts of cuspidal automorphic 
representation of $\GL (2 , \A)$. The second construction is via automorphic induction from Hecke
characters of quadratic extensions of $\Q(\sqrt{-d})$. In many cases, the part of the cuspidal cohomology 
thus obtained is already contained in the part obtained from the base-change construction. But it is not always the case --- see \cite[Corollary 4.16]{FGT} for a precise criterion --- and they in particular prove the 
following:

\begin{prop} 
Suppose there exists a real quadratic field $L$ such that 
$$\Q(\sqrt{-d}) L / \Q(\sqrt{-d})$$ is unramified and the narrow class number $h_L^+$ is strictly bigger than the corresponding of genera $g_L^+ = 2^{\mathcal{R}(L)|-1}$, where $\mathcal{R}(L)$ denotes the set of primes ramified in $L$. Then: there exists a non base-change cohomology class in $H^1_{\rm cusp} (\SL(2 , \mathcal{O}_{-d})$.
\end{prop}

The smallest discriminant of a real quadratic field $L$ such that $h_L^+ > g_L^+$ is $d_L=136$. 
Note that $\Q(\sqrt{-d}) L / \Q(\sqrt{-d})$ is unramified if and only if $d_L$ divides $D$ and $d_L$ and 
$D/d_L$ are coprime.

\medskip

As a corollary (using propositions \ref{prop:13.5} and \ref{prop:bc}) we get the following:

\begin{prop} \label{prop:dim3}
There exists a cohomology class in $H^1_{\rm cusp} (\SL(2 , \mathcal{O}_{-34}))$ which is not a linear
combination of classes of totally geodesic cycles.
\end{prop}

\newpage

\part*{Appendix}

In this appendix we prove Proposition \ref{Prop:appendix}. The strategy --- as suggested by Lemma \ref{Lem:S} --- relies on the general principle
stated by Clozel according to which {\it in an Arthur packet we find the representations belonging to the
Langlands packet and more tempered representations.} We address this question through 
the study of exponents.

\addtocounter{section}{1}

\setcounter{subsection}{0}

\subsection{} Here again we let $p$ and $q$ two non-negative integers with $p+q =m$ and let  $G= \SO (p,q)$. We set $\ell = [m/2]$ and $N=2\ell$. 

The goal of this appendix is to prove the following.

\begin{prop} \label{Prop:appendix2}
Let $\pi$ be an irreducible unitary cohomological representation of $G$ associated to a Levi subgroup $L=\SO(p-2r , q) \times \mathrm{U}(1)^r$. Assume that $\pi$ is the local component of an automorphic representation with associated (global) Arthur parameter $\Psi$. Assume that $3(m-2r-1)> m-1$. Then the parameter $\Psi$ contains a factor $\eta \boxtimes R_a$ where $\eta$ is a quadratic character and $a \geq m-2r-1$.
\end{prop}  

\medskip
\noindent
{\it Remark.} The inequality $3(m-2r-1)> m-1$ is equivalent to $m-1 >3r$. Proposition \ref{Prop:appendix} therefore indeed follows from Proposition \ref{Prop:appendix2}.

\medskip

Before entering into the proof we first recall some basic facts about the local classification of representations of $\GL (N ,\R)$ and $G$.

\subsection{Discrete series of the linear groups\label{seriediscrete}}
Let $k$ be a positive integer. If $k \geq 2$ we let $\delta(k)$ be the tempered irreducible representation of $\GL(2 , \R)$ obtained
as the unique irreducible quotient of the induced representation 
$$\mathrm{ind} ( | \cdot |^{\frac{k}{2}} \otimes | \cdot |^{-\frac{k}{2}})$$
(normalized induction from the Borel). It is the unique irreducible representation with trivial central character which
restricted to the subgroup $\SL^{\pm} (2 , \R)$ of elements $g$ such that
$|\det (g) | = 1$ is isomorphic to $\mathrm{ind}_{\SL(2 , \R)}^{\SL (2 , \R)^{\pm}} (D_k )$ where $D_k $ is the (more
standard) discrete series representation of $\SL (2 , \R)$ as considered in e.g. \cite[Chapter II, \S 5]{Knapp}. Recall that by definition the parameter of $\delta (k)$ 
is the half integer $(k-1)/2$

If $k=1$ we denote by $\delta (k)$ the trivial character of $\R^* = \GL (1 , \R)$. 

We note that if $\mu$ is a tempered irreducible representation of $\GL (d , \R)$ that is square integrable modulo the center then 
$d=1$ or $2$ and $\mu$ is obtained by tensoring some $\delta (k)$ with a unitary character $\nu$ of $\R^*$. More precisely:
\begin{itemize}
\item If $d=1$ either $\nu = 1 \otimes | \cdot |^{it}$ or $\nu = \mathrm{sgn} \otimes | \cdot |^{i t}$. Here 
$1$ denotes the trivial representation and $\mathrm{sgn}$ the sign character of $\R^*$ and $t \in \R$. Then $\mu = \nu$.
\item If $d=2$, $\nu = | \det (\cdot )|^{it}$  ($t \in  \R$). Then $\mu = \delta (k) \otimes | \det (\cdot )|^{it}$.
\end{itemize}
Hereafter we denote by $\mu (k , \nu)$  the representation obtained by tensoring $\delta (k)$ with $\nu$. And we simply denote by
$\mu(k , \nu) | \cdot |^{s}$ ($k \geq 1$, $s \in \C$) 
the representation $\mu (k , \nu) \otimes | \det (\cdot )|^{s}$. 

\subsection{Admissible representations of $\GL(N , \R)$}
Let $r$ be a positive integer and, for each $i= 1 , \ldots , r$, fix $k_i$ a positive integer and $\nu_i$ a unitary character of $\R^*$. 
We let $d_i= 1$ if $k_i = 1$ and let $d_i = 2$ if $k_i \geq 2$.
We assume that $N = d_1 + \ldots + d_r $.

Now let 
$\mathbf{x} = (x_1 , \ldots , x_r ) \in \R^r$ be such that $x_1 \geq \ldots \geq x_r$. Consider the induced representation of $\GL (N , \R)$:
$$\mathrm{I} ((k_i , \nu_i,  x_i)_{i=1 , \ldots , r}) =  \mathrm{ind} (\mu (k_1 , \nu_1 ) | \cdot |^{x_1} , \ldots , \mu (k_r , \nu_r ) | \cdot |^{x_r} )$$
(normalized induction from the standard parabolic of type $(d_1 , \ldots , d_r)$). 
We call such an induced representation a {\it standard module}. These generate\footnote{To have a basis we still have to take care of possible permutations of the indices $\{ 1 , \ldots , r \}$.} the Grothendieck group of the smooth admissible representation of $\GL (N, \R)$.

According to Langlands \cite{Langlands} $\mathrm{I} ((k_i , \nu_i , x_i)_{i=1 , \ldots , r})$ has a unique 
irreducible quotient. We note (see e.g. \cite[Chap. 3]{BC}) that --- when restricted to $\C^*$ --- 
the $L$-parameter of this representation is conjugate into the diagonal torus of $\GL (N, \C)$ and
each $(k_j , \nu_j , x_j)$ contributes in the following way:
\begin{itemize}
\item if $k_j = 1$ and $\nu_j = 1 \otimes | \cdot |^{it_j}$ or $\nu_j = \mathrm{sgn} \otimes | \cdot |^{i t_j}$, 
it contributes by $z \mapsto (z \overline{z})^{\frac{x_j +i t_j}{2}}$, and
\item if $k_j \geq 2$ and $\nu_j = | \det (\cdot )|^{it_j}$, it contributes by 
$$z \mapsto \left( \begin{array}{cc}
z^{x_j +it_j+ \frac{k_j}{2}} \overline{z}^{x_j+it_j- \frac{k_j}{2}} & \\
& z^{x_j+it_j- \frac{k_j}{2}} \overline{z}^{x_j +it_j + \frac{k_j}{2}} 
\end{array} \right).$$
\end{itemize}

\subsection{Arthur parameters}
Consider the outer automorphism:
$$\theta : x \mapsto J {}^t \- x^{-1} J  \ \  \ (x \in \GL(N , \R)).$$

Recall that a local Arthur parameter $\Psi$ is a formal sum of formal tensor products $\mu_j \boxtimes R_j$ 
where each $\mu_j$ is a (unitary) representation of say $\GL (a_j , \R)$ that is square integrable modulo the center, $R_j$ is an irreducible representation of $\SL_2 (\C)$ of dimension $b_j$ and $N = \sum_j a_j b_j$. We furthermore request that $\Psi^{\theta} = \Psi$. 

Note that we shall only be interested in local Arthur parameters which have the same (regular) infinitesimal character than a finite dimensional representation. This implies in particular that each  $\mu_{j}$ is a discrete series or a quadratic character.

As is now standard, for $j$ as above, we denote by $\mathrm{Speh}(\mu_{j},b_j)$ the (unique) irreductible quotient of the standard module obtained by inducing the representation 
$$\mu_{j} |\cdot |^{\frac12 (b_j-1)} \otimes \mu_{j} |\cdot |^{\frac12 (b_j-3)}  \otimes \ldots \otimes \mu_{j} |\cdot |^{\frac12 (1-b_j)}$$ 
as in \S \ref{S35}. Recall that the representation of $\GL (N , \R)$ associated to $\Psi$ is the induced representation of  $\otimes_{j} \mathrm{Speh} (\mu_{j},b_{j})$; it is irreducible and unitary representation.

We note that each $\mu_j$ is isomorphic to some $\mu (k , \nu)$ ($k \geq 1$, $\nu$ unitary character of $\R^*$). 
By Langlands' classification $\Pi_{\Psi}$ can then be realized as the unique irreducible sub-quotient of some standard module $\mathrm{I} ((k_i , \nu_i , x_i )_{i= 1 , \ldots ,r})$
where $(k_i , \nu_i , x_i)_{i= 1 , \ldots , r}$ are obtained as follows: We write $b_1 \geq b_2 \geq \ldots$. Then 
\begin{multline*}
| \cdot |^{x_1} \mu (k_1 , \nu_1) = | \cdot |^{\frac{b_1-1}{2}} \mu_1 , | \cdot |^{x_2} \mu (k_2 , \nu_2) = | \cdot |^{\frac{b_1 -1}{2}} \mu_2 \ (\mbox{if } b_2 = b_1) \\
\ldots  \  | \cdot |^{x_j} \mu (k_j , \nu_j ) = | \cdot  |^{\frac{b_1 -1}{2}} \mu_j \ (\mbox{if } b_j = b_1) .
\end{multline*}
We then put the characters of smaller absolute value, and so on. As $\Pi_{\Psi}$ is $\theta$-stable, we may furthermore arrange the $(k_i , \nu_i , x_i)$ 
so that there exists a integer $r_+ \in [0 , r/2]$ such that:
\begin{itemize}
\item For each $i= 1 , \ldots , r_+$, we have $k_i = k_{r-i+1}$, $\nu_i = \nu_{r-i+1}^{-1}$ and $x_i=-x_{r-i+1}$.
\item For any $j = r_+ +1 , \ldots , r-r_+$, $\nu_j$ is a quadratic character (trivial if $k_j >1$ with our convention) and $x_j=0$.
\item For any $i,j \in \{r_+ +1 , \ldots , r-r_+ \}$ with $i\neq j$, we have either $k_i \neq k_j$ or $\nu_i \neq \nu_j$. 
\end{itemize}
A standard module satisfying the above conditions is called a {\it $\theta$-stable standard module}.\footnote{Note that standard modules are generally not irreducible $\theta$-stability has thus to be defined.} 

To ease notations we will denote by $I(\lambda)$ a general standard module of $\GL(N ,{\mathbb R})$; then $\lambda$ has to be understood as $(k_i , \nu_i , x_i)_{i= 1 , \ldots , r}$.

\subsection{Twisted traces}
Let $\Pi$ be a $\theta$-stable irreducible admissible representation of $\GL (N , \R)$. We fix an action of $\theta$ on the space of $\Pi$ that is: an operator 
$A_{\theta}$ ($A_{\theta}^2 =1$) intertwining $\Pi$ and $\Pi^{\theta}$. For any 
test function $f \in C_c^{\infty} (\GL (N , \R))$ we can then form the $\theta$-trace 
$\mathrm{trace}_{\theta} \Pi (f) = \mathrm{trace} (\Pi (f) A_{\theta})$. As $A_{\theta}^2=1$ the action
of $\theta$ is well defined up to a sign. The $\theta$-trace of $\Pi$ is thus well defined --- independently of the choice of $A_{\theta}$ --- but only 
up to a sign. We shall fix the sign following Arthur's normalization (the so called Whittaker's normalization): Fix a character of ${\mathbb R}$. It defines a character of the upper diagonal unipotent subgroup of  $\GL(N,{\mathbb R})$. Recall that a Whittaker model is a map (continuous in a certain sense) between a smooth representation and the continuous function on $\GL(N,{\mathbb R})$ which transforms on the left under the unipotent subgroup through the character we have just defined. Any standard module has a unique Whittaker model (this is due to Shalika  and Hashizume). And any standard module whose parameter is invariant under $\theta$ has a unique action of $\theta$ which stabilizes the Whittaker model and varies holomorphicaly in the parameter, see \cite[\S 2.2]{Arthur} for more details. When $\Pi$ is $\theta$ invariant, its standard module has a Whittaker model and an action of $\theta$. This action restricts to an action of $\theta$ on the space of $\Pi$ which gives our choice of normalization for the operator $A_\theta$. Using this normalization $\mathrm{trace}_{\theta} (\Pi )$ is now well defined. 

More generally if $I(\lambda)$ is a $\theta$-stable standard module, we shall denote by $\mathrm{trace}_{\theta}I(\lambda)$ its twisted trace normalized, as above, using the theory of Whittaker models.

\subsection{Stabilization of the trace formula} \label{stabilisation}
We now turn to the global study. Let $F$ be a number field and $v_{0}$ a real place of $F$. We denote by $G$ an orthogonal group defined over $F$ such that $G(F_{v_0}) = 
\SO (p,q)$, our local $G$ above. We also fix $\nu$ an infinitesimal character of $\SO (p,q)$ which is the character of a finite dimensional representation of $\SO (p,q)$.\footnote{This will be the infinitesimal character of the cohomological representation $\pi$ of the proposition.} We denote by $R_{{\rm disc},\nu}^{G}$ the subspace of the square integrable funtions on $G(F)\backslash G({\mathbb A})$ with infinitesimal character $\nu$ at the place $v_{0}$. 
Let $V$ be a finite set of places of $F$ big enough so that it contains all the Archimedean places and for $v\notin V$ the group $G(F_v)$ is quasi-split and splits over a finite unramified extension of $F_v$; in particular $G(F_v)$ contains a hyperspecial compact subgroup $K_v$, see \cite[1.10.2]{Tits}. As usual we shall denote by $F^V$, resp. $F_V$, the restricted product of all completions $F_v$ for $v \notin V$, resp. $v \in V$. Similarly we shall use the corresponding standard notations $G(F_V)$, $G(F^V)$ and $K^V= \prod_{v \notin V} K_v$.

Denote by $R_{{\rm disc},\nu}^{G,V}$ the subspace of $R_{{\rm disc},\nu}^G$ which consists of the representations of $G$ that are unramified outside $V$. Let $c^V$ be a character of the spherical Hecke algebra of $G(F^V)$ and denote by $R_{{\rm disc} \nu,  c^V }^G$ the subspace of $R_{{\rm disc},\nu}^{G,V}$ where this spherical algebra acts through the character $c^V$.

It is a consequence of Arthur's work that $R_{{\rm disc},\nu,c^V}^G$ is a representation of finite length. Given a test function $f_V$ on $G(F_V)$ the product $f_V 1_{K^V}$ defines a test function on $G(\mathbb{A})$. We shall denote by $I_{{\rm disc},\nu,c^V}^G$ the distribution
$$f_V \mapsto \mathrm{trace} \ R_{{\rm disc},\nu,c^V}^G ( f_V 1_{K^V} )$$
on $G(F_V)$.

Denote by $G^*$ the quasi-split inner form of $G$. Arthur defines a distribution $S_{{\rm disc},\nu,c^V}^{G^*}$ (supported on character) of $G^* (F_{V})$ inductively by the following formula: for any $f_{V}$ a test function on $G^*(F_{V})$
$$
S_{{\rm disc},\nu,c^V}^{G^*} (f_{V}) := I_{{\rm disc},\nu,c^V}^{G^*} (f_{V}1_{K^V}) - \sum_{H} \iota (G^* ,H) S_{{\rm disc},\nu,c^V}^{H} (f^H_{V}),
$$
where $H$ runs through the set of all elliptic endoscopic data of $G^*$ different of $G^*$ itself and the coefficients $\iota (G^*,H)$ are certain positive rational numbers. We shall not give details, the only thing that matters to us is that such an $H$ is a product of two non trivial special orthogonal group and is quasi-split. Then the function $f^H_{V}$ is the Langlands-Shelstad transfer of $f_{V}$. We also have to explain the meaning of $\nu$ and $c^V$ for $H$: these are respectively the sums over the $\nu'$'s and the $(c')^{V}$'s, resp. an infinitesimal character of $H(F_{v_{0}})$ and a character of the spherical algebra of $H(F^V)$ which transfer to $\nu$ and $c^V$ through the Langlands functoriality.

We can now state the {\it stabilization of the trace formula}:
\begin{enumerate}
\item the distribution $S_{{\rm disc},\nu,c^V}^{G^*}$ is stable;
\item for any test function $f_{V}$ on $G(F_{V})$, we have
\begin{equation} \label{eq:app1}
I_{{\rm disc},\nu,c^V}^G (f_V) = \sum_{H} \iota (G,H) S_{{\rm disc},\nu,c^V}^{H} (f^H_{V}),
\end{equation}
where now $H$ runs through the set of {\it all} (including $G^*$) elliptic endoscopic data of $G$ and the coefficients $\iota (G,H)$ are again positive rational numbers.\footnote{Note that \eqref{eq:app1} --- applied to $G^*$ rather than to $H$, and inductively to its endoscopic subgroup --- uniquely defines the distributions $S_{{\rm disc},\nu,c^V}^{H}$.}
\end{enumerate} 
 
\subsection{Stabilization of the twisted trace formula} \label{stabilisationtordue}
Keep notations as above. Arthur has given a way to compute the distribution $S_{{\rm disc},\nu,c^V}^{G^*}$. First recall that such a distribution only depends on the image of $f_{V}$ modulo the function whose stable orbital integral are zero. Thanks to the good property of the twisted transfer, such a $f_{V}$ has the same image in the quotient that the twisted transfer to $G(F_{V})$ of a function $f^{\GL}_{V}$ and Arthur (see \cite[\S 3.4]{Arthur}) has proved: there exists a (global) Arthur parameter $\Psi$ such that
\begin{equation} \label{eq:app2}
\mathrm{trace}_{\theta}(\Pi_\Psi )(f^{\GL}_{V}1_{\tilde{K}^V})=S_{{\rm disc},\nu,c^V}^{G^*} (f_{V}) ; 
\end{equation}
here $\tilde{K}^V$ is the product of the hyperspecial compact subgroups $\tilde{K}_v$ of $\GL (N, F_v)$. In particular: outside $V$, the representation $\Pi_\Psi$ is unramified and the character of the twisted spherical algebra is obtained from $c^V$ by functoriality.

\subsection{Two invariants and the statement}
Let $\pi_{0}$ be a cohomological representation of $G(F_{v_0}) \cong \SO (p,q)$ associated to the Levi subgroup
$$L=\mathrm{U}(p_1 ,q_1) \times \ldots \times \mathrm{U}(p_r , q_r) \times \mathrm{SO}(p_0 , q_0 )$$
with $p_0+2\sum_{j}p_{j} = p$ and $q_0+2\sum_{j} q_{j}=q$. Set 
$$m_{0} (= m_0 (\pi_{0})) = p_0 +q_0.$$ 
Note that this only depends on $\pi_0$ and not on a particular choice of $L$. 
 
We now define $m^{\GL} (\pi_{0})$ in the following way. Let $\pi$ be a square integrable representation of $G({\mathbb A})$ which localize in $\pi_{0}$ at the place $v_{0}$. The global representation $\pi$ determines two characters $\nu$ and $c^V$ as in \S \ref{stabilisation}. It then follows from \S \ref{stabilisationtordue} that these two characters determine some global Arthur parameter $\Psi$. Next we localize $\Psi$ at the place $v_{0}$; in the subspace where $W_{{\mathbb R}}$ acts through characters in this localization, we look at the action of $\SL(2,{\mathbb C})$ and denote by $m (\Psi )$ the biggest dimension for an irreducible representation of $\SL(2,{\mathbb C})$ in this space. Finally we set:
$$
m^{\GL}(\pi_{0}):=\min_{\pi} m (\Psi ).
$$

\begin{prop} \label{P:app}
We have the inequality: $m_{0}(\pi_{0})-1 \leq m^{\GL}(\pi_{0})$.
\end{prop}

In fact the inequality is an equality but we do not need it here. 

\subsection{Proposition \ref{P:app} implies Proposition \ref{Prop:appendix2}}
Indeed take $\pi_0 = \pi$ as in Proposition \ref{Prop:appendix2} so that $m_0 = p+q-2r=m-2r$. Note that the representation $\pi_0$ being cohomological its infinitesimal character is regular. So if $\pi_0$ is the local component of an automorphic representation $\pi$ with associated Arthur parameter $\Psi$, then $\Psi$ is very particular: writing 
$$\Psi = \mu_1 \boxtimes R_1 \boxplus \ldots \boxplus \mu_r \boxtimes R_r$$
the regularity of the infinitesimal character has the following consequences.
\begin{enumerate}
\item If $m$ is odd. Then there is at most one $j$ such that $\mu_{j}$ is a quadratic character. If there exists such a $j$ we enumerate so that $j=1$. 
\item If $m$ is even there is no such $j$ or there is two such $j$. In the later case, we enumerate so that these two $j$'s are $1$ and $2$ and $R_{1}$ is maximal. We then have $R_{2}=1$. 
\end{enumerate}
In particular, $m (\Psi ) =n_1$ and we conclude from Proposition \ref{P:app} that 
$$n_1 \geq m_0-1 = m-2r-1.$$
Finally if $m$ is odd, $3n_1  \geq 3 (m-2r-1) >m-1=N$, and if $m$ is even $3n_1 +1 \geq 3 (m-2r-1)+1 >m=N$, it therefore follows from 1. and 2. above that the local character $\mu_1$ can only occur as the localization of a global character. This concludes the proof of Proposition \ref{P:app}.

\medskip

We shall prove Proposition \ref{P:app} by decomposing the traces and twisted traces as sums of traces and twisted traces of standard modules. We shall therefore decompose cohomological representations in the Grothendieck group using standard modules.

\subsection{Standard modules and their exponents}
The only representations of $\SO (p,q)$ we will be interested in have real infinitesimal character. The relevant standard modules $I_P (\sigma)$ of $\SO(p,q)$ are induced from a standard parabolic subgroup $P=MN$ of $\SO(p,q)$ with
$$M \cong \GL(d_1, \mathbb{R}) \times \cdots \times \GL(d_t , \mathbb{R}) \times \SO(p-n_0 , q-n_0),$$
where $t \in \mathbb{N}$, $d_i = 1$ or $2$, and $n_{0}= d_1 + \ldots +d_t$,
and
$$\sigma = \delta_1 | \cdot |^{x_1} \otimes \ldots \otimes \delta_t |\cdot |^{x_t} \otimes \pi_0.$$
Here each $\delta_{i}$ is either a discrete series of $\mathrm{GL}(2,{\mathbb R})$, if $d_i =2$, or a quadratic character of ${\mathbb R}^*$, if $d_i=1$, the $x_{i}$'s are positive real numbers with $x_{1}\geq x_{2} \geq \cdots \geq x_{t}$, and $\pi_{0}$ is a tempered representation of the group $\SO(p-n_{0},q-n_{0})$. 

The set of {\it exponents} of such a standard module is the set $\{\pm x_1 , \ldots , \pm x_t \}$. We define the set of {\it character exponents} as
$$\mathrm{CarExp} (I_P (\sigma)) = \{\pm x_{i} \; : \;  i\in [1,t] \mbox{ and } d_{i}=1\}.$$


\begin{lem} \label{L:app1}
Let $\pi_{0}$ be a cohomological representation of $\SO (p,q)$. There exists a finite set $\mathcal{E}$ of data $(P, \sigma)$ such that the corresponding standard modules $I_P (\sigma )$ satisfy
$$\mathrm{CarExp} (I_P (\sigma)) \subset \left[-\frac{m_{0}-2}{2} , \frac{m_{0}-2}{2} \right]$$ 
and, in the Grothendieck group,
$$
\pi= \sum_{(P, \sigma) \in \mathcal{E}} m(P, \sigma ) I_P(\sigma), \quad (m(P, \sigma)\in {\mathbb Z}).
$$
\end{lem}
\begin{proof} Since $\pi_0$ is cohomological, it is obtained by cohomological induction of a character of the associated Levi subgroup $L$. To decompose $\pi_0$ in the Grothendieck group using standard modules, we take the resolution of the character of $L$ with standard modules and cohomologically induce it; see Johnson's thesis \cite{Johnson}. The part coming from the unitary group gives only exponent of the form $\delta | \cdot |^x$ where $\delta$ is a discrete series of $\GL (2 , \mathbb{R})$. And the exponent coming form the $\SO (p_0, q_0)$-part are of the form $|\cdot |^{-x}$ with $0 < x \leq \frac12 (m_0-2)$. 
\end{proof}

\subsection{The twisted case} \label{app:tc}
Let $m(\Psi)$ We look only at very particular local Arthur parameters $\Psi$ since the infinitesimal character of the corresponding representation is regular. Writing 
$$\Psi = \mu_1 \boxtimes R_1 \boxplus \ldots \boxplus \mu_r \boxtimes R_r$$
recall that:
\begin{itemize}
\item If $m$ is odd. Then there is at most one $j$ such that $\mu_{j}$ is a quadratic character. If there exists such a $j$ we numerate so that $j=1$. 
\item If $m$ is even there is no such $j$ or there is two such $j$. In the later case, we numerate so that these two $j$'s are $1$ and $2$ and $R_{1}$ is maximal. We then have $R_{2}=1$. 
\end{itemize}
Set $m(\Psi)=n_1$ (the dimension of $R_1$). 

Standard modules $I(\lambda)$ for $\mathrm{GL}(N,{\mathbb R})$, as for $\SO (p,q)$, are representations induced from a tempered representation modulo a character of a Levi subgroup. The tempered representation contains twist of discrete series and twist of quadratic character (since we only consider representation with real infinitesimal character). 
We similarly denote by $\mathrm{CarExp}(I(\lambda))$ the set of real number $\pm x_{i}$  such that $x_{i}>0$ and $x_{i}$ is the absolute value of a character occuring in the definition of  $I(\lambda)$:
$$\mathrm{CarExp}(I( (k_i , \nu_i , x_i)_{i= 1 , \ldots , r})) = \{ \pm x_{i} \; : \;  i\in [1,r], \ x_i \neq 0 \mbox{ and } k_{i}=1\}.$$

Recall that we denote by $\Pi_\Psi$ the representation of $\mathrm{GL}(n,{\mathbb R})$ associated to $\Psi$. We shall proof the following analogue of Lemma \ref{L:app1}.

\begin{lem} \label{L:app2}
There exists a finite set $\mathcal{E}_\theta$ of data $\lambda$ such that the standard modules $I(\lambda)$ are $\theta$-stable and 
$$\mathrm{CarExp}I(\lambda) \subset [-(m(\Psi)-1)/2,(m(\Psi)-1)/2]$$ 
and 
$$
\mathrm{trace}_{\theta} \ \Pi_\Psi = \sum_{\lambda \in \mathcal{E}_\theta } m(\pi,\lambda)\mathrm{trace}_{\theta} \ I(\lambda), \quad (m(\pi,\lambda)\in {\mathbb Z}).
$$
\end{lem}

The proof requires a bit of preparation. We first write each module $\mathrm{Speh}(\mu,b)$ in the Grothendieck group using the basis made by standard module:
\begin{equation} \label{add1}
\mathrm{Speh}(\mu , b)=\sum_{\lambda }m(\lambda )I(\lambda).
\end{equation}
Let $k$ the positive integer such that $\mu=\delta (k)$. The condition that the infinitesimal character of $\mathrm{Speh} (\mu, b)$ is regular is equivalent to the fact that $k>b$.

\begin{lem} \label{L:add1}
We assume that $k\geq b$. In \eqref{add1}, if $m(\lambda ) \neq 0$ then $\mathrm{CarExp}I(\lambda) = \emptyset$.
\end{lem} 
\begin{proof} The module $\mathrm{Speh}(\mu,b)$ is cohomologicaly induced by a unitary character $(z/\overline{z})^{k/2}$ of $\GL(b,{\mathbb C})$. To obtain \eqref{add1} we take the analogous resolution of the trivial character of $\GL(b ,{\mathbb C})$ tensor it by the previous character and then make the cohomological induction. In our situation cohomological induction is an exact functor which sends irreducible representations to irreducible representations because the infinitesimal character is regular. Moreover this functor sends a standard module to a standard module. So it's enough to prove the analogous lemma for the trivial character of $\GL(b,{\mathbb C})$. But the resolution in this simple situation is perfectly known: the analogous of \eqref{add1} is (up to a sign)
\begin{equation} \label{add2}
\sum_{\sigma\in {\mathfrak S}_{b}}(-1)^{\ell(\sigma)}I(\lambda,\sigma(\lambda)),
\end{equation}
where $\lambda=((b-1)/2, \cdots, -(b-1)/2)$ is seen as a set of $b$ elements ($\lambda_{1},\cdots, \lambda_{b})$ and $I(\lambda,\sigma(\lambda))$ is the principal series of $\GL(b,{\mathbb C})$ induced by the character $\otimes_{i\in [1,b]}z^{\lambda_{i}}\overline{z}^{\sigma(\lambda_{i})}$ of the diagonal torus. The length $\ell(\sigma)$ is the ordinary length in ${\mathfrak S}_{b}$. The exponents of $I(\lambda,\sigma(\lambda))$ are obtained as products of a unitary character with an absolute value $(z\overline{z})^y$ with
$$
y \in \{ (\lambda_{i}+\sigma(\lambda_{i}))/2 \; : \; i\in [1,b] \} \quad (\mbox{and } \vert y\vert  \leq (b-1)/2.
$$
The unitary part of the character is of the form $(z/\overline{z})^c$ with $c= (\lambda_{i}-\sigma(\lambda_{i}))/2$. In particular $c> -k/2$. Recall that we have to cohomologically induce this character tensorised by the character $(z/\overline{z})^{k/2}$. This is necessarilly a non trivial unitary character and the induced representation is a discrete series.
\end{proof}

We shall need a more technical result which is a corollary of the proof:

\begin{lem} \label{L:add1bis}
Let $\rho$ be an irreducible subquotient of any standard module appearing in \eqref{add1} with non-zero coefficient. Then there exists some $\lambda$ such that $I(\lambda )$ appears in \eqref{add1} with $m(\lambda ) \neq 0$ and $\rho$ is the Langlands quotient of it.
\end{lem}
\begin{proof} We go back to the previous proof. Let $\rho$ be as in the statement of the lemma. By the previous proof, $\rho$ is cohomologicaly induced by an irreducible subquotient $\rho_{{\mathbb C}}$ appearing as subquotient of one of the principal series of \eqref{add2}. But in \eqref{add2} all possible principal series appear (the infinitesimal character is fixed, of course). So $\rho_{{\mathbb C}}$ is a Langlands quotient of such a principal series and $\rho$ is the Langlands quotient of the standard module obtained from it by cohomological induction.
\end{proof}

Now fix $b_{0}\in {\mathbb N}$ and denote by $\epsilon_{b_{0}}$ either the trivial representation of $\GL(b_{0},{\mathbb R})$ or the sign representation of this group. Denote by $\pi_{0}$ either $\epsilon_{b_0}$ or the representation of $\GL(b_{0}+1,{\mathbb R})$ obtained by inducing the tensor product of the representation $\epsilon_{b_{0}}$ of $\GL(b_{0},{\mathbb R})$ and the trivial representation of ${\mathbb R}^\times$ from the maximal parabolic corresponding to the partition $(b_{0},1)$ of $b_{0}+1$. In any case $\pi_0$ is a representation of $\GL(b_{0}+\eta,{\mathbb R})$ where $\eta=0$ or $1$. Write the analogue of \eqref{add1} for $\pi_{0}$:
\begin{equation} \label{add3}
\pi_{0}=\sum_{\lambda} m(\lambda) I(\lambda ).
\end{equation}

\begin{lem} \label{L:add2}
In \eqref{add3}, the exponents of $\lambda$ are of absolute value less than or equal to $(b_0-1)/2$ and this is also true for any irreducible subquotient $\rho$ of the $I(\lambda)$ appearing in \eqref{add3}. This last property is true for any irreducible representation of $\GL(b_{0}+\eta ,{\mathbb R})$ with the same infinitesimal character as $\pi_{0}$.
\end{lem}
\begin{proof} It is not obvious to write down explicitly \eqref{add3} but the second assertion of the lemma is more general than the first one and do not need to know \eqref{add3} to prove it. We have to prove that any representation of a Levi subgroup of $\GL(b_{0}+\eta,{\mathbb R})$, tempered modulo center, with the same infinitesimal character than $\pi_{0}$ satisfies the lemma. Such a representation can be written as a tensor product $\otimes \delta (k_j) | \cdot |^{x_j}$ with notations as above. Then the condition on the infinitesimal character implies that 
$$k_{j}-1+x_{j}\in [(b_0-1)/2,-(b_0-1)/2] \quad \mbox{ and } \quad -k_{j}+1+x_{j} \in [(b_0-1)/2,-b_0-1)/2].$$ 
In particular $k_{j}-1+\vert x_{j}\vert  \in [(b_0 -1)/2,-(b_0 -1)/2]$ and $\vert x_{j}\vert \leq (b_0-1)/2$. This proves the lemma.
\end{proof}

We now come back to the local Arthur parameter $\Psi$ and the corresponding representation $\Pi_{\Psi}$. Recall that $\Pi_\Psi$ is the induced representation of $\otimes_j \mathrm{Speh} (\mu_j , b_j)$. Since $\Psi$ satisfies (1) and (2) of \S \ref{app:tc} we may rewrite $\Pi_{\Psi}$ as the induced representation of 
$$\mathrm{Speh} (\delta (k_1 ) ,  b_1) \otimes \ldots \mathrm{Speh} (\delta (k_t ) ,  b_t) \otimes \pi_0$$
with $k_j\geq 2$, for $j=1, \ldots , t$, and either $\pi_0$ is as above or does not appear. In the latter case we will put $b_0 =0$. Moreover: the infinitesimal character of $\Pi_\Psi$ is ``almost'' regular, meaning that it is regular if $\pi_{0}$ is a character and that, if $\pi_{0}$ is not a character, the infinitesimal character of the representation induced from  
$\otimes_{j} \mathrm{Speh} (\delta (k_j ) , b_{j}) \otimes \epsilon_{b_{0}}$ is regular. Note that $N=2(b_1 + \ldots + b_t) + b_0$ or $N=2(b_1 + \ldots + b_t) + b_0+1$ according to the parity of $b_0$ (or equivalently the parity of $m$) and $b_0 = m(\Psi)$. 

Now we decompose each representation $ \mathrm{Speh} (\delta (k_j ) , b_{j})$ as in \eqref{add1}. And for any $j\in \{ 1 , \ldots , t \}$ we let $\rho_{j}$ be a subquotient of a standard modules appearing non-trivially, i.e. with non-zero coefficient $m(\lambda )$, in this decomposition. Finally we let $\rho_{0}$ be an irreducible representation with the same infinitesimal character as $\pi_{0}$. 

\begin{lem} 
The representation of $\GL (N, \mathbb{R})$ induced from $\rho_{1} \otimes \ldots \otimes \rho_t \otimes \rho_0$ is irreducible.
\end{lem}
\begin{proof} This lemma will appear in the thesis of N. Arancibia but, for the ease of the reader we briefly include a proof. It follows from the properties of the infinitesimal character and (now classical) results of Speh on the irreducibility of induced representations of $\GL$~: suppose that $\delta (k) | \cdot |^x$ appears in one the standard module for $\rho_{j}$  and $\delta' (k') | \cdot |^{x'}$ appears in a standard module for $\rho_{j'}$ with $j\neq j'$. Then we have 
$$
\frac{k-1}{2}+x \in \left[\frac{k_{j}+b_{j}-2}{2} ,\frac{|b_{j}-a_{j}|}{2} \right]
$$
and
$$
-\frac{k-1}{2}+x \in \left[-\frac{|b_{j}-a_{j}|}{2},- \frac{b_{j}+a_{j}-2}{2} \right];
$$
and similarly with $'$. But the two sets for $j$ are symmetric to $0$ and the corresponding set for $j'$ have the same property and are disjoint from the sets for $j$.  After Speh (see also \cite[Lemma 1.7]{MW}) this is a enough to conclude that the induced representation is irreducible. 
\end{proof}

Denote by ${\mathcal F}$ the set of irreducible representations $\sigma$ as in the previous lemma. This set has a length, Vogan's length. 
The tempered representations are of length $0$. Let $\sigma\in {\mathcal F}$ be a self-dual representation and denote by $I(\lambda_{\sigma})$ the standard module of $\sigma$. Whittaker normalization provides choices of actions on both $\sigma$ and $I(\lambda_\sigma )$ that are compatible; the twisted characters $\mathrm{trace}_{\theta}(\sigma)$ and $\mathrm{trace}_{\theta}I(\lambda_{\sigma})$ are taken with respect to these actions.
 
\begin{lem} \label{L:add3}
Let $\sigma$ be as above. Then there exists a finite set subset ${\mathcal F}_{\sigma}$ of $ {\mathcal F}$ containing only self-dual representations of length strictly less than the length of  $\sigma$ such that for suitable $m(\tau,\sigma) \in {\mathbb Z}-\{0\}$
$$
\mathrm{trace}_{\theta}(\sigma)-\mathrm{trace}_{\theta} I(\lambda_{\sigma})=\sum_{\tau}m(\tau,\sigma) \mathrm{trace}_{\theta}(\tau).
$$
\end{lem}
\begin{proof} If $\sigma$ is tempered we can take ${\mathcal F}_{\sigma}=\emptyset$. In general we prove the lemma by induction. We shall first prove that any irreducible subquotient of $I(\lambda_{\sigma})$ is in ${\mathcal F}$. By definition $\sigma$ is an induced representation of the $\rho_{j}$. So $I(\lambda_{\sigma})$ is the standard module induced from the standard modules of the $\rho_{j}$. Let $\tau$ be a subquotient of $I(\lambda_{\sigma})$. For any $j$ there exists an irreducible subquotient $\tau_{j}$ of the standard module of $\rho_{j}$ such that $\tau$ is a subquotient of the induced representation of the $\tau_{j}$. But we have seen that such an induced representation is irreducible and has to coincide with $\tau$. This proves that $\tau \in {\mathcal F}$. If the length of $\tau$ equals the length of $\sigma$ then $\tau=\sigma$ and otherwise the length of $\tau$ is strictly less than the length of $\sigma$. This proves that $\mathrm{trace}_{\theta}(\sigma)-\mathrm{trace}_{\theta}I(\lambda_{\sigma})$ is a linear combination of the  $\mathrm{trace}_{\theta}\tau$ for those $\tau$ which are self-dual (up to a sign wich depends of the choices).
\end{proof}

\subsection{Proof of Lemma \ref{L:app2}} 
We apply Lemma \ref{L:add3} to $\sigma = \Pi_\Psi$. The twisted trace of $\sigma$  can be written as a linear combination of twisted trace of standard modules $\mathrm{trace}_{\theta} (I_{\lambda_{\tau}})$ where $\tau$ runs in a subset of self-dual representations in ${\mathcal F}$ of length smaller or equal than that of $\sigma$. Moreover, it 
follows from the description of the exponent for the representation inducing the element in ${\mathcal F}_\sigma$ (see Lemma \ref{L:add1}, \ref{L:add1bis} and \ref{L:add2}) that if one of the representation $\tau$ has an exponent which is a character with absolute value $x$ then $b_{0}\geq 2\vert y\vert+1$. Since $b_0 = m (\Psi)$ we are done. \qed

\subsection{Transfer, local version}\label{transferlocal}
We now come back to the setting (and notations) of \S \ref{stabilisation} and \S \ref{stabilisationtordue}. The distribution $S_{{\rm disc}, \nu , c^V}^{G^*}$ on $G^* (F_{V})$ is a product of local stable distributions with a global coefficient: fix $v\in V$, there exists a finite set $\prod_v = \prod (\Psi_v)$ of representations $\pi_v$ of $G^* (F_v )$ and some multiplicities $m(\pi_v ) >0$ and signs $\varepsilon (\pi_v )$ ($\pi_v \in  \prod_v )$ such that such that 
$$
S_{{\rm disc}, \nu , c^V}^{G^*} (f_{V})= x(c^V) \prod_{v\in V} \left( \sum_{\pi_v \in \prod_v } \varepsilon (\pi_v) m(\pi_v ) \mathrm{trace} \  \pi_v (f_{v}) \right), 
$$

where $x(c^V)$ is a global constant. And the local packets $\prod_v$ (with the corresponding multiplicities and signs) are determined (see Proposition \ref{P37}) by:
\begin{equation} \label{eq:app3}
\sum_{\pi_v \in \prod_v } \varepsilon (\pi_v) m(\pi_v ) \mathrm{trace} \  \pi_v (f_{v}) = \mathrm{trace}_{\theta} \ \Pi_{\Psi_v} (\tilde{f}_{v}), 
\end{equation}
where $\tilde{f}_{v}$ is a function on $\GL(N,F_{v})$ whose twisted transfer to $G(F_{v})$ is $f_{v}$ modulo unstable functions. We will denote by $\mathrm{trace} ( \prod (\Psi_v))$ the (local stable) distribution on the left hand side. Similarly we denote by $\mathrm{trace} ( \prod (\Psi_v)^H)$ the local stable distribution on $H(F_v)$ associated to $S_{{\rm disc}, \nu , c^V}^{H} (f_{V})$.

Now we turn to the stabilization of the trace formula for $G$. The left hand side of \eqref{eq:app1} is (the character of) a linear combination with positive coefficients of representations of $G(F_{V})$. Because the right hand side is of finite length as a representation of endoscopic groups, we conclude that this also holds for the left hand side. 

Assume now that $\pi_{0}$ is a cohomological representation occuring as the local component of at least one of these representations. We first fix $f_{V}$ outside of $v_{0}$. For each choice, $f_{V}^{v_{0}}$, we have a distribution 
$$f_{v_{0}}\mapsto I_{{\rm disc} , \nu , c^V}^{G} (f_{v_{0}}f_{V}^{v_{0}}1_{K^V});$$ 
it is a finite linear combination of traces of representations of $G(F_{v_0} ) = \SO (p,q)$. Because of the positivity of the coefficients in the decomposition of $R_{{\rm disc} , \nu , c^V}^{G}$ and the linear independence of characters, we may find a test function $f_{V}^{v_{0}}$ such that this linear combination contains the trace of $\pi_{0}$. 

In this way we define a finite set ${\mathcal F}$ of representations of $G(F_{v_0} ) $ containing $\pi_{0}$ such that for any test function $f_{v_0}$ on $G(F_{v_0} ) $ we have:
\begin{equation} \label{eq:app4}
\sum_{\pi \in {\mathcal F}} c(\pi) \mathrm{trace} \ \pi(f_{v_{0}})= \sum_{H} x(H) \mathrm{trace} ( \prod (\Psi_{v_{0}})^{H}) (f^{H}_{v_{0}}), 
\end{equation}
for suitable coefficients $c(\pi)$ and $x(H)$. What is important for us is that $c(\pi_{0})\neq 0$. We recall that any $\pi \in {\mathcal F}$ is a unitary representation with infinitesimal character $\nu$, this implies (see \cite{Sala}) that any $\pi \in {\mathcal F}$ is cohomological.

\subsection{End of the proof of Proposition \ref{P:app}} \label{findelapreuve}
We first decompose the left hand side of \eqref{eq:app4} in terms of standard modules. It follows from Lemma \ref{L:app1} that there is certainly one standard module $I_P (\sigma)$ of $G(F_{v_0}) = \SO(p,q)$ that contributes non trivially to the left hand side of \eqref{eq:app4} and whose associated set $\mathrm{CarExp} (I_P (\sigma ))$ contains a term $(m-2)/2$ with $m\geq m_{0}(\pi_{0})$.  Now we decompose the right hand side of \eqref{eq:app4} in terms of ``stable'' standard modules (instead of inducing a tempered representation modulo the center, we induce a stable tempered representation modulo the center (see \cite{arthurselecta} and \cite[\S 3.3 and \S 3.4]{stabilisationIV}). At least one of these stable standard modules is such that its associated set $\mathrm{CarExp}$ contains a term $(m-2)/2$ with the same $m$ as above but with a suitable $H$. To ease the understanding we assume that $H=G^*$ (otherwise we have to introduce product of orthogonal groups and a decomposition of $\Psi$ into a product). We now use \eqref{eq:app3}: by the results of Mezo \cite{mezo}
a stable standard module for $G^*$ has a twisted transfer to $\GL(N,{\mathbb R})$ which is a $\theta$-stable standard module with the predicted Langlands parameter. In particular the left hand side of \eqref{eq:app3} contains a $\theta$-stable standard module with in $\mathrm{CarExp}$ a term $(m-2)/2$. It then follows from Lemma \ref{L:app2} that $m-1\leq m (\Psi)$ and Proposition \ref{P:app} follows.

\bibliography{bibli}

\bibliographystyle{plain}

\end{document}